\documentclass[a4paper, 11pt]{article}

\usepackage[T1]{fontenc}

\usepackage[plain]{fullpage}

\usepackage[dvipsnames,svgnames,table]{xcolor}

\usepackage[width=0.8\textwidth]{caption}

\usepackage{hyperref}
\hypersetup{
    pdftitle={Excluding a rectangular grid},
    pdfauthor={Clément Rambaud},
    colorlinks,
    linkcolor={RoyalBlue},
    citecolor={RubineRed},
    urlcolor={blue!80!black}
}

\usepackage{float}

\usepackage{amsmath,amssymb}

\usepackage[inline]{enumitem}
\setlist[enumerate,1]{label=(\roman*)}
\newlength{\marginenumerate}
\setlist[enumerate]{leftmargin=\marginenumerate}

\usepackage{graphicx}

\usepackage{amsthm}

\newtheorem{theorem}{Theorem}
\newtheorem{corollary}[theorem]{Corollary}
\newtheorem{problem}[theorem]{Problem}
\newtheorem*{problem*}{Problem}
\newtheorem{proposition}[theorem]{Proposition}

\newtheorem{lemma}[theorem]{Lemma}
\newtheorem{claim}{Claim}[theorem]
\newtheorem{observation}[theorem]{Observation}
\newtheorem*{observation*}{Observation}

\newenvironment{proofclaim}[1][{\it Proof of the claim. \hspace{0.066cm}}]%
	{\noindent {}{#1}{}}{ \strut\hfill $\lozenge$\vspace{2ex}}

\DeclareMathOperator{\td}{td}
\DeclareMathOperator{\vc}{vc}
\DeclareMathOperator{\fvs}{fvs}
\DeclareMathOperator{\pd}{pd}
\DeclareMathOperator{\tw}{tw}
\DeclareMathOperator{\pw}{pw}
\DeclareMathOperator{\torso}{torso}
\DeclareMathOperator{\dist}{dist}
\DeclareMathOperator{\term}{term}
\DeclareMathOperator{\init}{init}
\DeclareMathOperator{\GMh}{GMh}
\DeclareMathOperator{\p}{p}

\newcommand{\bigO}{\mathcal{O}}

\renewcommand{\phi}{\varphi}
\renewcommand{\epsilon}{\varepsilon}

\newcommand{\NN}{\mathbb{N}}

\newcommand{\ind}[1]{\left[#1\right]}

\renewcommand{\emph}{\textbf}

\renewcommand{\leq}{\leqslant}
\renewcommand{\geq}{\geqslant}

\renewcommand{\hat}{\widehat}

\newcommand{\cstref}[1]{\ref{#1}}

\title{Excluding a rectangular grid}
\author{
Cl\'ement Rambaud\thanks{Research supported by the French Agence Nationale de la Recherche under contract Digraphs ANR-19-CE48-0013.} \\
\small Universit\'e C\^ote d'Azur, CNRS, Inria, I3S\\[-0.8ex]
\small Sophia-Antipolis, France\\
\small\tt clement.rambaud@normalesup.org\\
~
}

\begin{document}

\maketitle

\begin{abstract}
    For every positive integer $k$,
    we define the $k$-treedepth
    as the largest graph parameter $\td_k$ satisfying
    \begin{enumerate*}
        \item $\td_k(\emptyset)=0$;
        \item $\td_k(G) \leq 1+ \td_k(G-u)$ for every graph $G$ and every vertex $u \in V(G)$; and 
        \item if $G$ is a $(<k)$-clique-sum of $G_1$ and $G_2$, then $\td_k(G) \leq \max \{\td_k(G_1),\td_k(G_2)\}$, for all graphs $G_1,G_2$.
    \end{enumerate*}
    This parameter coincides with treedepth if $k=1$, and with treewidth plus $1$ if $k \geq |V(G)|$.
    We prove that for every positive integer $k$,
    a class of graphs $\mathcal{C}$ has bounded $k$-treedepth if and only if there is a positive integer $\ell$
    such that for every tree $T$ on $k$ vertices, no graph in $\mathcal{C}$ contains $T \square P_\ell$ as a minor.

    This implies, for $k=1$,
    that a minor-closed class of graphs has bounded treedepth if and only if it excludes a path,
    for $k=2$, that a minor-closed class of graphs has bounded $2$-treedepth if and only if it excludes as a minor a ladder (Huynh, Joret, Micek, Seweryn, and Wollan; Combinatorica, 2021), and
    for large values of $k$, that a minor-closed class of graphs has bounded treewidth if and only if it excludes a grid (Grid-Minor Theorem, Robertson and Seymour; JCTB, 1986).
    
    As a corollary, we obtain the following qualitative strengthening of the Grid-Minor Theorem in the case of bounded-height grids.
    For all positive integers $k, \ell$, every graph that does not contain the $k \times \ell$ grid as a minor has
    $(2k-1)$-treedepth at most a function of $k$ and $\ell$.

    \medskip
    \noindent{}{\bf Keywords:} Graph minor, grid graph, tree decomposition, treedepth, treewidth.
\end{abstract}

\sloppy

\bigskip

The Grid-Minor Theorem, proved by Robertson and Seymour~\cite{Robertson1986},
characterizes classes of graphs having bounded treewidth in terms of excluded minors.
More precisely, it asserts that a class of graphs $\mathcal{C}$ has bounded treewidth if and only if there is a positive integer $\ell$ such that, 
for every $G \in \mathcal{C}$,  the $\ell \times \ell$ grid is not a minor of $G$
(see Section~\ref{sec:notation} for the definition of grid).
The ``only if'' part can be easily obtained from the fact that treewidth is monotone for the graph minor relation, and that grids have unbounded treewidth.
Hence the core of this theorem is the fact that graphs excluding the $\ell \times \ell$ grid as a minor have treewidth at most $f(\ell)$, for some function $f$.
A lot of research has been carried out to obtain good bounds on the smallest such function $f$.
While the original proof of Robertson and Seymour gives a superpolynomial upper bound on $f$,
Chekuri and Chuzhoy proved a first polynomial upper bound in~\cite{CC16}, 
which was later lowered down to $\ell^{9}\log(\ell)^{\bigO(1)}$ by Chuzhoy and Tan~\cite{Chuzhoy2021}.

Treedepth is a graph parameter larger than treewidth which
can be characterized as the largest graph parameter $\td$ satisfying
\begin{enumerate*}
    \item $\td(\emptyset)=0$\footnote{We denote by $\emptyset$ both the empty set and the graph with no vertices.}, \label{item:def_td_i}
    \item $\td(G) \leq \td(G-u)+1$ for every graph $G$ and for every $u \in V(G)$, and \label{item:def_td_ii}
    \item $\td(G) \leq \max\{\td(G_1),\td(G_2)\}$ if $G$ is the disjoint union of $G_1$ and $G_2$, for all graphs $G_1,G_2$. \label{item:def_td_iii}
\end{enumerate*}
A class of graphs $\mathcal{C}$ has bounded treedepth if and only if there is a positive integer $\ell$ such that
for every $G \in \mathcal{C}$, $P_\ell$, the path on $\ell$ vertices, is not a minor of $G$.
This analog of the Grid-Minor Theorem for treedepth can be easily obtained by elementary arguments, see~\cite[Proposition 6.1]{Nesetril2012}.

Huynh,  Joret,  Micek,  Seweryn, and Wollan~\cite{Huynh2021} recently introduced a variant of treedepth, called $2$-treedepth and denoted by $\td_2$,
which can be characterized as the largest graph parameter satisfying \ref{item:def_td_i}-\ref{item:def_td_iii} of the definition of treedepth, together with the property
\begin{enumerate*}[resume]
    \item $\td_2(G) \leq \max\{\td_2(G_1),\td_2(G_2)\}$ if $G$ can be obtained from the disjoint union of $G_1$ and $G_2$ by identifying one vertex in $G_1$ with one vertex in $G_2$,
        for all graphs $G_1,G_2$.
\end{enumerate*}
Huynh, Joret,  Micek,  Seweryn, and Wollan
proved that a class of graphs has bounded $2$-treedepth if and only if for some positive integer $\ell$,
every graph in this class excludes the $2 \times \ell$ grid, as a minor.
While initially defined to describe the structure of graphs excluding a $2 \times \ell$ grid as a minor,
it turns out that this parameter plays an important role in the study of weak coloring numbers and centered chromatic numbers
in minor-closed classes of graphs, see~\cite{wcol_paper,pcentered_coloring_of_Kt_minor_free_graphs,hodor2026centeredcoloringsweakcoloring}.

It is noteworthy that treewidth has a similar characterization:
the treewidth is the largest graph parameter $\tw$ satisfying
\begin{enumerate*}
    \item $\tw(\emptyset)=-1$, \label{item:def_tw_i}
    \item $\tw(G) \leq \tw(G-u)+1$ for every graph $G$ and for every $u \in V(G)$, and \label{item:def_tw_ii}
    \item $\tw(G) \leq \max\{\tw(G_1),\tw(G_2)\}$ if $G$ is a clique-sum of $G_1$ and $G_2$, for all graphs $G_1,G_2$. \label{item:def_tw_iii}
\end{enumerate*}
Here, a clique-sum of two graphs $G_1$ and $G_2$ is any graph $G$ obtained from the disjoint union of $G_1$ and $G_2$ 
by identifying a clique in $G_1$ with a clique in $G_2$ of the same size.
Furthermore, for every positive integer $k$, if these cliques have size less than $k$, then we say that $G$ is a $(<k)$-clique-sum of $G_1$ and $G_2$.
In particular, a $(<1)$-clique-sum of $G_1$ and $G_2$ is always the disjoint union of $G_1$ and $G_2$.
This notion leads to the following generalization of treewidth, which is the subject of this paper.

\medskip

For every $k \in \mathbb{N}_{>0} \cup \{+\infty\}$,
we define the \emph{$k$-treedepth} as the largest graph parameter $\td_k$ satisfying
\begin{enumerate}
    \item $\td_k(\emptyset)=0$, \label{item:def_tdk_i}
    \item $\td_k(G) \leq 1+ \td_k(G-u)$ for every graph $G$ and every vertex $u \in V(G)$, and \label{item:def_tdk_ii}
    \item $\td_k(G) \leq \max \{\td_k(G_1),\td_k(G_2)\}$ if $G$ is a $(<k)$-clique-sum of $G_1$ and $G_2$, for all graphs $G_1,G_2$. \label{item:def_tdk_iii}
\end{enumerate}
This gives a well-defined parameter because if $\mathcal{P}$ is the family of all the graph parameters satisfying \ref{item:def_tdk_i}-\ref{item:def_tdk_iii},
then $\td_k \colon G \mapsto \max_{\mathrm{p} \in \mathcal{P}} \mathrm{p}(G)$ also satisfies \ref{item:def_tdk_i}-\ref{item:def_tdk_iii}.
A more explicit definition is given in Section~\ref{subsec:k_dismantable_td}.

Note that the $1$-treedepth is the usual treedepth, the $2$-treedepth coincides with the homonymous parameter $\td_2$
introduced by Huynh,  Joret,  Micek,  Seweryn, and Wollan in~\cite{Huynh2021},
and that the $+\infty$-treedepth coincides with the treewidth plus $1$.
Hence, for every graph $G$,
\[
\td(G) = \td_1(G) \geq \td_2(G) \geq \dots \geq \td_{+\infty}(G) = \tw(G) + 1.
\]

\begin{figure}[ht]
    \centering
    \includegraphics{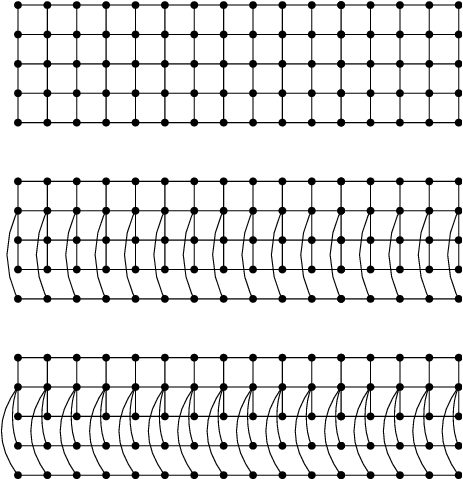}
    \caption{The three families of obstructions for $5$-treedepth given by Theorem~\ref{thm:main}.}
    \label{fig:obstruct_td5}
\end{figure}

Our main result is a characterization of classes of graphs having bounded $k$-treedepth in terms of excluded minors, 
which generalizes the aforementioned results on treedepth, $2$-treedepth, and treewidth.
\begin{theorem}\label{thm:main}
    Let $k$ be a positive integer.
    A class $\mathcal{C}$ of graphs has bounded $k$-treedepth if and only if there exists an integer $\ell$ such that
    for every tree $T$ on $k$ vertices, 
    no graph in $\mathcal{C}$ contains $T \square P_\ell$ as a minor.
\end{theorem}
Here, $\square$ denotes the Cartesian product: for all graphs $G_1,G_2$, 
$G_1 \square G_2 = \big(V(G_1) \times V(G_2), \{(u,v)(u',v') \mid (u=u' \text{ and } vv' \in E(G_2)) \text{ or } (uu' \in E(G_1) \text{ and } v=v')\}\big)$;
and $P_\ell$ denotes the path on $\ell$ vertices.
See Figure~\ref{fig:obstruct_td5} for some examples of graphs of the form $T \square P_\ell$.

Theorem~\ref{thm:main} gives for every positive integer $k$ a characterization of classes of graphs having bounded $k$-treedepth in terms of excluded minors, 
which was previously known only for $k=1$ and $k=2$~\cite{Huynh2021}.
For example, Theorem~\ref{thm:main}
applied to $k=5$ implies that for every fixed positive integer $\ell$, 
graphs excluding $T_1 \square P_\ell, T_2 \square P_\ell$, and $T_3 \square P_\ell$ as minors have bounded $5$-treedepth,
where $T_1,T_2,T_3$ are the three trees on $5$ vertices up to isomorphism (see Figure~\ref{fig:obstruct_td5}).
This is optimal since each of the three families $\{T_1 \square P_\ell\}_{\ell \geq 1},\{T_2 \square P_\ell\}_{\ell \geq 1},\{T_3 \square P_\ell\}_{\ell \geq 1}$
has unbounded $k$-treedepth.

A similar result was obtained by Geelen and Joeris~\cite{geelen2016_a_generalization_fo_the_grid_minor_theorem} for the parameter called $k$-treewidth,
which is the minimum width of a tree decomposition with adhesion smaller than $k$ (see Section~\ref{sec:prelim:minors_and_tree_decompositions} for the definition of tree decomposition and adhesion).
However, their proof relies on very different techniques.
Indeed, the argument of~\cite{geelen2016_a_generalization_fo_the_grid_minor_theorem} is a generalization of a proof of the Grid-Minor Theorem,
while here, the structure we are looking for is more complex, but we use the Grid-Minor Theorem as a black box.

\medskip

It is not difficult to see that, 
for every tree $T$ with $|V(T)|$ big enough compared to $k$, 
graphs in $\{T \square P_\ell \mid \ell \in \NN\}$ contain arbitrarily long 
grids of height $\ell$ as minors for some unbounded function $f\colon \NN \to \NN$: 
If $|V(T)|$ large enough,
then $T$ contains either $P_k$ of $K_{1,k}$ as a minor, 
and in both cases it follows that graphs in $\{T \square P_\ell \mid \ell \in \NN\}$ 
contain arbitrarily long grids of height $k$ as minors.
Therefore,
Theorem~\ref{thm:main} implies that graphs excluding a fixed grid of height $k$
as a minor have bounded $f(k)$-treedepth, for some function $f \colon \NN \to \NN$.

Refining the argument, we obtain that taking $f(k) = 2k-1$ is enough, 
and so we deduce the following corollary of Theorem~\ref{thm:main},
which is a structural property for the graphs excluding the $k \times \ell$ grid as a minor, when $k$ is small compared to $\ell$.
This is a qualitative strengthening of the Grid-Minor Theorem by Robertson and Seymour~\cite{Robertson1986}.

\begin{corollary}\label{corollary:excluding_only_the_kl_grid}
    There is a function $f\colon \NN^2 \to \NN$~\footnote{We denote by $\NN$ the set of nonnegative integers.} such that for all positive integers $k,\ell$,
    for every graph $G$, if the $k \times \ell$ grid is not a minor of $G$, then
    \[
    \td_{2k-1}(G) \leq f(k,\ell).
    \]
\end{corollary}

This corollary is tight up to the factor of $2$ and the value of $f(\cdot,\cdot)$:
the family of the grids of height $k-1$ have unbounded $(k-1)$-treedepth but
exclude the $k \times \ell$ grid for every $\ell \geq k$.
However, we do not know whether there exists $c < 2$ such that
graphs excluding the $k \times \ell$ grid have bounded $(ck-1)$-treedepth 
for every fixed positive integer $\ell$.

In~\cite{Biedl2019}, Biedl, Chambers,  Eppstein, De Mesmay, and Ophelders defined the \emph{grid-major height} of a planar graph $H$, denoted by $\GMh(H)$, 
as the smallest integer $k$ such that $H$ is a minor of the $k \times \ell$ grid for some integer $\ell$. 
With this definition, Corollary~\ref{corollary:excluding_only_the_kl_grid} implies the following.

\begin{corollary}
    There is a function $f\colon \NN \to \NN$ such that the following holds.
    For every planar graph $H$ with at least one vertex and for every $H$-minor-free graph $G$,
    \[
    \td_{2\GMh(H)-1}(G) \leq f(|V(H)|).
    \]
\end{corollary}
Moreover, $\GMh$ is the right parameter in such a result since,
for every graph $H$ with $\GMh(H) \geq 2$,
$H$-minor-free graphs have unbounded $\td_{\GMh(H)-1}$,
as witnessed by the family of the $(\GMh(H)-1) \times \ell$ grids for $\ell>0$.
This gives a better understanding of $H$-minor-free graphs when $H$ is planar and has a specific structure.
In the same spirit, Dujmovi\'c, Hickingbotham, Hodor, Joret, La, Micek, Morin, Rambaud, and
Wood~\cite{Dujmovi2024} recently gave a product structure for $H$-minor-free graphs when $H$ is planar and has bounded treedepth.

\medskip

We also define an analogue of $k$-treedepth for path decompositions, instead of tree decompositions, that we call $k$-pathdepth and denoted by $\pd_k$.
Informally, the \emph{$k$-pathdepth} $\pd_k$ is the largest parameter satisfying \ref{item:def_td_i} and \ref{item:def_td_ii} of the definition of $k$-treedepth,
and an analog of \ref{item:def_td_iii} ensuring that the clique-sums are made in the way of a path. 
See Section~\ref{sec:outline} for a formal definition.
This gives a family of parameters satisfying for every graph $G$
\[
    \td(G) = \pd_1(G) \geq \pd_2(G) \geq \dots \geq \pd_{+\infty}(G) = \pw(G) + 1
\]
and
\[
    \pd_k(G) \geq \td_k(G)
\]
for every positive integer $k$.
Hence, if a class of graphs has bounded $k$-pathdepth, then it has bounded $k$-treedepth and bounded pathwidth.
Quite surprisingly, these two necessary conditions are also sufficient.

\begin{theorem}\label{thm:main_pdk}
    Let $k$ be a positive integer and let $\mathcal{C}$ be a class of graphs.
    The following are equivalent.
    \begin{enumerate}[label=(\arabic*)]
        \item $\mathcal{C}$ has bounded $k$-pathdepth. \label{item:main_pdk_i}
        \item $\mathcal{C}$ has bounded pathwidth and bounded $k$-treedepth. \label{item:main_pdk_ii}
        \item There is an integer $\ell$ such that for every $G \in \mathcal{C}$, $\pw(G) \leq \ell$, and for every tree $T$ on $k$ vertices,
            $T \square P_\ell$ is not a minor of $G$. \label{item:main_pdk_iii}
    \end{enumerate}
\end{theorem}
Our proof actually shows the equivalence between Items~\ref{item:main_pdk_i} and~\ref{item:main_pdk_iii},
while the equivalence between Items~\ref{item:main_pdk_ii} and~\ref{item:main_pdk_iii} follows from Theorem~\ref{thm:main}.

\paragraph*{Outline of the paper.}

In Section~\ref{sec:outline} the main ideas behind the proof 
of Theorem~\ref{thm:main} are presented, 
and in particular a proof of the implication \ref{item:main_pdk_iii}$\Rightarrow$\ref{item:main_pdk_i} of Theorem~\ref{thm:main_pdk}.
In Section~\ref{sec:preliminaries}, we introduce some notation and basic tools to study graph minors and tree decompositions.
Section~\ref{sec:kladders_have_large_k_treedepth} contains the proof that graphs of the form $T \square P_\ell$ have unbounded $k$-treedepth as $\ell \to +\infty$,
for every tree $T$ on $k$ vertices. This constitutes the easy direction in Theorem~\ref{thm:main}.
Sections~\ref{sec:finding_a_kladder} to \ref{sec:finding_a_suitable_tree_decomposition} are devoted to the proof of technical lemmas used
in the final proof of Theorem~\ref{thm:main}, which is in Section~\ref{sec:graphs_excluding_k_ladders_have_bounded_tdk}.
Then, in Section~\ref{sec:token_sliding}, we show that every graph of the form $T \square P_L$ for some tree $T$ on $2k-1$ vertices
contains the $k \times \ell$ grid as a minor, if $L$ is large enough compared to $\ell$.
Together with Theorem~\ref{thm:main}, this gives Corollary~\ref{corollary:excluding_only_the_kl_grid}.
Finally, we conclude in Section~\ref{sec:conclusion} by introducing a wide family of graph parameters,
including notably treewidth, $k$-treedepth, vertex cover number, feedback vertex set number,
and we pose the problem of finding universal obstructions for these parameters.

\section{Outline of the proof}\label{sec:outline}

In this section, we present the main ideas behind the proof of Theorem~\ref{thm:main}.
The definition of $k$-treedepth given in the introduction relies on the notion
of clique-sum.
So first, let us define exactly what is a clique-sum.
Given two graphs $G_1,G_2$, a \emph{clique-sum of $G_1$ and $G_2$} is a graph $G$ obtained by choosing a clique $C_1$ in $G_1$, a clique $C_2$ in $G_2$,
and a bijection $\phi\colon V(C_1) \to V(C_2)$, and by identifying in the disjoint union of $G_1$ and $G_2$ the vertex $u$ with $\phi(u)$, for every $u \in V(C_1)$,
and possibly removing some edges in the resulting clique.
Given an integer $k$, a \emph{$(<k)$-clique-sum of $G_1$ and $G_2$} is such a graph $G$, with the additional condition $|C_1|=|C_2|<k$.

The notion of clique-sum is closely related to the notion of tree decomposition.
This allows us to give an alternative definition of $k$-treedepth of a graph $G$,
as the minimum width (plus one) of a ``$k$-dismantable'' tree decomposition of $G$.
This definition will also be recursive, and should be seen as the translation
of the definition of $k$-treedepth in terms of tree decompositions.
Recall that a \emph{tree decomposition} of a graph $G$ is a pair $\mathcal{D} = \big(T,(W_x \mid x \in V(T))\big)$, where $T$ is a tree and $W_x \subseteq V(G)$ for every $x \in V(T)$, satisfying
\begin{enumerate}
    \item for every $u \in V(G)$, $\{x \in V(T) \mid u \in W_x\}$ induces a non-empty connected subtree of $T$, and \label{item:def_tree_decomposition:i}
    \item for every $uv \in E(G)$, there exists $x \in V(T)$ such that $u,v \in W_x$. \label{item:def_tree_decomposition:ii}
\end{enumerate}
More generally, when no graph is specified, 
we call tree decomposition any pair $\big(T,(W_x \mid x \in V(T))\big)$ satisfying \ref{item:def_tree_decomposition:i}.
For every vertex $x \in V(T)$, we call the set $W_x$ the \emph{bag} of $x$ in this tree decomposition.
For every edge $xy \in E(T)$, we call the set $W_x \cap W_y$ the \emph{adhesion} of $xy$ in this tree decomposition.
The \emph{width} of $\mathcal{D}$ is $\max_{x \in V(T)} |W_x|-1$.
The \emph{adhesion} of $\mathcal{D}$ is $\max_{xy \in E(T)} |W_x \cap W_y|$.
The \emph{treewidth} of $G$ is the minimum width of a tree decomposition of $G$.
For an introduction to tree decompositions and treewidth, see~\cite[Chapter 12]{Diestel2017}.

Let $\mathcal{D} = \big(T,(W_x \mid x \in V(T))\big)$ be a tree decomposition of a graph $G$.
Given a connected subgraph $F$ of $G$, its \emph{projection} on $T$ (according to $\mathcal{D}$) is the subtree of $T$ induced by 
$\{x \in V(T) \mid V(F) \cap W_x \neq \emptyset\}$.
By the properties of tree decompositions, the projection of a connected subgraph is connected.
By the Helly property for subtrees,
we obtain the following fundamental lemma.

\begin{lemma}[Statement~(8.7) in~\cite{Robertson1986}]\label{lemma:helly_property}
    Let $d$ be a positive integer.
    For every graph $G$, for every tree decomposition $\mathcal{D}$ of $G$, for every family $\mathcal{F}$ of connected subgraphs of $G$, one of the following holds:
    \begin{enumerate}[label=(\arabic*)]
        \item there are $d$ members of $\mathcal{F}$ whose projections are pairwise disjoint, or
        \item there is a set of at most $d-1$ bags of $\mathcal{D}$ whose union intersects every member of $\mathcal{F}$.
    \end{enumerate}
\end{lemma}

\medskip

A tree decomposition $\big(T,(W_x \mid x \in V(T))\big)$ is \emph{$k$-dismantable} if one of the following holds.
\begin{enumerate}
    \item $T$ has a single vertex $x$ and $W_x = \emptyset$; or
    \item there is a vertex $v \in \bigcup_{x \in V(T)} W_x$ such that $v \in W_x$ for every $x \in V(T)$, and $\big(T ,(W_x \setminus \{v\} \mid x \in V(T))\big)$ is $k$-dismantable; or
    \item there is an edge $xy \in E(T)$ such that $|W_x \cap W_y| < k$ and  
         the tree decompositions $\big(T_{x \mid y}, (W_z \mid z \in V(T_{x \mid y}))\big)$ and $\big(T_{y \mid x}, (W_z \mid z \in V(T_{y \mid x}))\big)$ are $k$-dismantable.
\end{enumerate}
Here, $T_{x \mid y}$ denotes the connected component of $x$ in $T \setminus xy$, for every edge $xy$ of $T$.
A straightforward induction shows the following observation.
See Section~\ref{subsec:k_dismantable_td} for a proof.
\begin{observation*}
    Let $k$ be a positive integer, let $G$ be a graph.
    The $k$-treedepth of $G$ is $1$ plus the minimum width of a $k$-dismantable tree decomposition of $G$.
\end{observation*}

This definition of $k$-treedepth based on tree decompositions will be useful in particular to show
that $\td_k(T \square P_\ell) \to +\infty$ as $\ell \to +\infty$ for every tree $T$ on $k$ vertices.
This proof does not represent any particular difficulty, see Section~\ref{sec:kladders_have_large_k_treedepth}.

Hence, to prove Theorem~\ref{thm:main}, it remains to show that there is a function $f\colon \NN^2 \to \NN$ such that
for every graph $G$ and every positive integer $\ell$,
if $T \square P_\ell$ is not a minor of $G$ for every tree $T$ on $k$ vertices, then $\td_k(G) \leq f(k,\ell)$.
A first ingredient is the following lemma, which shows how to extract a minor of the form $T \square P_\ell$ for some tree $T$ on $k$ vertices,
from $k$ pairwise disjoint paths $Q_1, \dots, Q_k$ and sufficiently many pairwise disjoint connected subgraphs, each intersecting $Q_i$ for every $i \in [k]$.

\begin{lemma}\label{lemma:finding_kladder_from_a_mess_outline}    
    There exists a function $f_{\cstref{lemma:finding_kladder_from_a_mess_outline}}\colon \NN^2 \to \NN$ such that the following holds.
    Let $k$ be a positive integer.
    For every graph $G$, for all $k$ pairwise disjoint paths $Q_1,\dots, Q_k$ in $G$,
    if there exist $f_{\cstref{lemma:finding_kladder_from_a_mess_outline}}(k,\ell)$ pairwise disjoint connected subgraphs of $G$
    each of them intersecting $V(Q_i)$ for every $i \in [k]$,
    then there is a tree $T$ on $k$ vertices such that
    $T \square P_\ell$ is a minor of $G$.
\end{lemma}
The proof of this lemma is given in Section~\ref{sec:finding_a_kladder}.
The fact that trees on $k$ vertices appear in this statement, and so in Theorem~\ref{thm:main}, can be easily explained.
The proof of Lemma~\ref{lemma:finding_kladder_from_a_mess_outline} actually yields a minor $H$ of $G$
with vertex set $[k] \times [k^{k-2} (\ell-1)+1]$ such that each ``row'' $\{(i,1), \dots, (i,k^{k-2} (\ell-1)+1)\}$ for $i \in [k]$ induces a path in $H$,
and each ``column'' $\{(1,j), \dots, (k,j)\}$ for $j \in [k^{k-2} (\ell-1)+1]$ induces a connected subgraph of $H$.
Now, each of these columns contains a spanning tree. By Cayley's formula and the pigeonhole principle, at least $\ell$ of these columns
have the same spanning tree $T$ when projecting on the first coordinate.
By ignoring the other columns, we deduce that $H$ contains $T \square P_\ell$ as a minor.

To illustrate how Lemma~\ref{lemma:finding_kladder_from_a_mess_outline} can be used, we prove Theorem~\ref{thm:main_pdk}.
The \emph{$k$-pathdepth} of a graph $G$, noted $\pd_k(G)$, is $1$ plus the minimum width of a $k$-dismantable path decomposition of $G$,
that is a $k$-dismantable tree decomposition $\big(T,(W_x \mid x \in V(T))\big)$ of $G$ for which $T$ is a path.
Together with the following unpublished result of Robertson and Seymour~\cite{RS_linked_path_decomposition} (see~\cite{Erde2018} for a proof), Lemma~\ref{lemma:finding_kladder_from_a_mess_outline} provides a simple proof of Theorem~\ref{thm:main_pdk}.
For a path $P$ and two vertices $x,y$ along $P$, we denote by $P[x,y]$
the subpath of $P$ with endpoints $x$ and $y$.

\begin{lemma}[Robertson and Seymour~\cite{RS_linked_path_decomposition}]\label{lemma:linked_path_decomposition}
    For every graph $G$,
    there exists a path decomposition $\mathcal{D} = \big(P,(W_x \mid x \in V(P))\big)$ of $G$ of width at most $\pw(G)$ such that for every distinct $x,y \in V(P)$, either
    \begin{enumerate}[label=(\arabic*)]
        \item there are $k$ pairwise disjoint $(W_x,W_y)$-paths in $G$, or \label{item:linked_path_decomposition:i}
        \item there is an edge $zz'$ in $P[x,y]$ with $|W_z \cap W_{z'}| < k$. \label{item:linked_path_decomposition:ii}
    \end{enumerate}
\end{lemma}

Again, knowing that $\pd_k(T \square P_\ell) \geq \td_k(T \square P_\ell)$ and $\td_k(T \square P_\ell) \to +\infty$ as $\ell \to +\infty$ for every tree $T$ on $k$ vertices,
and $\pd_k(G) \geq \pw(G)+1$ for every graph $G$, it remains to show the following.

\begin{theorem}\label{thm:pdk_technical}
    There is a function $f_{\cstref{thm:pdk_technical}}\colon \NN^3 \to \NN$ such that for every positive integer $\ell$ and every graph $G$, 
    if $T \square P_\ell$ is not a minor of $G$ for every tree $T$ on $k$ vertices, then
    \[
        \pd_k(G) \leq f_{\cstref{thm:pdk_technical}}(k, \ell, \pw(G)).
    \]
\end{theorem}

We will use the notation $G \cup F$ for $G$ a graph $G$ and
$F$ a subset of (unordered) pairs of vertices of $G$ 
to denote the graph obtained from $G$ by adding as edges all the pairs in $F$.
Also, for a set $X$, we denote by $\binom{X}{2}$ the set of all the (unordered)
pairs of elements of $X$.

\begin{figure}[ht]
    \centering
    \includegraphics{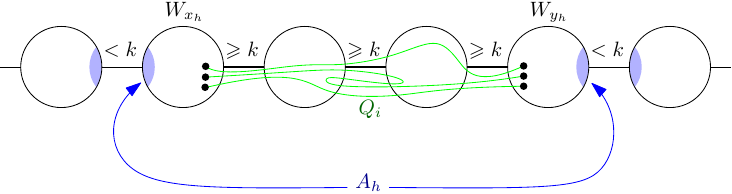}
    \caption{Illustration for the proof of Theorem~\ref{thm:pdk_technical}.}
    \label{fig:proof_pdk}
\end{figure}

\begin{proof}[Proof of Theorem~\ref{thm:pdk_technical}]
    Let $k, \ell$ be positive integers.
    For every nonnegative integer $p$,
    let
    \[
    f_{\cstref{thm:pdk_technical}}(k,\ell,p) = 
    \begin{cases}
        k &\textrm{if $p < k$,} \\
        f_{\cstref{thm:pdk_technical}}(k,\ell,p-1) + (f_{\cstref{lemma:finding_kladder_from_a_mess_outline}}(k,\ell)-1)(p+1) + 2(k-1) & \textrm{if $p \geq k$.}
    \end{cases}
    \]
    We proceed by induction on $\pw(G)$.
    Suppose that $T \square P_\ell$ is not a minor of $G$ for every tree $T$ on $k$ vertices.
    If $\pw(G) < k$, then every path decomposition of $G$ of width less than $k$ is $k$-dismantable, and so $\pd_k(G) \leq k = f_{\cstref{thm:pdk_technical}}(k,\ell,\pw(G))$.
    Now assume $\pw(G)\geq k$.
    Then, by Lemma~\ref{lemma:linked_path_decomposition},
    there exists a path decomposition $\mathcal{D} = \big(P,(W_x \mid x \in V(P))\big)$ of $G$ of width at most $\pw(G)$ 
    satisfying \ref{item:linked_path_decomposition:i} and \ref{item:linked_path_decomposition:ii} for every distinct $x,y \in V(P)$.
    Let $P^1, \dots, P^m$ be the connected components of $P \setminus \{xy \in E(P) \mid |W_x \cap W_y| < k \}$,
    in this order along $P$.

    Fix some $h \in [m]$. Let $x_h, y_h$ be respectively the first and last vertex of $P^h$.
    Let $A_h$ be the union of the at most two adhesions neighboring $P^h$,
    that is the set 
    \[
         \bigcup_{zz' \in E(P), z \in V(P^h), z' \not\in V(P^h)} W_z \cap W_{z'}.
    \]
    Note that $|A_h| \leq 2(k-1)$.
    Finally, let $G_h = G\ind{\bigcup_{z \in V(P^h)} W_z} \cup \binom{A_h}{2}$.
    Since $G$ is a $(<k)$-clique-sum of $G_1, \dots, G_m$, we want to bound $\pd_k(G_h)$.
    If $x_h=y_h$, then $|V(G_h)| = |W_{x_h}|$ and so the path decomposition of $G_h$ consisting of a single bag has width less than $\pw(G)+1 \leq f_{\cstref{thm:pdk_technical}}(k,\ell,\pw(G))$ and is $k$-dismantable.
    Now suppose that $x_h \neq y_h$.
    By the properties of $\mathcal{D}$, there are $k$ pairwise vertex-disjoint paths $Q_1, \dots, Q_k$ from $W_{x_h}$ to $W_{y_h}$ in $G$.
    See Figure~\ref{fig:proof_pdk}.
    Let $\mathcal{F}$ be the family of all the connected subgraphs $H$ of $G_h - A_h$ such that $V(H) \cap V(Q_i) \neq \emptyset$ for every $i \in [k]$.
    By Lemma~\ref{lemma:finding_kladder_from_a_mess_outline}, 
    there are no $f_{\cstref{lemma:finding_kladder_from_a_mess_outline}}(k,\ell)$ pairwise disjoint members of $\mathcal{F}$.
    Now, by Lemma~\ref{lemma:helly_property},
    there are $f_{\cstref{lemma:finding_kladder_from_a_mess_outline}}(k,\ell)-1$ bags of $\mathcal{D}$ whose union intersects every member of $\mathcal{F}$.
    Hence there is a set $X^h_0$ of at most $(f_{\cstref{lemma:finding_kladder_from_a_mess_outline}}(k,\ell)-1)(\pw(G)+1)$ 
    vertices in $G$ that intersects every member of $\mathcal{F}$.
    Let $X^h_1 = (X^h_0 \cap V(G_h)) \cup A_h$.
    For every connected component $C$ of $G_h-X^h_1$, $V(C)$ is disjoint from $X^h_0$, and so $C \not \in \mathcal{F}$.
    This implies that one the the paths in $\{Q_1, \dots, Q_k\}$ is disjoint from $V(C)$,
    and so there is a path in $G$ from $W_{x_h}$ to $W_{y_h}$ disjoint from $V(C)$.
    Since this path intersects $W_z$ for all $z \in V(P^h)$,
    $\big(P^h, (W_z \cap V(C) \mid z \in V(P^h))\big)$ is a path decomposition of $C$ of width at most $\pw(G)-1$.
    By induction, $\pd_k(C) \leq f_{\cstref{thm:pdk_technical}}(k,\ell,\pw(G)-1)$ and so 
    \begin{align*}
        \pd_k\left(G_h - X^h_1\right) 
        & \leq \max \{\pd_k(C) \mid \text{$C$ connected component of $G_h - X_1$}\}  \\
        & \leq f_{\cstref{thm:pdk_technical}}(k,\ell,\pw(G)-1).
    \end{align*}

    Consider now a $k$-dismantable path decomposition $\mathcal{D}_h = \big(R^h,(W^h_{0,x} \mid x \in V(R^h))\big)$ 
    of $G_h - X^h_1$ of width at most $f_{\cstref{thm:pdk_technical}}(k,\ell,\pw(G)-1)$.
    Then $\big(R^h,(W^h_{0,x} \cup X^h_1 \mid x \in V(R^h))\big)$ is a $k$-dismantable path decomposition of $G_h$ of width at most
    $f_{\cstref{thm:pdk_technical}}(k,\ell,\pw(G)-1) + (f_{\cstref{lemma:finding_kladder_from_a_mess_outline}}(k,\ell)-1)(\pw(G)+1) + 2(k-1) = f_{\cstref{thm:pdk_technical}}(k,\ell,\pw(G))$
    such that $A_h$ is included in every bag.
    In both cases ($x_h=y_h$ and $x_h \neq y_h$), we obtained a $k$-dismantable path decomposition $\big(R^h, (W^h_x \mid x \in V(R^h))\big)$ of $G_h$
    whose bags have size at most $f_{\cstref{thm:pdk_technical}}(k,\ell,\pw(G))$, and with $A_h$ included in every bag.
    Let $R$ be the path obtained from the disjoint union of $R^1, \dots, R^m$, and adding an edge between the last vertex of $R^h$ and the first vertex of $R^{h+1}$, 
    for every $h \in [m-1]$, and let $W'_x = W^h_x$ for every $x \in V(R^h)$, for every $h \in [m]$.
    Then $\big(R,(W'_x \mid x \in V(R))\big)$ is a $k$-dismantable path decomposition of $G$ whose bags have size at most
    $f_{\cstref{thm:pdk_technical}}(k,\ell,\pw(G))$.
    This proves the theorem.
\end{proof}

\medskip

This proof already contains several of the main ideas behind the proof of Theorem~\ref{thm:main}.
The general approach to bound the $k$-treedepth of a graph $G$ is to decompose $G$ as a $(<k)$-clique-sum of graphs $G_1, \dots, G_m$
such that for every $i \in [m]$, there exists $X_i \subseteq G_i$ of bounded size such that we can bound the $k$-treedepth of $G_i-X_i$ (typically by induction).
Hence a first step is to identify these graphs $G_i$.
In the bounded pathwidth case, this was done using Lemma~\ref{lemma:linked_path_decomposition}, 
the clique-sums being then given by the adhesions of size less than $k$ in the obtained path decomposition.
In the case $k=2$, there is a very natural way to find these clique-sums:
$G$ is a $(<2)$-clique-sum of its blocks\footnote{The \emph{blocks} of a graph $G$ are the maximal subgraphs of $G$ which consists of either a single vertex, a graph consisting of a single edge, or a $2$-connected subgraph.}, as used in~\cite{Huynh2021}.
For the general case, we will use a version of Lemma~\ref{lemma:linked_path_decomposition} for tree decompositions,
which is inspired by a classical result by Thomas~\cite{Thomas1990} on the existence of ``lean'' tree decompositions. 
This initial tree decomposition will be built in Section~\ref{sec:finding_a_suitable_tree_decomposition}. 
The techniques used in this section are inspired by a simple proof of Thomas' result by Bellenbaum and Diestel~\cite{BELLENBAUM2002}.

\medskip

An issue arising from this approach when considering tree decompositions instead of path decompositions is the following.
When we say that $G$ is a $(<k)$-clique-sums of $G_1, \dots, G_m$,
the graphs $G_i$ are not necessarily minors of $G$, and so they may contain a graph of the form $T \square P_\ell$ for some tree $T$ on $k$ vertices.
In the bounded pathwidth case, the graphs $G_h$ for $h \in [m]$ had just two additional cliques, corresponding to the right and left adhesions, 
and so we just removed these cliques (the set $A_h$) to obtain a subgraph of $G$.
This is not possible in the general case since we can have arbitrarily many such cliques.
To solve this issue, we will build $G_1, \dots, G_m$ with $V(G_h) \subseteq V(G)$ for every $i \in [m]$,
such that, while $G_h$ might not be a minor of $G$, the ``connectivity'' in $G_h$ is not higher than in $G$,
meaning that for every $Z_1, Z_2 \subseteq V(G_h)$ both of size $i$, if there are $i$ pairwise disjoint paths between $Z_1$ and $Z_2$ in $G_h$, then
there are $i$ pairwise disjoint paths between $Z_1$ and $Z_2$ in $G$.
Hence, in the main proof, we will work on such an ``almost minor'' of $G$, and every time we will need the assumption that there is no minor of the form $T \square P_\ell$
for a tree $T$ on $k$ vertices, we will come back to the original graph using this connectivity property.
This notion is developed in Section~\ref{sec:nice_sets}.

\medskip

Another important idea in the proof of Theorem~\ref{thm:pdk_technical} is the induction on $\pw(G)$.
For the general case, we want to proceed by induction on $\tw(G)$, which is bounded by the Grid-Minor Theorem.
However, it is in practice harder to decrease the treewidth than the pathwidth.
To circumvent this issue, we need a notion of ``partial tree decomposition''.
We note that some very similar notions were recently used, for example, in~\cite{Eiben2021,quickly_excluding_an_apex_forest,JansenSwennnenhuis2024}.
Note also that this notion is, while closely related, different from the
``bidimensionality'' as defined by Thilikos and Wiederrecht in~\cite{GMthroughBidim}.

\begin{figure}[ht]
    \centering
    \includegraphics{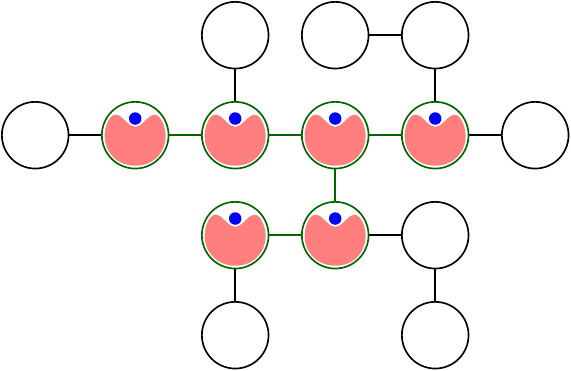}
    \caption{A tree decomposition of $(G-u,S)$ of width at most $\tw(G)-1$.
        The blue vertex depicts $u$,
        the tree $T' = T[\{x \in V(T) \mid u \in W_x\}]$ is depicted in green,
        and $S = \bigcup_{x \in V(T')} W_x \setminus\{u\}$ is the union of the red sets.}
    \label{fig:tree_decomposition_GS}
\end{figure}

Let $G$ be a graph and let $S \subseteq V(G)$.
A \emph{tree decomposition of $(G,S)$} is a tree decomposition $\big(T,(W_x \mid x \in V(T))\big)$ of $G[S]$
such that for every connected component $C$ of $G-S$, there exists $x \in V(T)$ such that $N_G(V(C)) \subseteq W_x$ 
(the notation $N_G(V(C))$ refers to the set $\{u \in V(G)\setminus V(C) \mid \exists v \in V(C), uv \in E(G)\}$).
Now, we will just assume that there is a tree decomposition of $(G,S)$ of small width,
and build by induction on this tree decomposition (actually on its adhesion) a $k$-dismantable tree decomposition of $(G,S')$ of bounded width, for some $S' \subseteq V(G)$ containing $S$.
The point of reinforcing the induction hypothesis in this way is that it is much easier to decrease the width or adhesion of a tree decomposition of a subset of $V(G)$ than for the full $V(G)$.
A crucial example is the following:
let $G$ be a graph, 
let $ \mathcal{D} = \big(T,(W_x \mid x \in V(T))\big)$ be a tree decomposition of $G$,
let $u \in V(G)$,
and let $S = \bigcup_{x \in V(T), u \in W_x} W_x \setminus \{u\}$.
Now, for $T' = T[\{x \in V(T) \mid u \in W_x\}]$,
$\big(T', (W_x \setminus \{u\} \mid x \in V(T'))\big)$ is a tree decomposition of $(G-u,S)$
whose width (resp. adhesion) is smaller than the width (resp. adhesion) of $\mathcal{D}$.
See Figure~\ref{fig:tree_decomposition_GS}.
We can now call the induction hypothesis on $(G-u,S)$, and thus decompose a superset of $N_G(u)$ in $G-u$.
This idea is used in the final proof of Theorem~\ref{thm:main}
to build a ``path partition'' (also called layering) satisfying some properties 
ensuring that any long enough such path partition contains one of the obstructions as a minor.
Morally, each part will give a column of a grid-like minor $T \square P_\ell$, while the rows will be given by $k$ disjoint paths from the first part to the last part.
This part of the argument is presented in Section~\ref{sec:graphs_excluding_k_ladders_have_bounded_tdk}.

\section{Preliminaries}\label{sec:preliminaries}

\subsection{Notation}\label{sec:notation}
We will use the following standard notation.
Let $k$ be a positive integer and let $X$ be a set.
We denote by $[k]$ the set $\{1, \dots, k\}$,
and by $\binom{X}{k}$ the family of all the subsets of $X$ of size $k$.

Graphs considered in this paper are finite, undirected, with no loops, and no multiple edges.
Let $G$ be a graph.
We denote by $V(G)$ its vertex set and by $E(G)$ its edge set.
If $u \in V(G)$, then $N_G(u) = \{v \in V(G) \mid uv \in E(G)\}$. 
If $U \subseteq V(G)$, then $N_G(U) = \bigcup_{u \in U} (N_G(u) \setminus U)$.
If $U \subseteq V(G)$, we say that $U$ is \emph{connected in $G$} if there is a connected component of $G$ containing $U$.
If $X \subseteq V(G)$, we denote by $G-X$ the subgraph of $G$ induced by $V(G)\setminus X$.
If $F \subseteq \binom{V(G)}{2}$, we denote by $G \cup F$ the graph obtained from $G$ by adding an edge between every pair of vertices in $F$,
and by $G \setminus F$ is the graph obtained from $G$ by removing every edge in $F \cap E(G)$.
A \emph{path partition} of $G$ is a sequence of $(V_0, \dots, V_{\ell})$ of pairwise disjoint nonempty subsets of $V(G)$,
such that $\bigcup_{0 \leq i \leq \ell} V_i = V(G)$, and for every edge $uv$ of $G$, if $i,j \in \{0, \dots, \ell\}$ are the indices
such that $u \in V_i$ and $v \in V_j$, then $|i-j| \leq 1$.

Let $G$ be a graph.
A \emph{path} in $G$ is a sequence $P=(p_1, \dots, p_\ell)$ of pairwise distinct vertices of $G$ such that $p_i p_{i+1} \in E(G)$ for every $i \in [\ell-1]$.
Then we say that $P$ is a $(p_1, p_\ell)$-path in $G$, or simply a path from $p_1$ to $p_\ell$.
The \emph{length} of $P$ is $\ell-1$.
We denote by $\init(P)$ its \emph{initial vertex} $p_1$, and by $\term(P)$ its \emph{terminal vertex} $p_\ell$.
If $p_1 \in U$ and $p_\ell \in U'$ for some $U,U' \subseteq V(G)$, we say that $P$ is a $(U,U')$-path, or a path from $U$ to $U'$.
Note that paths are thus oriented, but we will sometimes identify such a path $P$ with the subgraph $(V(P), \{p_i p_{i+1} \mid i \in [\ell-1]\})$.
A \emph{$(U,U')$-cut} is a set $X$ of vertices in $G$ such that every $(U,U')$-path intersects $X$.
Recall that Menger's Theorem asserts that the minimum size of a $(U,U')$-cut is equal to the maximum number of pairwise vertex-disjoint $(U,U')$-paths.

Let $G_1,G_2$ be two graphs.
The \emph{Cartesian product} of $G_1$ and $G_2$, denoted by $G_1 \square G_2$, is the graph with vertex set $V(G_1) \times V(G_2)$ and edge set $\{(u_1,u_2)(v_1,v_2) \mid (u_1=v_1 \text{ and } u_2v_2 \in E(G_2)) \text{ or } (u_1v_1 \in E(G_1) \text{ and } u_2=v_2)\}$.
For every positive integer $\ell$, we denote by $P_\ell$ the path $(1, \dots, \ell)$, seen as the graph $([\ell], \{(i,i+1) \mid i \in [\ell-1]\})$.
We call the graph $P_k \square P_\ell$ the \emph{$k \times \ell$ grid}.

Let $T$ be a tree and let $u,v \in V(T)$.
We denote by $T[u,v]$ the (unique) path from $u$ to $v$ in $T$.
Similarly, $T]u,v] = T[u,v] - \{u\}, T[u,v[ = T[u,v] - \{v\}, T]u,v[ = T[u,v] - \{u,v\}$.
Moreover, if $uv$ is an edge of $T$, then we denote by $T_{u \mid v}$ the connected component of $u$ in $T \setminus \{uv\}$.

Throughout the paper, we say that two graphs $G_1,G_2$ are \emph{disjoint} if they are vertex-disjoint, that is $V(G_1) \cap V(G_2) = \emptyset$.

\subsection{Graph minors and tree decompositions}\label{sec:prelim:minors_and_tree_decompositions}
Let $H, G$ be two graphs.
We say that $H$ is a \emph{minor} of $G$ if $H$ can be obtained from $G$ by successive
vertex deletions, edge deletions, edge contractions.
Equivalently, $H$ is a minor of $G$ if and only if $G$ contains
a \emph{model of $H$}, that is a family $(B_x \mid x \in V(H))$ of pairwise disjoint subsets of $V(G)$ such that
\begin{enumerate}
    \item $G\ind{B_x}$ is connected for every $x \in V(H)$, and
    \item there is an edge in $G$ between $B_x$ and $B_y$ for every $xy \in E(H)$.
\end{enumerate}
If $H$ is not a minor of $G$, then we say that $G$ is \emph{$H$-minor-free}.

The notion of clique-sum behaves well with respect to the graph minor relation. 
In particular a minor of a clique-sum of $G_1$ and $G_2$ is a clique-sum of a minor of $G_1$ and a minor of $G_2$. 
Moreover, the notion of clique-sum gives the following characterization of treewidth.

\begin{proposition}[Folklore]
    The treewidth is the largest graph parameter $\tw$ satisfying
    \begin{enumerate}[label=(\alph*)]
        \item $\tw(\emptyset)=-1$,
        \item $\tw(G) \leq \tw(G-u) + 1$ for every graph $G$ and every $u \in V(G)$, and
        \item $\tw(G) \leq \max\{\tw(G_1),\tw(G_2)\}$ if $G$ is a clique-sum of $G_1$ and $G_2$.
    \end{enumerate}
\end{proposition}

Robertson and Seymour~\cite{Robertson1986} characterized minor-closed classes of graphs of bounded treewidth in terms of excluded minors.
This fundamental result, known as the Grid-Minor Theorem, will be used several times in this paper.
\begin{theorem}[Robertson and Seymour~\cite{Robertson1986}]\label{thm:grid_minor_thm}
    There is a function $f_{\cstref{thm:grid_minor_thm}}\colon \NN \to \NN$ such that for every integer $\ell$, for every graph $G$, if the $\ell \times \ell$ grid is not a minor of $G$,
    then
    \[
    \tw(G) < f_{\cstref{thm:grid_minor_thm}}(\ell).
    \]
\end{theorem}

Let $G$ be a graph and let $S$ be a non empty subset of $V(G)$.
A \emph{tree decomposition of $(G,S)$} is a pair $\mathcal{D} = \big(T, (W_x \mid x \in V(T))\big)$ such that
\begin{enumerate}
    \item $\mathcal{D}$ is a tree decomposition of $G[S]$, and
    \item for every connected component $C$ of $G-S$, there exists $x \in V(T)$ such that $N_G(V(C)) \subseteq W_x$.
\end{enumerate}
We denote by $\tw(G,S)$ the minimum width of a tree decomposition of $(G,S)$.

\subsection{\texorpdfstring{$k$}{k}-dismantable tree decompositions}\label{subsec:k_dismantable_td}
In this section, we recall the definition of $k$-dismantable tree decompositions, and we prove some of their properties.
Let $k$ be a positive integer.
A tree decomposition $\big(T,(W_x \mid x \in V(T))\big)$ is \emph{$k$-dismantable} if one of the following holds.
\begin{enumerate}[label={\normalfont(D\arabic*)}]
    \item $T$ has a single vertex $x$ and $W_x = \emptyset$; or
    \item there is a vertex $v \in V(G)$ such that $v \in W_x$ for every $x \in V(T)$, and $\big(T ,(W_x \setminus \{v\} \mid x \in V(T))\big)$ is $k$-dismantable; or
    \item there is an edge $xy \in E(T)$ such that $|W_x \cap W_y| < k$ and  
         the tree decompositions $\big(T_{x \mid y}, (W_z \mid z \in V(T_{x \mid y}))\big)$ and $\big(T_{y \mid x}, (W_z \mid z \in V(T_{y \mid x}))\big)$ are $k$-dismantable. \label{item:kdismantble_def:k-cliquesum}
\end{enumerate}

A straightforward induction shows the following.
\begin{observation}\label{obs:equivalence_ktreedepth_kdecomposable_treedecomposition}
    A graph $G$ has $k$-treedepth at most $t$ if and only if $G$ admits a $k$-dismantable tree decomposition of width less than $t$.
\end{observation}

\begin{proof}
    Let $\td'_k(G)$ be $1$ plus the minimum width of a $k$-dismantable tree decomposition of $G$, for every graph $G$.
    We want to prove that $\td'_k(G)=\td_k(G)$ for every graph $G$.
    First we show that $\td'_k(G) \leq \td_k(G)$.
    It is enough to show that 
    \begin{enumerate*}
        \item $\td'_k(\emptyset) = 0$;
        \item $\td'_k(G) \leq 1+\td'_k(G-u)$ for every $u \in V(G)$; and
        \item $\td'_k(G) \leq \max\{\td'_k(G_1), \td'_k(G_2)\}$ if $G$ is a $(<k)$-clique-sum of $G_1$ and $G_2$.
    \end{enumerate*}
    The first point is clear.
    The second one comes from the fact that 
    if $\big(T,(W_x \mid x \in V(T))\big)$ is a $k$-dismantable tree decomposition of $G-u$, 
    then $\big(T,(W_x \cup \{u\} \mid x \in V(T))\big)$ is a $k$-dismantable tree decomposition of $G$.
    For the third point, suppose that $G$ is obtained from the $(<k)$-clique-sum of two graphs $G_1$ and $G_2$ with $C = V(G_1) \cap V(G_2)$
    inducing a clique of size at most $k-1$ in both $G_1$ and $G_2$.
    Let $\big(T_i, (W_x \mid x \in V(T_i))\big)$ be a $k$-dismantable tree decomposition of $G_i$,
    for $i \in \{1,2\}$.
    Then, since $C$ is a clique in $G_i$,
    the projections of $G_i[\{u\}]$ for $u \in C$ pairwise intersect.
    Hence by the Helly property, there exists $x_i \in V(T_i)$ intersecting all these projections, that is $C \subseteq W_{x_i}$.
    Now, $\big(T_1 \cup T_2 \cup \{x_1x_2\}, (W_x \mid x \in V(T_1) \cup V(T_2))\big)$ is a $k$-dismantable tree decomposition of $G$.
    This proves that $\td'_k(G) \leq \td_k(G)$ for every graph $G$.

    Now we show by induction on $|V(G)| + |V(T)|$ that $\td_k(G) \leq \td'_k(G)$ for every graph $G$.
    Let $\big(T, (W_x \mid x \in V(T))\big)$ be a $k$-dismantable tree decomposition of $G$ of width $t-1$ for $t = \td'_k(G)$.
    If $T$ has a single vertex $x$ and $W_x = \emptyset$, then $t=0$ and $\td_k(G) = 0 \leq t$.
    If $v \in V(G)$ is such that $v \in W_x$ for every $x \in V(T)$, and $\big(T ,(W_x \setminus \{v\} \mid x \in V(T))\big)$ is $k$-dismantable,
    then $\td_k(G-v) \leq \td'_k(G-v) \leq t-1$ by the induction hypothesis, and so $\td_k(G) \leq 1 + \td_k(G-v) \leq t$.
    If there is an edge $xy \in E(T)$ such that $|W_x \cap W_y| < k$ and  
    the tree decompositions $\big(T_{x \mid y}, (W_z \mid z \in V(T_{x \mid y}))\big)$ and $(T_{y \mid x}, (W_z \mid z \in V(T_{y \mid x}))\big)$ are $k$-dismantable,
    then let $G_1 = G\left[\bigcup_{z \in V(T_{x \mid y})} W_z\right] \cup \binom{W_x \cap W_y}{2}$ 
    and $G_2 =G\left[\bigcup_{z \in V(T_{y \mid x})} W_z\right] \cup \binom{W_x \cap W_y}{2}$.
    Observe that $W_x \cap W_y = V(G_1) \cap V(G_2)$ induces a clique of size less than $k$ in both $G_1$ and $G_2$,
    and $G$ is a $(<k)$-clique-sum of $G_1$ and $G_2$.
    Moreover, the $k$-dismantable tree decompositions $\big(T_{x \mid y}, (W_z \mid z \in V(T_{x \mid y}))\big),\big(T_{y \mid x}, (W_z \mid z \in V(T_{y \mid x}))\big)$ 
    witness the fact that $\td'_k(G_1),\td'_k(G_2) \leq t$.
    Hence, since by the induction hypothesis $\td'_k(G_1),\td'_k(G_2) \leq t$, we have $\td_k(G) \leq \max\{\td_k(G_1),\td_k(G_2)\} \leq t$.
    This proves that $\td_k(G) \leq \td'_k(G)$.
\end{proof}

We now prove that $k$-treedepth is monotone for the graph minor relation.
\begin{lemma}\label{lemma:ktd_minor_monotone}
    Let $k$ be a positive integer and let $G, H$ be two graphs.
    If $H$ is a minor of $G$, then 
    \[
        \td_k(H) \leq \td_k(G).
    \]
\end{lemma}

\begin{proof}
    Given a class of graphs $\mathcal{C}$,
    we denote by $\mathbf{A}(\mathcal{C})$ the class of all the graphs $G$ such that $G \in \mathcal{C}$ or there exists $u \in V(G)$ such that
    $G-u \in \mathcal{C}$.
    We claim that
    \begin{equation*}
        \text{if $\mathcal{C}$ is minor-closed, then $\mathbf{A}(\mathcal{C})$ is minor-closed.}
    \end{equation*}
    Let $\mathcal{C}$ be a minor-closed class of graphs.
    If $G$ is a graph and $u \in V(G)$ such that $G-u \in \mathcal{C}$,
    then for every minor $H$ of $G$, if $(B_x \mid x \in V(H))$ is a model of $H$ in $G$, then
    either $u \not \in \bigcup_{x \in V(H)} B_x$, and so $(B_x \mid x \in V(H))$ is a model of $H$ in $G-u$, which implies $H \in \mathcal{C} \subseteq \mathbf{A}(\mathcal{C})$;
    or $u \in \bigcup_{x \in V(H)} B_x$, and then
    $(B_x \mid x \in V(H-x_0))$ is a model of $H - x_0$ in $G-u$, where $x_0$ is the unique vertex of $V(H)$ such that $u \in B_{x_0}$,
    which implies $H-x_0 \in \mathcal{C}$ and so $H \in \mathbf{A}(\mathcal{C})$.

    Let $k$ be a positive integer.
    We denote by $\mathbf{S}_k(\mathcal{C})$ the class of all the graphs $G$ such that $G$ is a $(<k)$-clique-sum
    of two graphs in $\mathcal{C}$.
    We claim that
    \begin{equation*}
        \text{if $\mathcal{C}$ is minor-closed, then $\mathbf{S}_k(\mathcal{C})$ is minor-closed.}
    \end{equation*}
    To see that, consider a graph $G$ which is a $(<k)$-clique-sum of two graphs $G_1$ and $G_2$ in $\mathcal{C}$.
    Namely $V(G) = V(G_1) \cup V(G_2)$, $E(G) \subseteq E(G_1) \cup E(G_2)$, and $V(G_1) \cap V(G_2)$ induces a clique
    of size at most $k-1$ in both $G_1$ and $G_2$.
    Let $H$ be a minor of $G$, and let $(B_x \mid x \in V(H))$ be a model of $H$ in $G$.
    For $i \in \{1,2\}$, let $X_i = \{x \in V(H) \mid B_x \cap V(G_i) \neq \emptyset\}$,
    and let $H_i = H[X_i] \cup \binom{X_1 \cap X_2}{2}$.
    Then, for each $i \in \{1,2\}$,
    $(B_x \cap V(G_i) \mid x \in X_i)$ is a model of $H_i$ in $G_i$, and so $H_i$ is a minor of $G_i$.
    In particular, $H_i \in \mathcal{C}$.
    Finally, $H$ is a $(<k)$-clique-sum of $H_1$ and $H_2$, and so $H \in \mathbf{S}_k(\mathcal{C})$.

    For every nonnegative integer $i$, let 
    \[
        \mathbf{S}_k^i(\mathcal{C}) = 
        \begin{cases}
            \mathcal{C} & \text{if $i=0$} \\
            \mathbf{S}_k(\mathbf{S}^{i-1}_k(\mathcal{C})) & \text{if $i>0$.}
        \end{cases}
    \]
    By induction on $i$, $\mathbf{S}_k^i(\mathcal{C})$ is minor-closed if $\mathcal{C}$ is minor-closed.
    Then $\big(\mathbf{S}_k^i(\mathcal{C})\big)_{i \geq 0}$ is an inclusion-wise monotone sequence of minor-closed classes of graphs,
    and so its union $\overline{\mathbf{S}}_k(\mathcal{C})$ is minor-closed.

    We can now show by induction that for every nonnegative integer $t$, the class $\{\td_k \leq t\}$ of all the
    graphs of $k$-treedepth at most $t$ is minor-closed.
    Indeed, when $t=0$, $\{\td_k \leq 0\} = \{\emptyset\}$ is clearly minor-closed.
    Moreover, if $t \geq 1$ and $\{\td_k \leq t-1\}$ is minor-closed, then the class
    $ \overline{\mathbf{S}}_k\big(\mathbf{A}(\{\td_k \leq t-1\})\big)$ is also minor-closed by the previous observations.
    To conclude, it is enough to show that
    \begin{equation}\label{eq:ktd_and_class_operations}
      \{\td_k \leq t\} = \overline{\mathbf{S}}_k\big(\mathbf{A}(\{\td_k \leq t-1\})\big).
    \end{equation}
    
    The inclusion $\overline{\mathbf{S}}_k\big(\mathbf{A}(\{\td_k \leq t-1\})\big) \subseteq \{\td_k \leq t\}$ follows from the definition of $k$-treedepth.
    We now prove the other inclusion.
    Let $G$ be a graph of $k$-treedepth at most $t$ and let $\big(T,(W_x \mid x \in V(T))\big)$ be a $k$-dismantable tree decomposition of $G$,
    which exists by Observation~\ref{obs:equivalence_ktreedepth_kdecomposable_treedecomposition}.
    We want to show that $G \in \overline{\mathbf{S}}_k\big(\mathbf{A}(\{\td_k \leq t-1\})\big)$.
    We proceed by induction on $|V(G)|+|V(T)|$. 
    If $V(G) = \emptyset$, then the result is clear.
    Now assume $V(G) \neq \emptyset$.
    By the definition of $k$-dismantable tree decompositions, one of the following holds.
    \begin{description}
        \item[Case 1.] There is a vertex $v \in V(G)$ such that $v \in W_x$ for every $x \in V(T)$, and $\big(T ,(W_x \setminus \{v\} \mid x \in V(T))\big)$ is $k$-dismantable.
            Then by Observation~\ref{obs:equivalence_ktreedepth_kdecomposable_treedecomposition}, $\td_k(G-u) \leq t-1$ and so $G \in \mathbf{A}(\{\td_k \leq t-1\}) \subseteq \overline{\mathbf{S}}_k\big(\mathbf{A}(\{\td_k \leq t-1\})\big)$.
        \item[Case 2.] There is an edge $xy \in E(T)$ such that $|W_x \cap W_y| < k$ and  
             the tree decompositions $\mathcal{D}_{x \mid y} = \big(T_{x \mid y}, (W_z \mid z \in V(T_{x \mid y}))\big)$ and $\mathcal{D}_{y \mid x} = \big(T_{y \mid x}, (W_z \mid z \in V(T_{y \mid x}))\big)$ are $k$-dismantable.
             Let $G_{x \mid y}$ (resp. $G_{y \mid x}$) be the graph $G[\bigcup_{z \in V(T_{x \mid y})} W_z] \cup \binom{W_x \cap W_y}{2}$ (resp. $G[\bigcup_{z \in V(T_{y \mid x})} W_z] \cup \binom{W_x \cap W_y}{2}$).
             The $k$-dismantable tree decompositions $\mathcal{D}_{x\mid y}$ and $\mathcal{D}_{y \mid x}$ witness the fact that respectively
             $G_{x \mid y}$ and $G_{y \mid x}$ have $k$-treedepth at most $t$.
             By the induction hypothesis, we deduce that $G_{x \mid y}, G_{y \mid x} \in \overline{\mathbf{S}}_k\big(\mathbf{A}(\{\td_k \leq t-1\})\big)$.
             Since $G$ is a $(<k)$-clique-sum of $G_{x \mid y}$ and $G_{y \mid x}$, this implies that $G \in \overline{\mathbf{S}}_k\big(\mathbf{A}(\{\td_k \leq t-1\})\big)$.
             This proves \eqref{eq:ktd_and_class_operations}, and concludes the proof of the lemma. \qedhere
    \end{description}
\end{proof}

Similarly, for every positive integer $k$, for every graph $G$ and for every nonempty set $S \subseteq V(G)$,
we denote by $\td_k(G,S)$ the integer $1$ plus the minimum width of a $k$-dismantable tree decomposition of $(G,S)$.

\begin{lemma}\label{lemma:combine_tree_decompositions_of_S}
    Let $m,t,k$ be positive integers.
    For every graph $G$ and for every disjoint sets $S_1, \dots, S_m \subseteq V(G)$,
    if $S = \bigcup_{i \in [m]} S_i$, then
    \begin{enumerate}[label=(\alph*)]
        \item $\tw(G,S) + 1 \leq \sum_{i \in [m]} (\tw(G-(S_1 \cup \dots \cup S_{i-1}),S_i)+1)$, and
        \item $\td_k(G,S) \leq \sum_{i \in [m]} \td_k(G-(S_1 \cup \dots \cup S_{i-1}),S_i)$.
    \end{enumerate}
\end{lemma}

\begin{figure}[ht]
    \centering
    \includegraphics{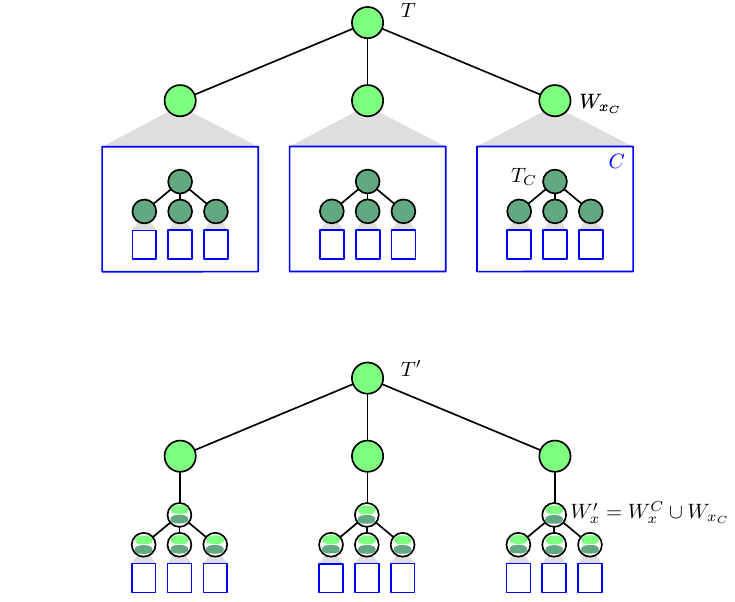}
    \caption{Illustration for the proof of Lemma~\ref{lemma:combine_tree_decompositions_of_S}.}
    \label{fig:combine_tree_decompositions}
\end{figure}

\begin{proof}
    Since for every $S \subseteq V(G)$,
    $\td_{|V(G)|}(G,S) = \tw(G,S)+1$, it is enough to show the second item.
    We proceed by induction on $m$.
    The result is clear for $m=1$.
    Now suppose $m>1$ and $\td_k(G,\bigcup_{i \in [m-1]}S_i) \leq \sum_{i \in [m-1]} \td_k(G-(S_1 \cup \dots \cup S_{i-1}),S_i)$.
    Let $\mathcal{D} = \big(T,(W_x \mid x \in V(T))\big)$ be a $k$-dismantable tree decomposition of $(G,\bigcup_{i \in [m-1]}S_i)$
    of width less than $\sum_{i \in [m-1]} \td_k(G-(S_1 \cup \dots \cup S_{i-1}),S_i)$.
    For every connected component $C$ of $G-\bigcup_{i \in [m-1]} S_i$,
    we have $\td_k(C,S_m \cap V(C)) \leq \td_k(G-(S_1 \cup \dots \cup S_{m-1}),S_m)$.
    Let $x_C \in V(T)$ be such that $N_G(V(C)) \subseteq W_{x_C}$.
    Hence there is a $k$-dismantable tree decomposition $\big(T_C, (W^C_x \mid x \in V(T_C))\big)$ of $(C,S_m \cap V(C))$ of width
    less than $\td_k(G,S_m)$.
    We suppose that the trees $T_C$ for $C$ connected components of $G-\bigcup_{i \in [m-1]} S_i$ and $T$ are pairwise disjoint.
    Let $T'$ be the tree obtained from the union of $T$ with all the $T_C$, by adding an edge between $x_C$ and an arbitrary vertex in $T_C$,
    for each connected component $C$ of $G-\bigcup_{i \in [m-1]}S_i$.
    For every $x \in V(T)$, let $W'_x = W_x$, and for every $x \in V(T_C)$, let $W'_x = W^C_x \cup W_{x_C}$,
    for each connected component $C$ of $G-\bigcup_{i \in [m-1]}S_i$.
    See Figure~\ref{fig:combine_tree_decompositions}.
    
    We claim that $\big(T', (W'_x \mid x \in V(T'))\big)$ is a $k$-dismantable tree decomposition of $(G,S)$.
    First, for every $u \in \bigcup_{i \in [m-1]} S_i$, $\{x \in V(T') \mid u \in W'_x\} = \{x \in V(T) \mid u \in W_x\} \cup \bigcup_{C} V(T_C)$ where $C$ ranges over all the connected components of $G-\bigcup_{i \in [m-1]} S_i$ such that $u \in W_{x_C}$. 
    Since $T[\{x \in V(T) \mid u \in W_x\}]$ is connected and by the definition of $T'$,
    we deduce that $T'[\{x \in V(T') \mid u \in W'_x\}]$ is connected.
    Moreover, for every $u \in S_m$, 
    if $C$ is the component of $u$ in $G-\bigcup_{i \in [m-1]} S_i$,
    then $T'[\{x \in V(T') \mid u \in W'_x\}] = T_C[\{x \in V(T_C) \mid u \in W^C_x\}]$, which is connected.
    This shows that for every $u \in S$, we have $T'[\{x \in V(T') \mid u \in W'_x\}]$ connected.
    For every $uv \in E(G[S])$, either $uv \in E(G[\bigcup_{i \in [m-1]} S_i])$ and
    so there is $x \in V(T)$ such that $u,v \in W_x = W'_x$,
    or there is a component $C$ of $G-\bigcup_{i \in [m-1]} S_i$ such that one of $u,v$, say $u$, is in $C$.
    Then, either $v \in V(C)$, and so there is $x \in V(T_C)$ such that $u,v \in W^C_x \subseteq W'_x$, or $v \in W_{x_C}$ and so $u,v \in  W^C_x \cup W_{x_C} = W'_x$ for any $x \in V(T_C)$ such that $u \in W^C_x$.
    Finally, for every connected component $C'$ of $G-S$,
    there is a connected component $C$ of $G-\bigcup_{i \in [m-1]} S_i$ such that
    $C' \subseteq C$.
    Then, $C'$ is a connected component of $C - (V(C) \cap S_m)$,
    and since $\big(T_C, (W^C_x \mid x \in V(T_C))\big)$ is a tree decomposition of $(C, V(C) \cap S_m)$,
    we deduce that there is $x \in V(T_C)$ such that $N_C(V(C')) \subseteq W^C_x$.
    It follows that $N_G(V(C')) \subseteq N_G(V(C)) \cup N_C(V(C')) \subseteq W_{x_C} \cup W^C_x = W'_x$.
    This proves that $\big(T', (W'_x \mid x \in V(T'))\big)$ is a tree decomposition of $(G,S)$.
    The fact that it is $k$-dismantable can be seen by induction on the structure of $\mathcal{D}$, using the fact that $\mathcal{D}$ is $k$-dismantable.
    
    Finally, $\big(T', (W'_x \mid x \in V(T'))\big)$ has width less than
    \begin{align*}
        \td_k(G-(S_1 \cup \dots \cup S_{m-1})&,S_m) + \sum_{i \in [m-1]} \td_k(G-(S_1 \cup \dots \cup S_{i-1}),S_i) \\
        &= \sum_{i \in [m]} \td_k(G-(S_1 \cup \dots \cup S_{i-1}),S_i).
    \end{align*}
    This proves the lemma.
\end{proof}

Note that Lemma~\ref{fig:combine_tree_decompositions} implies the following:
for every positive integer $k$, for every graph $G$,
for every $S_1, \dots, S_m \subseteq V(G)$,
\[
    \td_k\left(G,\ \textstyle\bigcup_{i \in [m]} S_i\right) \leq \sum_{i \in [m]} \td_k(G,S_i).
\]

\section{\texorpdfstring{$\{T \square P_\ell \mid \ell \geq 1\}$}{\{T □ P\_l | l > 0\}} has unbounded \texorpdfstring{$k$}{k}-treedepth}\label{sec:kladders_have_large_k_treedepth}

In this section, we prove that graphs of the form $T \square P_\ell$ for some tree $T$ on $k$ vertices 
have unbounded $k$-treedepth when $\ell$ tends to $+\infty$.
This proves the ``only if'' part of Theorem~\ref{thm:main}.

\begin{proposition}\label{prop:lower_bound}
    Let $k$ be a positive integer.
    For every tree $S$ on $k$ vertices,
    for every positive integer $N$,
    there is a positive integer $\ell$ such that
    \[
      \td_k(S \square P_\ell) > N.
    \]
\end{proposition}

\begin{proof}
    Fix a tree $S$ on $k$ vertices.
    Suppose for a contradiction that there exists a positive integer $N$ such that $\td_k(S \square P_\ell) \leq N$ for every positive integer $\ell$.
    By taking $N$ minimal, we can assume that, for some positive integer $\ell_0$,
    \[
        \td_k(S \square P_\ell) = N
    \]
    for every $\ell \geq \ell_0$.

    Consider $G = S \square P_{(2\ell_0 + 2)N}$.
    We denote by $R_1, \dots, R_k$ the copies of $P_{(2\ell_0+1)N}$ in $G$, that is the subgraphs $G\big[\{(x,i) \mid i \in [(2\ell_0+2)N]\}\big]$ for $x \in V(S)$;
    and by $C_1, \dots, C_{(2\ell_0+2)N}$ the copies of $S$ in $G$, that is the subgraphs $G\big[\{(x,i) \mid x \in V(S)\}\big]$ for $i \in [(2\ell_0+2)N]$.
    By assumptions, $\td_k(G) = N$, and so by Observation~\ref{obs:equivalence_ktreedepth_kdecomposable_treedecomposition},
    $G$ admits a $k$-dismantable tree decomposition $\mathcal{D} = \big(T,(W_x \mid x \in V(T))\big)$ of width $N-1$.

    For every $x \in V(T)$, the bag $W_x$ has size at most $N$, and so intersects at most $N$ of the subgraphs $C_i$ for $i \in [(2\ell_0+2)N]$.
    Therefore, there is no set of $2\ell_0+2$ vertices of $T$ that intersect all
    the projections of $C_i$ for $i \in [(2\ell_0+2)N]$.
    By Lemma~\ref{lemma:helly_property} applied to the graph $T$,
    there are $2\ell_0+2$ pairwise disjoint such projections.
    More precisely, there are indices $1 \leq i_1< \dots< i_{2\ell_0+2} \leq (2\ell_0+2)N$ such that for every distinct $j,j' \in [2\ell_0+2]$,
    $T_j=T\big[\{x \in V(T) \mid W_x \cap C_{i_j} \neq \emptyset\}\big]$ and $T_{j'}=T\big[\{x \in V(T) \mid W_x \cap C_{i_{j'}} \neq \emptyset\}\big]$ are pairwise vertex-disjoint.
    For every $a \in [k]$ and $j \in [2\ell_0+1]$, let $R_{a,j}$ be the subpath of $R_a$ with first vertex in $C_{i_j}$ and last vertex in $C_{i_{j+1}}$.
    Moreover, let $u_{a,j}$ be the last vertex belonging to $\bigcup_{x \in V(T_j)} W_x \cap V(R_{a,j})$ along $R_{a,j}$.
    Note that for every $j \in [2\ell_0+1]$, $C'_j = \{u_{a,j} \mid a \in [k]\}$ is included in a bag of $\mathcal{D}$,
    namely $W_{z_j}$ where $z_j$ is the first node of the (unique) path from $V(T_j)$ to $V(T_{j+1})$ in $T$.
    Moreover, the vertices $u_{a,j}$ for $a \in [k]$ and $j \in [2\ell_0 + 1]$ are pairwise distinct.

    Let $F$ be the forest obtained from $T$ by removing every edge $xy \in E(T)$ such that $|W_x \cap W_y| < k$.
    We claim that $\{z_j \mid j \in [2\ell_0+1]\}$ is in a single connected component of $F$.
    Indeed, otherwise there is a cut of size less than $k$ between $W_{z_{j}}$ and $W_{z_{j'}}$ in $G$, for some distinct $j,j' \in [2\ell_0+1]$.
    But this is a contradiction since $R_1, \dots, R_k$ induce $k$ pairwise vertex-disjoint paths from $C'_j \subseteq W_{z_j}$ to $C'_{j'} \subseteq W_{z_{j'}}$.
    Let $T'$ be this connected component.

    Let $\hat{G}$ be the graph with vertex set $\bigcup_{z \in V(T')} W_z$ and edge set  $\bigcup_{z \in V(T')} \binom{W_z}{2}$.
    For every $a \in [k]$, let $r_{a,1}, \dots, r_{a, m_a}$ be the vertices of $V(R_a) \cap V(\hat{G})$, in this order along $R_a$.
    Observe that $\hat{R}_a = (r_{a,1}, \dots, r_{a, m_a})$ is a path in $\hat{G}$.
    Moreover, for every $a \in [k]$, $\hat{R}_a$ intersects each of $C'_1, \dots, C'_{2\ell_0+1}$ in this order along $\hat{R}_a$.
    Since $C'_j$ induces a clique in $\hat{G}$ for every $j \in [2\ell_0+1]$, 
    we deduce that there is a model of a $K_k \square P_{2\ell_0 +1}$ in $\hat{G}$.
    Moreover, $\hat{G}$ has a $k$-dismantable tree decomposition $\mathcal{D}' = \big(T', (W_x \mid x \in V(T'))\big)$ of width less than $N$ 
    and such that every adhesion has size at least $k$.
    It follows that there is a vertex $v \in V(\hat{G})$ such that 
    $v \in W_x$ for every $x \in V(T')$ and $\big(T', (W_x \setminus \{v\} \mid x \in V(T'))\big)$ is a $k$-dismantable tree decomposition.
    This implies that $\td_k(\hat{G}-v) < N$.
    But since $\hat{G}$ contains a model of $K_k \square P_{2\ell_0+1}$,
    $\hat{G} - v$ contains a model of $K_k \square P_{\ell_0}$.
    By Lemma~\ref{lemma:ktd_minor_monotone}, we deduce that 
    \[
        \td_k(S \square P_{\ell_0}) \leq \td_k(K_k \square P_{\ell_0}) \leq \td_k(\hat{G}-v) < N,
    \]
    contradicting the definition of $\ell_0$.
    This proves the lemma.
\end{proof}

\section{Finding a \texorpdfstring{$k$}{k}-ladder}\label{sec:finding_a_kladder}

In this section, we show how to extract a graph of the form $T \square P_\ell$ with a tree $T$ on $k$ vertices, 
from a graph having $k$ pairwise disjoint paths $Q_1, \dots, Q_k$, 
and sufficiently many pairwise vertex-disjoint connected subgraphs intersecting all paths $Q_i$.

\begin{figure}[ht]
    \centering
    \includegraphics{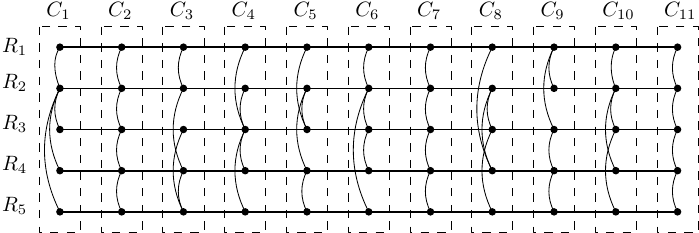}
    \caption{A $5$-ladder of length $10$.}
    \label{fig:def_kladder}
\end{figure}

Let $k,\ell$ be two positive integers.
A \emph{$k$-ladder of length $\ell-1$} is a graph $H$ with vertex set $[k] \times [\ell]$ such that
\begin{enumerate}
    \item for every $i \in [k]$, $R_i = H\big[\{(i,j) \mid j \in [\ell]\}\big]$ is a path with ordering $(i,1), \dots, (i,\ell)$, and
    \item for every $j \in [\ell]$, $C_j = H\big[\{(i,j) \mid i \in [k]\}\big]$ is connected.
\end{enumerate}
See Figure~\ref{fig:def_kladder}.
We call $R_1, \dots, R_k$ the \emph{rows} of $H$, and $C_1, \dots, C_\ell$ its \emph{columns}.

Since every $C_i$ is connected, we can assume without loss of generality that they are trees.
Then, using Cayley's formula stating that there are $k^{k-2}$ labelled trees on $k$ vertices, we deduce the following observation.
\begin{observation}\label{obs:TsquareP_from_a_long_ladder}
    For every $k,\ell$, if $H$ is a $k$-ladder of length $k^{k-2} (\ell-1)$,
    then there is a tree $T$ on $k$ vertices such that $T \square P_\ell$ is a minor of $H$.
\end{observation}
Hence it suffices to find a long enough $k$-ladder as a minor to find a minor of the form $T \square P_\ell$ for some tree $T$ on $k$ vertices.
Our main tool to do so is the following lemma, which is the main result of this section.

\begin{lemma}\label{lemma:finding_kladder_from_a_mess}    
    There is a function $f_{\cstref{lemma:finding_kladder_from_a_mess}}\colon \NN^2 \to \NN$ such that the following holds.
    Let $k$ be a positive integer, let $G$ be a graph,
    and let $Q_1,\dots, Q_k$ be $k$ pairwise vertex-disjoint paths in $G$.
    If there exists $f_{\cstref{lemma:finding_kladder_from_a_mess}}(k,\ell)$ pairwise vertex-disjoint connected subgraphs in $G$
    each intersecting $V(Q_i)$ for every $i \in [k]$,
    then $G$ has a $k$-ladder of length $\ell$ as a minor.
\end{lemma}

Together with Observation~\ref{obs:TsquareP_from_a_long_ladder}, this lemma implies Lemma~\ref{lemma:finding_kladder_from_a_mess_outline}.
The core of the proof of Lemma~\ref{lemma:finding_kladder_from_a_mess} is the following intermediary result.

\begin{lemma}\label{lemma:very_technical_lemma_to_find_a_kladder_from_a_mess}
    There is a function $f_{\cstref{lemma:very_technical_lemma_to_find_a_kladder_from_a_mess}}\colon \NN^2 \to \NN$ such that the following holds.
    Let $G$ be a graph, let $t, \ell$ be positive integers,
    let $P$ be a path in $G$, and let $A_1, \dots, A_{f_{\cstref{lemma:very_technical_lemma_to_find_a_kladder_from_a_mess}}(t,\ell)}$ 
    be pairwise vertex-disjoint connected subgraphs of $G$ such that
    for every $i \in [f_{\cstref{lemma:very_technical_lemma_to_find_a_kladder_from_a_mess}}(t,\ell)]$, $1 \leq |V(A_i) \cap V(P)|  \leq 2t - 1$.
    Then there are subgraphs $B_1, \dots, B_{\ell}$ in $\{A_i\}_{i \in [f_{\cstref{lemma:very_technical_lemma_to_find_a_kladder_from_a_mess}}(t,\ell)]}$,
    a subpath of $P$ of the form $P[\init(P),b]$, and
    for every $j \in [\ell]$ a subpath $P[a_j,b_j]$ of $P[\init(P),b]$,
    such that, for every $j \in [\ell]$,
    \begin{enumerate}[label=(\alph*)]
        \item $V(B_j) \cap V(P[\init(P),b]) \neq \emptyset$,
        \item $V(B_j) \cap V(P[\init(P),b]) \subseteq V(P[a_j,b_j])$, and
        \item $V(P[a_j,b_j]) \cap V(P[a_{j'},b_{j'}]) = \emptyset$ for every $j' \in [\ell]\setminus \{j\}$.
    \end{enumerate}
\end{lemma}

\begin{figure}
    \centering
    \includegraphics{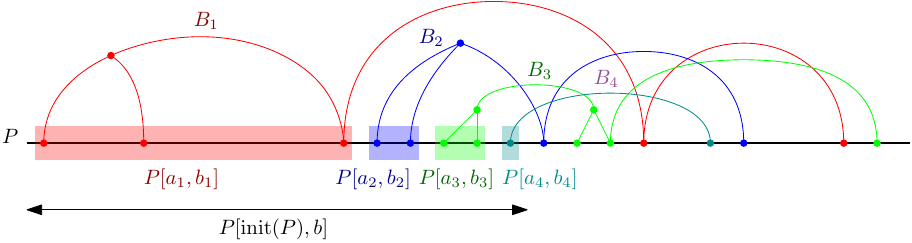}
    \caption{Illustration for Lemmas~\ref{lemma:very_technical_lemma_to_find_a_kladder_from_a_mess} and~\ref{lemma:technical_lemma_to_find_a_kladder_from_a_mess} with $\ell=4$. 
        There is a subpath $P[\init(P),b]$ such that each $B_j$ has a ``private'' interval $P[a_j,b_j]$ 
        that contains $V(B_j) \cap V(P[\init(P),b])$ and is disjoint from $V(P[a_{j'},b_{j'}])$ for every $j' \neq j$.}
    \label{fig:finding_a_k_ladder}
\end{figure}

Figure~\ref{fig:finding_a_k_ladder} gives an illustration for the outcome of this lemma.

\begin{proof}
    Let $f_{\cstref{lemma:very_technical_lemma_to_find_a_kladder_from_a_mess}}(t,\cdot)$ be the function defined inductively
    by $f_{\cstref{lemma:very_technical_lemma_to_find_a_kladder_from_a_mess}}(t,1)=1$ and 
    $f_{\cstref{lemma:very_technical_lemma_to_find_a_kladder_from_a_mess}}(t,\ell) = 1 + (2t-1)f_{\cstref{lemma:very_technical_lemma_to_find_a_kladder_from_a_mess}}(t,\ell-1)$ 
    for every $\ell \geq 2$.
    We proceed by induction on $\ell$.
    If $\ell=1$, then take for $B_1 = A_1$ and $P[\init(P),b] = P[a_1,b_1]=P$.

    Now suppose $\ell>1$.
    Without loss of generality $A_1$ is such that the first vertex along $P$ 
    in $\bigcup_{i \in [f_{\cstref{lemma:very_technical_lemma_to_find_a_kladder_from_a_mess}}(t,\ell)]} V(A_i)$ belongs to $A_1$.
    Let $x_1, \dots, x_r$ be the vertices in $V(A_1) \cap V(P)$, in this order along $P$.
    Note that $r \leq 2t - 1$ by assumptions.
    Let $x_{r+1}$ be the last vertex of $P$.
    For every $s \in [r]$, let $I_s$ be the set of indices $j \in \{2, \dots, f_{\cstref{lemma:very_technical_lemma_to_find_a_kladder_from_a_mess}}(t,\ell)\}$ 
    such that the first vertex in $V(A_j)$ along $P$
    belongs to $P[x_s, x_{s+1}]$.
    By the pigeonhole principle, there is an integer $s \in [r]$ such that 
    $|I_s| \geq \frac{f_{\cstref{lemma:very_technical_lemma_to_find_a_kladder_from_a_mess}}(t,\ell)-1}{2t-1} = f_{\cstref{lemma:very_technical_lemma_to_find_a_kladder_from_a_mess}}(t,\ell-1)$.

    Hence by the induction hypothesis applied to any family $\{A'_1, \dots, A'_{f_{\cstref{lemma:very_technical_lemma_to_find_a_kladder_from_a_mess}}(t,\ell-1)}\} \subseteq \{A_j\}_{j \in I_s}$ 
    and 
    \[
       P' = 
       \begin{cases} 
            P]x_s, x_{s+1}[ &\textrm{if $s<r$,} \\
            P]x_s, x_{s+1}] &\textrm{if $s=r$,}
       \end{cases}
    \]
    there are subgraphs $B_2, \dots, B_{\ell}$ in $\{A'_i\}_{i \in [f_{\cstref{lemma:very_technical_lemma_to_find_a_kladder_from_a_mess}}(t,\ell-1)]}$ and a subpath $P]x_s,b']$ of $P'$ such that
    for every $j \in \{2, \dots, \ell\}$, there is a subpath $P[a_j,b_j]$ of $P]x_s,b']$ such that
    \begin{enumerate}[label=(\alph*')]
        \item $V(B_j) \cap V(P[\init(P),b]) \neq \emptyset$, and
        \item $V(B_j) \cap V(P[\init(P),b]) \subseteq V(P[a_j,b_j])$, and
        \item $V(P[a_{j'},b_{j'}]) \cap V(P[a_{j''},b_{j''}]) = \emptyset$ for every distinct $j',j'' \in \{2, \dots, \ell\} \setminus \{j\}$.
    \end{enumerate}
    Then for $B_1 = A_1$, $P[\init(P),b] = P[\init(P), b']$ and $P[a_1,b_1] = P[\init(P),x_s]$, this proves the lemma.
\end{proof}

The following simple lemma will be useful.
\begin{lemma}\label{lemma:tree_kleaves_or_kpath}
    For every positive integer $k$, and for every tree $T$ on at least $(k-1)(k-2)+2$ vertices, at least one of the following holds
    \begin{enumerate}[label=(\arabic*)]
        \item $T$ has at least $k$ leaves, or
        \item $T$ has a path on $k$ vertices.
    \end{enumerate}
\end{lemma}

\begin{proof}
    If $k=1$ then the result is clear. Now assume $k \geq 2$.
    Root $T$ on an arbitrary vertex $r$, and suppose that $T$ has at most $k-1$ leaves.
    Then for every leaf $x$, consider the path $T[x,r[$.
    Since $\bigcup_{x \text{ leaf}} V(T[x,r[)$ covers $V(T-r)$, by the pigeon hole principle,
    there is a leaf $x$ such that $|V(T[x,r[)| \geq k-1$ and so $|V(T[x,r])| \geq k$.
\end{proof}

We now show that the condition $|V(A_i) \cap V(P)|  \leq 2t - 1$ in Lemma~\ref{lemma:very_technical_lemma_to_find_a_kladder_from_a_mess} can be removed,
assuming that $G$ has bounded treewidth.

\begin{lemma}\label{lemma:technical_lemma_to_find_a_kladder_from_a_mess}
    There is a function $f_{\cstref{lemma:technical_lemma_to_find_a_kladder_from_a_mess}}\colon\NN^2 \to \NN$ such that the following holds.
    Let $t, \ell$ be positive integers, and
    let $G$ be a graph with $\tw(G) < t$.
    If $P$ is a path in $G$, and if $\{A_1, \dots, A_{f_{\cstref{lemma:technical_lemma_to_find_a_kladder_from_a_mess}}(t,\ell)}\}$
    is a family of pairwise vertex-disjoint connected subgraphs of $G$, all intersecting $V(P)$,
    then there exists $\ell$ subgraphs $B_{1}, \dots, B_{\ell}$ among $\{A_1, \dots, A_{f_{\cstref{lemma:technical_lemma_to_find_a_kladder_from_a_mess}}(t,\ell)}\}$
    and a subpath of $P$ of the form $P[\init(P),b]$
    such that for every $j \in [\ell]$, there is a subpath $P[a_j,b_j]$ of $P[\init(P),b]$ such that
    \begin{enumerate}[label=(\alph*)]
        \item $V(B_{j}) \cap V(P[\init(P),b]) \subseteq V(P[a_j,b_j])$, and
        \item $V(P[a_j,b_j]) \cap V(P[a_{j'},b_{j'}]) = \emptyset$ for every $j' \in [\ell] \setminus \{j\}$.
    \end{enumerate}
\end{lemma}

\begin{proof}
    Let $f_{\cstref{lemma:technical_lemma_to_find_a_kladder_from_a_mess}}(t,\ell) = t \max\{2, f_{\cstref{lemma:very_technical_lemma_to_find_a_kladder_from_a_mess}}(t,\ell)\}^2$
    where $f_{\cstref{lemma:very_technical_lemma_to_find_a_kladder_from_a_mess}}$ is as in Lemma~\ref{lemma:very_technical_lemma_to_find_a_kladder_from_a_mess}.
    We can assume that if $u,v \in V(P) \cap V(A_i)$ are distinct and in this order along $P$, 
    then $P[u,v]$ intersect $A_j$ for some $j \neq i$, for every $i \in [f_{\cstref{lemma:technical_lemma_to_find_a_kladder_from_a_mess}}(t,\ell)]$. Otherwise we just contract $P[u,v]$ into a single vertex.

    Consider now a tree decomposition $\big(T,(W_x \mid x \in V(T))\big)$ of $G$ of width less than $t$.
    Since $A_1, \dots, A_{f_{\cstref{lemma:technical_lemma_to_find_a_kladder_from_a_mess}}(t,\ell)}$ are pairwise vertex-disjoint, 
    every set of vertices intersecting all of them must have size at least $f_{\cstref{lemma:technical_lemma_to_find_a_kladder_from_a_mess}}(t,\ell)$.
    Consider now for every $i\in [f_{\cstref{lemma:technical_lemma_to_find_a_kladder_from_a_mess}}(t,\ell)]$ the projection $T_i$ of $A_i$ on $T$.
    By Lemma~\ref{lemma:helly_property}, the family $(T_i)_{i \in [f_{\cstref{lemma:technical_lemma_to_find_a_kladder_from_a_mess}}(t,\ell)]}$ contains $N = \frac{f_{\cstref{lemma:technical_lemma_to_find_a_kladder_from_a_mess}}(t,\ell)}{t} = \max\{2, f_{\cstref{lemma:very_technical_lemma_to_find_a_kladder_from_a_mess}}(t,\ell)\}^2$ 
    pairwise vertex-disjoint
    subtrees $T_{i_1}, \dots, T_{i_N}$.

    Consider now the tree $T'$ obtained from $T$ by contracting every vertex outside $\bigcup_{j \in [N]} V(T_{i_j})$ 
    with one of its closest vertex in $\bigcup_{j \in [n]} V(T_{i_j})$,
    and by contracting $T_{i_j}$ into a single node $x_j$ for each $j \in [N]$.
    Hence $V(T') = \{x_j \mid j \in [N]\}$.
    By Lemma~\ref{lemma:tree_kleaves_or_kpath}, $T'$ has either $\sqrt{N} \geq f_{\cstref{lemma:very_technical_lemma_to_find_a_kladder_from_a_mess}}(t,\ell)$ leaves or
    a path with $\sqrt{N} \geq f_{\cstref{lemma:very_technical_lemma_to_find_a_kladder_from_a_mess}}(t,\ell)$ vertices. 
    Let $I$ be a subset of $[N]$ of size $\sqrt{N}$ such that $x_j$ for $j \in I$ is a leaf in $T'$, or such that a permutation of $(x_j)_{j \in I}$ is a path in $T'$.
    We claim that
    \begin{equation}\label{eq:lemma_technical_to_find_kladder_nice_Hj}
        \parbox{0.85\textwidth}{for every $j \in I$, the number of pairs $u,v \in V(P)$ in this order along $P$ such that $u \in V(A_{i_j})$, $V(P]u,v[) \cap \left(\bigcup_{j' \in I} V(A_{i_{j'}})\right) = \emptyset$, and $v \in \bigcup_{j' \in I \setminus \{j\}} V(A_{i_{j'}})$, is at most $2(t-1)$.}
    \end{equation}   
    Indeed, if $x_j$ is a leaf, then every time $P$ leaves $A_{i_j}$, $P$ crosses the adhesion between the root of $T_{i_j}$ 
    and its unique neighbor out of $V(T_{i_j})$, which has size at most $t-1$.
    And if a permutation of $(x_{j'})_{j' \in I}$ is a path in $T'$, then every time $P$ leaves $A_{i_j}$,
    $P$ must cross one of the at most two adhesions neighboring $T_{i_j}$ and leading to another $T_{i_{j'}}$ for some $j' \in I \setminus \{j\}$.
    Since these adhesions have size at most $t-1$, this proves \eqref{eq:lemma_technical_to_find_kladder_nice_Hj}.

    \medskip

    Then contract repeatedly every edge of $P$ that intersect at most one of the sets $V(A_{i_j})$ for $j \in I$.
    Let $G'$ be the resulting minor of $G$, let $P'$ be the path in $G'$ induced by $P$, 
    and for every $j \in I$, let $A'_j$ be the connected subgraph of $G'$ inherited from $A_{i_j}$.
    Then $|V(P') \cap A'_j| \leq 2t-1$ for every $j \in I$.
    Hence we can apply Lemma~\ref{lemma:very_technical_lemma_to_find_a_kladder_from_a_mess} to 
    $G', P', \{A_{i_j} \mid j \in I\}$.
    This is possible since $|I| = \sqrt{N} \geq f_{\cstref{lemma:very_technical_lemma_to_find_a_kladder_from_a_mess}}(t,\ell)$.
    We obtain subgraphs $B'_1, \dots, B'_\ell$ of $G'$ and a subpath $P'[\init(P'),b']$ of $P'$ corresponding respectively
    to subgraphs $B_1, \dots, B_\ell$ of $G$ 
    and a subpath $P[\init(P),b]$ of $P$ 
    as desired.
    This proves the lemma.
\end{proof}

To deduce Lemma~\ref{lemma:finding_kladder_from_a_mess}, we will use the celebrated Erd\H{o}s-Szekeres Theorem.
\begin{theorem}[Erd\H{o}s Szekeres~\cite{ErdosSzekeres}]\label{thm:erdos-szkeres}
    For all integers $r,s$, for every permutation $\sigma$ of $[(r-1)(s-1)+1]$,
    either
    there is a set $I \subseteq [(r-1)(s-1)+1]$ of size $r$ such that $\sigma\vert_I$ is increasing, or        
    there is a set $D \subseteq [(r-1)(s-1)+1]$ of size $s$ such that $\sigma\vert_D$ is decreasing.
\end{theorem}

\begin{proof}[Proof of Lemma~\ref{lemma:finding_kladder_from_a_mess}]
    Let $t= \max\{2,f_{\cstref{thm:grid_minor_thm}}(\max\{k,\ell+1\})\}$.
    If $\tw(G) \geq t$, then the $k \times (\ell+1)$ grid is a minor of $G$
    by Theorem~\ref{thm:grid_minor_thm}, and so $G$ contains a $k$-ladder of length $\ell$ as a minor.
    Now assume that $\tw(G)<t$.
    For every $i\in\{0,\dots, k\}$, let
    \[
        f_{i,k}(\ell) = 
        \left\{
        \begin{aligned}
            &(\ell+1)^{2^{k-1}}                                                                       & \text{if $i=0$,} \\
            &f_{\cstref{lemma:technical_lemma_to_find_a_kladder_from_a_mess}}(t, f_{i-1,k}(\ell)) & \text{if $i\geq 1$,} \\
        \end{aligned}
        \right.
    \]
    where $f_{\cstref{lemma:technical_lemma_to_find_a_kladder_from_a_mess}}$ is as in Lemma~\ref{lemma:technical_lemma_to_find_a_kladder_from_a_mess},
    and take $f_{\cstref{lemma:finding_kladder_from_a_mess}}(k,\ell) = f_{k,k}(\ell)$.

    Let $G$ be a graph, 
    let $Q_1,\dots, Q_k$ in $G$ be $k$ pairwise disjoint paths in $G$,
    and let $A_1, \dots, A_{f_{\cstref{lemma:finding_kladder_from_a_mess}}(k,\ell)}$
    be pairwise disjoint connected subgraphs in $G$ each intersecting $V(Q_i)$ for every $i \in [k]$,
    By applying Lemma~\ref{lemma:technical_lemma_to_find_a_kladder_from_a_mess} successively on the paths $Q_1, \dots, Q_k$ and contracting the subpaths of the form $Q_i[a_j,b_j]$,
    we obtain a minor $G'$ of $G$, connected subgraphs $B'_1, \dots, B'_{(\ell+1)^{2^{k-1}}}$ of $G'$,
    and paths $Q'_1, \dots, Q'_k$ in $G'$ such that
    \[
      |V(B_j) \cap V(Q'_a)| = 1
    \]
    for every $j \in \big[(\ell+1)^{2^{k-1}}\big]$ and $a \in [k]$.
    For every $j \in \big[(\ell+1)^{2^{k-1}}\big]$ and $a \in [k]$, let $u_{j,a}$ be the unique vertex in $V(B_j) \cap V(Q'_a)$.
    Finally, by applying Theorem~\ref{thm:erdos-szkeres} (Erd\H{o}s-Szekeres Theorem) $k-1$ times, 
    there are $\ell+1$ distinct indices $j_1, \dots, j_{\ell+1}$ in $\big\{1, \dots, (\ell+1)^{2^{k-1}}\big\}$ such that the orderings of 
    $(u_{j_i,a})_{i \in [\ell+1]}$ along $Q'_a$ for $a \in [k]$ are pairwise identical or reverse of each other.
    By possibly replacing $Q'_{a}$ by its reverse for every $a \in [k]$, and by contracting some edges, we obtain a $k$-ladder of length $\ell$ as a minor of $G'$
    with rows inherited from $Q'_1, \dots, Q'_k$ and columns inherited from $B_{j_1}, \dots, B_{j_{\ell+1}}$.
    This proves the lemma.
\end{proof}

\section{Nice pairs}\label{sec:nice_sets}

In order to show that every graph $G$ excluding all the $k$-ladders of length $\ell$ as minors have bounded $k$-treedepth, we need to find a 
decomposition of $G$ into graphs $G_1, \dots, G_m$, such that $G$ is a $(<k)$-clique-sum of $G_1, \dots, G_m$.
Typically, $G_i$ is the ``torso'' of a bag in a suitable tree decomposition.
The main issue with this approach is that we then need to decompose each $G_i$, while $G_i$ is not necessarily a minor of $G$ (since we possibly added some cliques).
The solution we develop in this section is to force $G_i$ to be ``nice'' in $G$, which roughly means that the cliques added to $G[V(G_i)]$ to obtain $G_i$
do not increase the connectivity between any two subsets of $V(G_i)$.
In this section, we define this notion of ``nice'' sets, and we prove several properties on them.
In the following section, we will find tree decompositions that will enable us to decompose our graph into nice sets in the final proof.

\begin{figure}[ht]
    \centering
    \includegraphics{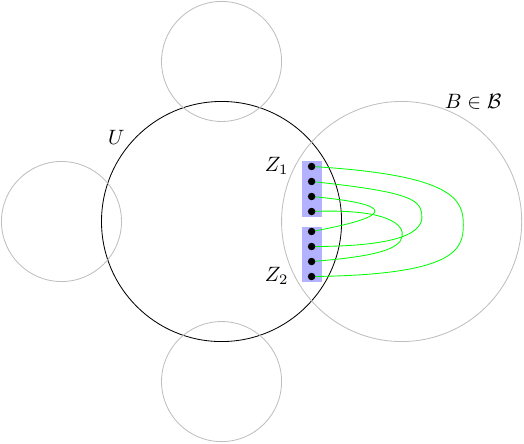}
    \caption{A good pair $(U,\mathcal{B})$ in $G$ is nice if for every $B \in \mathcal{B}$, for every $i \geq 0$,
        for every $Z_1,Z_2 \subseteq U \cap B$ both of size $i$, 
        there are $i$ pairwise disjoint $(Z_1,Z_2)$-paths in $G[B] \setminus \binom{U \cap B}{2}$.
        This property implies that the maximum number of pairwise disjoint paths between two given subsets of $U$ is the same in
        $G$ and in $\torso_G(U,\mathcal{B})$. See Lemma~\ref{lemma:finding_k_paths_between_sets_connected_in_the_torso_of_a_nice_set}.}
    \label{fig:nice_pair}
\end{figure}

Let $G$ be a graph.
A \emph{good pair} in $G$ is a pair $(U, \mathcal{B})$
where $U$ is a non-empty subset of $V(G)$, and
$\mathcal{B}$ is a family of subsets of $V(G)$ such that 
\begin{enumerate}[label=(g\arabic*)]
    \item $G = \bigcup_{B \in \mathcal{B} \cup \{U\}} G[B]$, and \label{item:def_nice_pair:i}
    \item the sets $B \setminus U$ for $B \in \mathcal{B}$ are pairwise disjoint. \label{item:def_nice_pair:ii}
\end{enumerate}
Equivalently, $\mathcal{B} \cup \{U\}$ is the family of the bags of a tree decomposition indexed by a star whose center bag is $U$.
Then, we say that $(U,\mathcal{B})$ is a \emph{nice pair} in $G$ if it is a good pair in $G$ and
\begin{enumerate}[resume*]
    \item for every $B \in \mathcal{B}$, for every $i \geq 0$, \label{item:def_nice_pair:iii}
        for every $Z_1,Z_2 \subseteq U \cap B$ both of size $i$, 
        there are $i$ pairwise disjoint $(Z_1,Z_2)$-paths in $G[B] \setminus \binom{B \cap U}{2}$.
\end{enumerate}
See Figure~\ref{fig:nice_pair}.

\medskip

In practice, we will build good pairs from tree decompositions as follows.
If $\big(T,(W_x \mid x \in V(T))\big)$ is a tree decomposition of a graph $G$,
and if $R$ is a subtree of $T$, then $(U, \mathcal{B})$ 
for $U = \bigcup_{x \in V(R)} W_x$ and $\mathcal{B} = \{\bigcup_{x \in V(S)} W_x \mid \text{$S$ connected component of $T-V(R)$}\}$ is a good pair in $G$.
To ensure that $(U,\mathcal{B})$ is a nice pair, we will need further assumptions on the tree decomposition.
Finding such decompositions will be the goal of Section~\ref{sec:finding_a_suitable_tree_decomposition},
but for now, we prove some properties on nice pairs.

First we give some notation.
Let $G$ be a graph and let $(U,\mathcal{B})$ be a good pair in $G$.
The \emph{torso of $U$ in $G$ with respect to $\mathcal{B}$},
denoted by $\torso_G(U,\mathcal{B})$, is the graph with vertex set $U$ and edge set $E(G[U]) \cup \bigcup_{B \in \mathcal{B}} \binom{U \cap B}{2}$.

\begin{lemma}\label{lemma:U_nice_impliesU-X_nice}
    Let $G$ be a graph and let $(U, \mathcal{B})$ be a nice pair in $G$.
    Then for every $X \subseteq U$, 
    there exists a nice pair in $G-X$ of the form $(U\setminus X, \mathcal{B}')$ with $\torso_{G-X}(U\setminus X, \mathcal{B}') = \torso_G(U,\mathcal{B}) - X$.
\end{lemma}

\begin{proof}
    We claim that $(U \setminus X, \{B \setminus X \mid B \in \mathcal{B}\})$ is nice in $G-X$.
    First, \ref{item:def_nice_pair:i} and \ref{item:def_nice_pair:ii} of the definition of good pair clearly hold,
    and so $(U \setminus X, \{B \setminus X \mid B \in \mathcal{B}\})$ is a good pair.
    
    Consider now $B \in \mathcal{B}$,
    a positive integer $i$,
    and two sets $Z_1,Z_2 \subseteq (U \setminus X) \cap (B \setminus X)$ both of size $i$.
    Let $Z_0 = B \cap U \cap X$.
    Since $(U, \mathcal{B})$ is a nice pair,
    there are $i + |Z_0|$ pairwise disjoint $(Z_0 \cup Z_1, Z_0 \cup Z_2)$-paths in $G[B] \setminus \binom{B \cap U}{2}$.
    Note that $B \cap U = ((U \setminus X) \cap (B \setminus X)) \cup Z_0$.
    Hence there are $i$ pairwise disjoint $(Z_1,Z_2)$-paths in $G[B \setminus X] \setminus \binom{(U \setminus X) \cap (B \setminus X)}{2}$.
    This proves that $\big(U \setminus X, \{B \setminus X \mid B \in \mathcal{B}\}\big)$ is a nice pair in $G-X$.
    Finally, by definition of $\torso$,
    $\torso_G\big(U \setminus X, \{B \setminus X \mid B \in \mathcal{B}\}\big) = \torso_G(U,\mathcal{B}) - X$.
\end{proof}

\begin{lemma}\label{lemma:technical_about_nice_sets}
    Let $G$ be a graph and let $(U, \mathcal{B})$ be a nice pair in $G$.
    For every $B \in \mathcal{B}$, $(V(G) \setminus (B \setminus U), \{B\})$ is a nice pair in $G$.
\end{lemma}

\begin{proof}
    Clearly, $(V(G) \setminus (B\setminus U), \{B\})$ is a good pair in $G$.
    Moreover, for every $i \geq 0$,
    for every $Z_1,Z_2 \subseteq U \cap B$ both of size $i$,
    there are $i$ pairwise disjoint $(Z_1,Z_2)$-paths in $G[B] \setminus \binom{B \cap U}{2}$.
    Hence $(V(G) \setminus (B\setminus U), \{B\})$ is a nice pair in $G$.
\end{proof}

\begin{lemma}\label{lemma:paths_in_G_implies_paths_in_torso}
    Let $k'$ be a positive integer.
    Let $G$ be a graph and let $(U, \mathcal{B})$ be a good pair in $G$.
    For every $Z_1,Z_2 \subseteq U$ with $|Z_1|=|Z_2|=k'$ such that there are $k'$ pairwise disjoint $(Z_1,Z_2)$-paths $Q_1, \dots, Q_{k'}$ in $G$,
    there are $k'$ pairwise disjoint $(Z_1,Z_2)$-paths $\hat{Q}_1, \dots, \hat{Q}_{k'}$ in $\torso_G(U,\mathcal{B})$. Moreover $V(Q_i) \cap U = V(\hat{Q}_i)$ for every $i \in [k']$.
\end{lemma}

\begin{proof}
    Let $Z_1,Z_2 \subseteq U$ with $|Z_1|=|Z_2|=k'$  and let $Q'_1, \dots, Q'_{k'}$ be $k'$ pairwise disjoint $(Z_1,Z_2)$-paths in $G$.
    Let $a \in [k']$ and let $q_{a,1}, \dots, q_{a,\ell_a}$ be the vertices in $V(Q_a) \cap U$ in this order along $Q_a$.
    Then, for every $i \in [\ell_a - 1]$, either $q_{a,i}q_{a,i+1}$ is an edge of $G$, or there exists $B \in \mathcal{B}$ such that $q_{a,i}q_{a,i+1} \in B$.
    In both cases $q_{a,i}q_{a,i+1}$ is an edge in $\torso_G(U,\mathcal{B})$.
    Hence $\hat{Q}_a = (q_{a,1}, \dots, q_{a,\ell_a})$ is a path in $\torso_G(U,\mathcal{B})$ with $V(\hat{Q}_a) = V(Q_a) \cap U$,
    and so $(\hat{Q}_1, \dots, \hat{Q}_{k'})$ is as claimed.
\end{proof}

The following lemma, which is the key property of nice pairs, informally says that if $(U,\mathcal{B})$ is nice, then the converse of the previous lemma also holds.

\begin{lemma}\label{lemma:finding_k_paths_between_sets_connected_in_the_torso_of_a_nice_set}
    Let $k'$ be a positive integer.
    Let $G$ be a graph and let $(U, \mathcal{B})$ be a nice pair in $G$.
    For every $Z_1,Z_2 \subseteq U$ with $|Z_1|=|Z_2|=k'$ such that there are $k'$ pairwise disjoint $(Z_1,Z_2)$-paths $Q_1, \dots, Q_{k'}$ in $\torso_G(U, \mathcal{B})$,
    there are $k'$ pairwise disjoint $(Z_1,Z_2)$-paths $Q'_1, \dots, Q'_{k'}$ in $G$. Moreover $\bigcup_{i \in [k']} V(Q'_i) \cap U \subseteq \bigcup_{i \in [k']} V(Q_i)$.
\end{lemma}

\begin{proof}
    We proceed by induction on $|\mathcal{B}|$.
    If $\mathcal{B} = \emptyset$, then $\torso_G(U,\mathcal{B}) = G$ and the result is clear.
    Now assume that $\mathcal{B} \neq \emptyset$ and consider some $B \in \mathcal{B}$.
    Note that $(U, \mathcal{B} \setminus \{B\})$ is a nice pair in $G' = \torso_G(V(G) \setminus (B \setminus U), \{B\})$.
    Hence by the induction hypothesis applied to $(U, \mathcal{B} \setminus \{B\})$ in $G'$, there are $k'$ pairwise disjoint $(Z_1,Z_2)$-paths $Q''_1, \dots, Q''_{k'}$ in $G'$
    with $\bigcup_{i \in [k']} V(Q''_i) \cap U \subseteq \bigcup_{i \in [k']} V(Q_i)$.
    Without loss of generality, assume that among $Q''_1, \dots, Q''_{k'}$, the paths intersecting $U \cap B$ in at least one vertex are $Q''_1, \dots, Q''_r$.
    Moreover, by taking $Q''_1, \dots, Q''_{k'}$ with minimum total length,
    we assume that these paths are induced paths in $G'$.
    In particular, since $B \cap U$ induces a clique in $G'$,
    we have $V(Q''_i) \cap B = \{z^i_1, z^i_2\}$ for some possibly equal $z_1^i, z_2^i \in U \cap B$,
    for every $i \in [r]$.
    Let $Z'_1 = \{z_1^i \mid i \in [r]\}$, $Z'_2 = \{z_2^i \mid i \in [r]\}$, and $Z'_0 = (U \cap B) \setminus (Z'_1 \cup Z'_2)$.
    By hypothesis, there are $r + |Z'_0|$ pairwise disjoint $(Z'_0 \cup Z'_1, Z'_0 \cup Z'_2)$-paths $Q'_1, \dots, Q'_r, (\{u\})_{u \in Z'_0}$ in $G[B] \setminus \binom{B \cap U}{2}$.
    Hence the union of $Q'_1, \dots, Q'_r$ with $Q''_1 \setminus z_1^1 z_2^1, \dots, Q''_r \setminus z_1^r z_2^r$, together with $Q''_{r+1}, \dots, Q''_{k'}$, 
    yields $k'$ pairwise vertex-disjoint $(Z_1,Z_2)$-paths in $G$.
    Moreover $\bigcup_{i \in [k']} V(Q'_i) \cap U \subseteq \bigcup_{i \in [k']} V(Q_i)$. This proves the lemma.
\end{proof}

The following lemma informally says that a nice pair in the torso of a nice pair in $G$, is itself a nice pair in $G$.

\begin{lemma}\label{lemma:nice_in_torso_of_nice_are_nice}
    Let $G$ be a graph and let $(U,\mathcal{B})$ be a nice pair in $G$.
    If $(U',\mathcal{B}')$ is a nice pair in $\torso_G(U,\mathcal{B})$,
    then there exists $\mathcal{B}''$ such that $(U', \mathcal{B}'')$ is a nice pair in $G$,
    and $\torso_G(U', \mathcal{B}'') = \torso_{\torso_G(U,\mathcal{B})}(U', \mathcal{B}')$.
\end{lemma}

\begin{proof}
    We proceed by induction on $|\mathcal{B}|$.
    If $\mathcal{B} = \emptyset$, then $\torso_{G}(U, \mathcal{B}) = G$ and so $(U', \mathcal{B}')$ is a nice pair in $G$ with $\torso_G(U', \mathcal{B}') = \torso_{\torso_G(U,\mathcal{B})}(U', \mathcal{B}')$.
    Now suppose that $\mathcal{B} \neq \emptyset$.

    Consider a member $B$ of $\mathcal{B}$.
    Let $G' = \torso_{G}(V(G) \setminus (B \setminus U), \{B\})$.
    Observe that $(U, \mathcal{B} \setminus \{B\})$ is a nice pair in $G'$.
    Hence, by the induction hypothesis,
    there is a nice pair in $G'$ of the form $(U', \mathcal{B}''_0)$
    with $\torso_{G'}(U', \mathcal{B}''_0) = \torso_{\torso_{G'}(U,\mathcal{B}\setminus \{B\})}(U', \mathcal{B}')$.
    By definition of $G'$,
    the set $B \cap U$ induces a clique in $G'$.
    As $\{U'\} \cup \mathcal{B}''_0$ is the family of the bags of a tree decomposition of $G'$,
    this implies that there exists $B''_0 \in \mathcal{B}''_0 \cup \{U'\}$ such that $B \cap U \subseteq B''_0$.

    \medskip

    First suppose that $B''_0 \neq U'$.
    Let 
    \begin{align*}
        B''_1 &= B \cup B''_0 \\
    \intertext{and}
        \mathcal{B}'' &= (\mathcal{B}''_0 \setminus \{B''_0\}) \cup \{B''_1\}.
    \end{align*}
    By construction, 
    \[
        \torso_G(U', \mathcal{B}'') 
        = \torso_{G'}(U', \mathcal{B}''_0) 
        = \torso_{\torso_{G'}(U,\mathcal{B}\setminus \{B\})}(U',\mathcal{B}')
        = \torso_{\torso_G(U,\mathcal{B})}(U', \mathcal{B}').
    \]
    It remains to show that $(U', \mathcal{B}'')$ is a nice pair in $G$.
    Observe that \ref{item:def_nice_pair:i} holds since $G' = \bigcup_{B' \in \mathcal{B}''_0 \cup \{U'\}} G'[B']$, $E(G) = (E(G') \setminus E(G'[B''_0])) \cup E(G[B''_1])$.
    Moreover, \ref{item:def_nice_pair:ii} holds since, for all $B' \in \mathcal{B}''_0 \setminus \{B''_0\}$, 
    the set $B' \setminus U$ is disjoint from $B$, and so from $B''_1$.
    We now prove \ref{item:def_nice_pair:iii}.
    Let $B'' \in \mathcal{B}''$, $i \geq 0$,
    and $Z_1, Z_2 \subseteq B'' \cap U'$ both of size $i$.
    By possibly replacing $Z_a$ by $Z_a\cup (B'' \cap U' \setminus (Z_1 \cup Z_2))$ for each $a \in \{1,2\}$ and increasing $i$ accordingly,
    we now assume that $Z_1 \cup Z_2 = B'' \cap U'$. 
    If $B'' \neq B''_1$, then there are $i$ pairwise disjoint $(Z_1 \cap B'', Z_2 \cap B'')$-paths 
    in $G'[B''] \setminus \binom{B'' \cap U'}{2} = G[B''] \setminus \binom{B'' \cap U'}{2}$ since
    $(U', \mathcal{B}''_0)$ is a nice pair in $G'$.
    If $B'' = B''_1$, then there are $|Z_1 \cap B''_0|$ pairwise disjoint $(Z_1,Z_2)$-paths $Q_1, \dots, Q_i$ in $G'[B''_0] \setminus \binom{B''_0 \cap U'}{2}$ since
    $(U', \mathcal{B}''_0)$ is a nice pair in $G'$.
    Without loss of generality, suppose that among $Q_1, \dots, Q_i$, the paths $Q_1, \dots, Q_j,$ are the ones intersecting $B$.
    For each $a \in [j]$, 
    let $z^a_1,z^a_2$ be respectively the first and last vertex of $Q_a$ in $U \cap B$
    (with possibly $z_a^1=z^a_2$).
    Since $(U, \mathcal{B})$ is a nice pair in $G$, there are $j$ pairwise disjoint $(\{z^a_1 \mid a \in [j]\}, \{z^a_2 \mid a \in [j]\})$-paths 
    in $G[B] \setminus \binom{U \cap B}{2}$.
    Combining these paths with $Q_1, \dots, Q_i$,
    we obtain $i$ pairwise disjoint $(Z_1,Z_2)$-paths in 
    $G[B \cup B''_0] \setminus \binom{(B \cup B''_0) \cap U'}{2}
    = G[B''_1] \setminus \binom{U' \cap B''_0}{2}$.

    \medskip

    Now assume $B''_0 = U'$. Let 
    \[
        \mathcal{B}'' = \mathcal{B}''_0 \cup \{B\}.
    \]
    By the definition of $G'$, 
    \[
        \torso_G(U', \mathcal{B}'') 
        = \torso_{G'}(U', \mathcal{B}''_0) 
        = \torso_{\torso_{G'}(U,\mathcal{B}\setminus \{B\})}(U',\mathcal{B}')
        = \torso_{\torso_G(U,\mathcal{B})}(U', \mathcal{B}').
    \]
    It remains to show that $(U', \mathcal{B}'')$ is a nice pair in $G$.
    Observe that \ref{item:def_nice_pair:i} holds since $G \subseteq G' \cup G[B]$ and $G' = \bigcup_{B' \in \mathcal{B}''_0 \cup \{U'\}} G'[B']$.
    Moreover, \ref{item:def_nice_pair:ii} holds since $B$ is disjoint from all $B'' \setminus U$ for $B'' \in \mathcal{B}''_0$.
    We now prove \ref{item:def_nice_pair:iii}.
    Let $B'' \in \mathcal{B}''$, $i \geq 0$,
    and $Z_1, Z_2 \subseteq B'' \cap U'$ both of size $i$.
    If $B'' \neq B$, then there are $i$ pairwise disjoint 
    $(Z_1 \cap B'', Z_2 \cap B'')$-paths in 
    $G'[B''] \setminus \binom{B'' \cap U'}{2} 
    = G[B''] \setminus \binom{B'' \cap U'}{2}$ since
    $(U', \mathcal{B}''_0)$ is a nice pair in $G'$.
    If $B'' = B$, then there are $i$ pairwise disjoint 
    $(Z_1,Z_2)$-paths $Q_1, \dots, Q_i$ in 
    $G[B] \setminus \binom{B \cap U}{2} 
    = G[B] \setminus \binom{B \cap U'}{2}$ since
    $(U, \mathcal{B})$ is a nice pair in $G$.
    Hence \ref{item:def_nice_pair:iii} holds.
    This concludes the proof of the lemma.
\end{proof}

\section{Finding a suitable tree decomposition}\label{sec:finding_a_suitable_tree_decomposition}

In this section, we prove the existence of ``well connected'' tree decompositions,
that will be used to decompose a graph into nice pairs in the final proof of Theorem~\ref{thm:main}.
This is inspired by the seminal result of Thomas~\cite{Thomas1990}, which asserts that every graph $G$
admits a tree decomposition $\big(T,(W_x \mid x \in V(T))\big)$ of width $\tw(G)$ which is \textit{lean}, 
that is such that for every $x_1,x_2 \in V(T)$, for every $Z_1\subseteq W_{x_1}$ and $Z_2 \subseteq W_{x_2}$ of same size $i$,
either there are $i$ pairwise disjoint $(Z_1,Z_2)$-paths in $G$, or there exists $z_1z_2 \in E(T[x_1,x_2])$ such that $|W_{z_1} \cap W_{z_2}|<i$.
Hence the connectivity between two bags is witnessed by the smallest adhesion between them.
Here we prove analogous results for tree decompositions of $(G,S)$, with some additional properties.
The techniques we use are strongly inspired by the short proof of Thomas' theorem by Bellenbaum and Diestel~\cite{BELLENBAUM2002}.

First we define a potential function on tree decompositions.
Let $k$ be a positive integer,
let $G$ be a graph, and
let $S \subseteq V(G)$.
Let $\mathcal{D} = \big(T,(W_x \mid x \in V(T))\big)$ be a tree decomposition of $(G,S)$.
For every positive integer $i$,
let $s_i(\mathcal{D}) = |\{x \in V(T) \mid |W_x| = i\}|$, that is the number of bags of size $i$ in $\mathcal{D}$, and let
\[
s^{(k)}_{i}(\mathcal{D}) =
\begin{cases}
    |\{x \in V(T) \mid |W_x| = i\}| & \textrm{if $i > \frac{3}{2}(k-1)$,} \\
    |\{x \in V(T) \mid |W_x| = i \text{ and } \exists y \in N_T(x), |W_x \cap W_y| \geq k\}| & \textrm{if $i \leq \frac{3}{2}(k-1)$.}
\end{cases}
\]
Then, let 
\begin{align*}
    s(\mathcal{D}) &= \big(s_{|V(G)|}(\mathcal{D}), \dots, s_0(\mathcal{D})\big) \\
\intertext{and}
    s^{(k)}(\mathcal{D}) &= \big(s^{(k)}_{|V(G)|}(\mathcal{D}), \dots, s^{(k)}_0(\mathcal{D})\big).
\end{align*}
We will consider these tuples in the lexicographic order.
The following lemma shows how to improve $\mathcal{D}$ by reducing $s(\mathcal{D})$ or $s^{(k)}(\mathcal{D})$ 
if some connectivity property is not satisfied.
Its proof follows step by step Bellenbaum and Diestel's argument~\cite{BELLENBAUM2002}.

\begin{lemma}\label{lemma:improvement_step_tree_decomposition}
    Let $a$ be a positive integer.
    Let $G$ be a graph and let $S$ be a nonempty subset of $V(G)$.
    Let $\mathcal{D} = \big(T,(W_x \mid x \in V(T))\big)$ be a tree decomposition of $(G,S)$ of adhesion at most $a$.
    If there exists
    a positive integer $i$,
    two vertices $x_1,x_2 \in V(T)$,
    and two sets $Z_1 \subseteq W_{x_1}, Z_2 \subseteq W_{x_2}$ with $|Z_1|=|Z_2|=i$ 
    such that there are no $i$ pairwise disjoint $(Z_1,Z_2)$-paths in $G$,
    but $|W_{z} \cap W_{z'}|\geq i$ for every $zz' \in E(T[x_1,x_2])$, then
    there exist $S' \subseteq V(G)$ with $S \subseteq S'$ and a tree decomposition $\mathcal{D}' = \big(T',(W'_z \mid z \in V(T'))\big)$ of $(G,S')$
    of adhesion at most $\max\{a,i-1\}$
    such that $s(\mathcal{D}')$ is smaller than $s(\mathcal{D})$ in the lexicographic order.
    More precisely, there exist $z \in V(T[x_1,x_2])$ and $j_0 \in \{|W_z|, \dots, |V(G)|\}$ such that
    \begin{enumerate}[label=(\alph*)]
        \item $s_j(\mathcal{D}') \leq s_j(\mathcal{D})$ for every $j \in \{j_0, \dots, |V(G)|\}$; \label{item:improvement_step_tree_decomposition_i}
        \item $s_{j_0}(\mathcal{D}') < s_{j_0}(\mathcal{D})$; and \label{item:improvement_step_tree_decomposition_ii}
        \item if $x_1 = x_2$, then there are adjacent vertices $z_1,z_2$ in $T'$ 
            such that $Z_1 \subseteq W'_{z_1}$, $Z_2 \subseteq W'_{z_2}$, and $W'_{z_1} \cap W'_{z_2}$ is a minimum 
            $(Z_1,Z_2)$-cut in $G$.
            In particular $|W'_{z_1} \cap W'_{z_2}| < i$. \label{item:improvement_step_tree_decomposition_iii}
    \end{enumerate}
    Moreover, for every positive integer $k$, if $|W_z| > \frac{3}{2}(k-1)$ or 
    there is $z' \in N_T(z)$ such that $|W_z \cap W_{z'}| \geq k$,
    then there exists $j_1 \in \{|W_z|, \dots, |V(G)|\}$ such that
    \begin{enumerate}[resume*]
        \item $s^{(k)}_j(\mathcal{D}') \leq s^{(k)}_j(\mathcal{D})$ for every $j \in \{j_1, \dots, |V(G)|\}$; and \label{item:improvement_step_tree_decomposition_i:bis}
        \item $s^{(k)}_{j_1}(\mathcal{D}') < s^{(k)}_{j_1}(\mathcal{D})$. \label{item:improvement_step_tree_decomposition_ii:bis}
    \end{enumerate}
\end{lemma}

\begin{proof}
    Let $x_1,x_2 \in V(T)$ and let $Z_1 \subseteq W_{x_1}, Z_2 \subseteq W_{x_2}$ with $|Z_1|=|Z_2|=i$ such that there are no $i$ pairwise disjoint $(Z_1,Z_2)$-paths in $G$
    but $|W_{z} \cap W_{z'}|\geq i$ for every $zz' \in E(T[x_1,x_2])$.
    By choosing such a pair $\{x_1,x_2\}$ minimizing $\dist_T(x_1,x_2)$,
    we assume that there are no $x'_1,x'_2 \in V(T[x_1,x_2])$ with $\{x'_1,x'_2\}\neq \{x_1,x_2\}$ also satisfying this property,
    that is such that there are $Z'_1 \subseteq W_{x'_1}, Z'_2 \subseteq W_{x'_2}$ with $|Z'_1|=|Z'_2|=i'$ such that there are no $i'$ pairwise disjoint $(Z'_1,Z'_2)$-paths in $G$.
    For every $u \in V(G)$, let $z_u \in V(T)$ be such that $u \in W_{z_u}$ if $u \in S$, and otherwise such that $W_{z_u}$ contains the neighborhood of the connected component
    of $u$ in $G-S$.
    In both cases, pick such a $z_u$ with $\dist_T(z_u,V(T[x_1,x_2]))$ minimum.
    For every $u \in V(G)$, let
    \[
        d_u = 
        \begin{cases}
            \dist_T\big(z_u, V(T[x_1,x_2])\big) & \text{if $u \in S$,} \\
            \dist_T\big(z_u, V(T[x_1,x_2])\big) + 1 & \text{otherwise.}
        \end{cases}
    \]
    Since there are no $i$ pairwise disjoint $(Z_1,Z_2)$-paths in $G$, by Menger's Theorem,
    there exists a set $X \subseteq V(G)$ of size at most $i-1$ such that every $(Z_1,Z_2)$-path in $G$ intersects $X$.
    Take such a set $X$ of minimum size, and among those,
    take the one minimizing $\sum_{u \in X} d_u$.

    Let $C_1$ be the union of all the connected components of $G-X$ that intersect $Z_1$, and let $C_2 = G - (X \cup V(C_1))$.
    Let $G_1 = G[V(C_1) \cup X]$ and $G_2 = G[V(C_2) \cup X]$.
    By Menger's Theorem, there is a family $(P_u \mid u \in X)$ of pairwise disjoint $(Z_1,Z_2)$-paths
    in $G$ such that $u \in V(P_u)$ for every $u \in X$.
    For every $u \in X$, $P_u$ is the union of two paths $P^1_u \subseteq G - V(C_1)$ and $P^2_u \subseteq G - V(C_2)$ 
    meeting in exactly one vertex, namely $u$.
    For each $b \in \{1,2\}$, let $T^b$ be a copy of $T$ with vertex set $\{z^b \mid z \in V(T)\}$ and edge set $\{z_1^b z_2^b \mid z_1 z_2 \in E(T)\}$,
    and let $\mathcal{D}_b = \big(T^b, (W^b_{z^b} \mid z^b \in V(T^b))\big)$ be the tree decomposition of $(G_b, (S \cap V(G_b)) \cup X)$ defined by
    \[
    W^b_{z^b} = \big(W_z \cap V(G_b)\big) \cup \big\{u \in X \mid z \in V(T[z_u, x_{3-b}])\big\}
    \]
    for every $z \in V(T)$.

    Finally, let $S' = S \cup X$, $T' = T^1 \cup T^2 \cup \{x^2_1 x^1_2\}$,
    and let $W'_{z^b} = W^b_{z^b}$ for every $z \in V(T)$ and for every $b \in \{1,2\}$.
    We claim that $\mathcal{D}' = \big(T', (W'_z \mid z \in V(T'))\big)$ satisfies the outcome of the lemma.
    
    \begin{claim}
        $\mathcal{D}'$ is a tree decomposition of $(G,S')$.
    \end{claim}

    \begin{proofclaim}
        Let $u \in S' = S \cup X$,
        and let $b \in \{1,2\}$ such that $u \in V(G_b)$.
        Then we have
        \[
            \{z \in V(T') \mid u \in W'_z\} 
            = 
            \begin{cases}
                \{z^b \mid z \in V(T), u \in W_z\} & \textrm{if $u \not\in X$,} \\
                \{z^b \mid z \in V(T), u \in W_z\} \cup \{z^b \mid z \in V(T[z_u, x_{3-b}])\} & \textrm{if $u \in X$.}
            \end{cases}
        \]
        In both cases, using the fact that $\mathcal{D}$ is a tree decomposition of $G[S]$, this set induces a connected subtree of $T'$.

        Let $C'$ be a connected component of $G-S'$.
        Since $C'$ is disjoint from $X = V(G_1) \cap V(G_2)$, there is $b \in \{1,2\}$
        such that $C' \subseteq G_b$.
        Moreover, 
        there is a connected component $C$ of $G-S$ such that $C' \subseteq C$.
        Since $\mathcal{D}$ is a tree decomposition of $(G,S)$,
        there is $z \in V(T)$ such that $N_G(V(C)) \subseteq W_z$.
        Moreover, for every $u \in N_G(V(C')) \setminus N_G(V(C))$,
        we have $u \in X$ and $z_u = z$.
        Altogether, 
        we get that $N_G(V(C')) \subseteq (N_G(V(C)) \cap V(G_b)) \cup \{u \in X \mid z_u = z\} = W'_{z^b}$.

        Let $uv \in E(G[S'])$,
        and let $b \in \{1,2\}$ such that $uv \in E(G_b)$.
        If $u,v \in S$, then $uv \in E(G[S])$,
        and so, since $\mathcal{D}$ is a tree decomposition of $(G,S)$,
        there is $z \in V(T)$ such that $u,v \in W_z$, and we deduce $u,v \in W'_{z^b}$.
        Now suppose, without loss of generality, 
        that $u \in S' \setminus S$.
        Let $C$ be the connected component of $u$ in $G-S$.
        Then, we have $u,v \in (X \cap V(C)) \cup N_G(V(C)) \subseteq \{w \in X \mid z_u=z_w\} \cup N_G(V(C)) \subseteq W'_{z^b_u}$.
        This concludes the proof that $\mathcal{D}'$ is a tree decomposition of $(G,S')$.
    \end{proofclaim}
    
    Then, we claim that the following property holds.
    \begin{equation}\label{eq:BD_0}
        \text{For all $z \in V(T)$ and $b \in \{1,2\}$, $|W^b_{z^b}| \leq |W_z|$.}
    \end{equation}
    To see this, observe that for every $u \in W^b_{z^b} \setminus W_z$, we have $u \in X$ and $z \in V(T[z_u,x_{3-b}])$.
    In particular, the subpath $P^{3-b}_u$ of $P_u$ starting at $u$ and ending
    in $Z_{3-b}$ intersects $W_z$, and so 
    we can fix some $v(u) \in V(P^{3-b}_u) \cap W_z$.
    Note that $v(u) \not\in V(G_b)\setminus X$.
    Moreover, $v(u) \neq u$ since $v(u) \in W_z$,
    and since $V(P_u) \cap X = \{u\}$, it implies $v(u) \not\in X$.
    Therefore, we have $v(u) \in W_z \setminus V(G_b) \subseteq W_z \setminus W^b_{z^b}$.
    Since the paths $P_u$ for $u \in X$ are pairwise disjoint,
    we deduce that $v(\cdot)$ is an injective mapping 
    from $W^b_{z^b} \setminus W_z$ to $W_z \setminus W^b_{z^b}$.
    This shows that $|W^b_{z^b} \setminus W_z| \leq |W_z \setminus W^b_{z^b}|$ and so $|W^b_{z^b}| = |W^b_{z^b} \setminus W_z| + |W^b_{z^b} \cap W_z|
    \leq |W_z \setminus W^b_{z^b}| + |W^b_{z^b} \cap W_z| = |W_z|$.

    The same argument also proves the following.
    \begin{equation}\label{eq:BD_0'}
        \text{For all $z_1 z_2 \in E(T)$ and $b \in \{1,2\}$, $|W^b_{z_1^b} \cap W^b_{z_2^b}| \leq |W_{z_1} \cap W_{z_2}|$.}
    \end{equation}
    In particular, $\mathcal{D}'$ has adhesion at most $\max\{a,i-1\}$.
    We now investigate the equality case in \eqref{eq:BD_0}.
    \begin{equation}\label{eq:BD_1}
        \text{For all $z \in V(T)$ and $b \in \{1,2\}$, if $|W^b_{z^b}| = |W_z|$, then $W^{3-b}_{z^{3-b}} \subseteq X$.}
    \end{equation}
    Suppose for a contradiction that $|W^b_{z^b}| = |W_z|$ and $W^{3-b}_{z^{3-b}} \not\subseteq X$ for some $z \in V(T)$ and $b \in \{1,2\}$.
    Then $W_z$ intersects $V(C_{3-b})$,
    and there is a bijection between $W_z \cap V(C_{3-b})$
    and vertices $v \in X$ such that $z \in V(T[z_{v},x_{3-b}])$.
    Let $Y = W^b_{z^b} \setminus W_z$ be the set of all these vertices.
    Then let
    \[
        X' = (X \setminus Y) \cup \big(W_z \cap V(C_{3-b})\big).
    \]
    Note that $|X'| = |X| < i$.
    We claim that $X'$ intersects every $(Z_b \cup W_z,Z_{3-b})$-path in $G$.
    Let $Q$ be a $(Z_b \cup W_z,Z_{3-b})$-path in $G$.
    Suppose for a contradiction that $V(Q) \cap X' = \emptyset$.
    The endpoint of $Q$ in $Z_b \cup W_z$ lies in $C_b$, $X$, or $C_{3-b}$.
    In all the first two cases, $Q$ intersects $X$, and in the former one
    $Q$ intersects $W_z \cap V(C_{3-b})$.
    Hence $V(Q) \cap (X \cup (W_z \cap V(C_{3-b}))) \neq \emptyset$,
    and since $V(Q) \cap X' = \emptyset$, 
    we have $V(Q) \cap X = V(Q) \cap Y \neq \emptyset$.
    Let $u_0$ be the last vertex in $V(Q) \cap Y$ along $Q$.
    Then the subpath $Q]u_0, \term(Q)]$ of $Q$ is included in $C_{3-b}$,
    and $Q[u_0, \term(Q)]$ contains a $(W_{z_{u_0}}, W_{x_{3-b}})$-path (in the case $u_0 \in S$, $Q[u_0, \term(Q)]$ is this subpath, while this subpath may be shorter in the case $u_0 \not\in S$).
    Since $z \in V(T[z_{u_0},x_{3-b}])$, $Q[u_0, \term(Q)]$ intersects $W_z$.
    As $u_0 \not\in W_z$, $Q$ intersects $W_z$ in $V(C_{3-b})$, and so $V(Q) \cap X' \neq \emptyset$.
    This proves that $X'$ intersects every $(Z_b \cup W_z, Z_{3-b})$-path in $G$.

    If $z \in V(T]x_b,x_{3-b}])$, then take any $Z \subseteq W_z$ of size $i$.
    Then $X'$ intersects every $(Z,Z_{3-b})$-path in $G$, and so
    there are no $i$ pairwise disjoint $(Z,Z_{3-b})$-paths in $G$.
    Therefore, the pair $\{z,x_{3-b}\}$ contradicts the choice of $\{x_1,x_2\}$.
    Now suppose $z \not\in V(T]x_b, x_{3-b}])$.
    We claim that $d_{u'} < d_u$ for every $u' \in X' \setminus X$ and for every $u \in X \setminus X' = Y$.
    Indeed, $z \in V(T[z_{u}, x_{3-b}])$ by definition of $W^b_{z^b}$, and as $z \not\in V(T]x_b,x_{3-b}])$,
    $z$ lies in the path from $z_{u}$ to $V(T[x_1,x_2])$ in $T$.
    Recall also that $u \not\in W_z$ so $z_u \neq z$.
    For both cases $u \in S$ and $u \not\in S$, this implies that $d_u > \dist_T(z, V(T[x_1,x_2]))$.
    On the other hand, since $u' \in W_z$ and
    $z_{u'}$ is such that $u' \in W_{z_{u'}}$ and $\dist_T(z_{u'}, V(T[x_1,x_2]))$ is minimum,
    we have $d_{u'} \leq \dist(z, V(T[x_1,x_2]))$.
    It follows that $d_u > \dist_T(z, V(T[x_1,x_2])) \geq d_{u'}$ as claimed. 
    Since $X' \setminus X$ and $X \setminus X'$ are nonempty, this contradicts the minimality of $\sum_{u \in X} d_u$.
    This proves \eqref{eq:BD_1}.

    \medskip

    We now prove the following property.
    \begin{equation}\label{eq:BD_2}
        \text{There exists $z \in V(T[x_1,x_2])$ such that $|W^1_{z^1}|, |W^2_{z^2}| < |W_z|$.}
    \end{equation}

    By \eqref{eq:BD_1}, it is enough to find $z \in V(T[x_1,x_2])$ such that $W_z$ intersects both $V(C_1)$ and $V(C_2)$.
    Suppose for a contradiction that for every $z \in V(T[x_1,x_2])$, $W_z \setminus X \subseteq V(C_1)$ or $W_z \setminus X \subseteq V(C_2)$.
    Then there is an edge $zz'$ in $T[x_1,x_2]$ such that $W_z  \setminus X \subseteq V(C_1)$ and $W_{z'} \setminus X \subseteq V(C_2)$.
    It follows that $W_z \cap W_{z'} \subseteq X$ and so $|W_z \cap W_{z'}| < i$, a contradiction.
    This proves \eqref{eq:BD_2}.

    To conclude, fix $z$ as in \eqref{eq:BD_2}.
    Let $j_0 \in \{|W_z|, \dots, |V(G)|\}$ be maximum such that there exists $z' \in V(T)$ satisfying $|W_{z'}|=j_0$ and $|W^1_{{z'}^1}|,|W^2_{{z'}^2}| < |W_{z'}|$.
    In particular, $j_0 \geq |W_z| \geq |X|+1$.
    Then, \eqref{eq:BD_0} and \eqref{eq:BD_1} imply that for every $z'' \in V(T)$ with $|W_{z''}|\geq j_0$,
    $W'_{z''^1}$ and $W'_{z''^2}$ have size at most $|W_{z''}|$,
    and at most one of them has size exactly $|W_{z''}|$.
    Hence $s_j(\mathcal{D}') \leq s_j(\mathcal{D})$ for every $j \in \{j_0, \dots, |V(G)|\}$,
    and $s_{j_0}(\mathcal{D}') < s_{j_0}(\mathcal{D})$.
    This proves \ref{item:improvement_step_tree_decomposition_i} and \ref{item:improvement_step_tree_decomposition_ii}.
    Finally, if $x_1=x_2$, then \ref{item:improvement_step_tree_decomposition_iii} follows from the definition of $\mathcal{D}'$ for $z_1 = x_1^2$ and $z_2 = x_2^1$.
    This proves the first part of the lemma.

    \medskip

    We now prove the ``moreover'' part of the lemma.
    Suppose that $|W_z| > \frac{3}{2}(k-1)$ or there is $z' \in N_T(z)$ with $|W_z \cap W_{z'}| \geq k$.
    Then let $j_1 \in \{|W_z|, \dots, |V(G)|\}$ be maximum such that there exists $z'' \in V(T)$
    satisfying $|W_{z''}|=j_1$, $|W^1_{{z''}^1}|,|W^2_{{z''}^2}| < |W_{z''}|$, and $|W_{z''}| > \frac{3}{2}(k-1)$ or there is $z' \in N_T(z'')$ with $|W_{z''} \cap W_{z'}| \geq k$.
    Just as before, \eqref{eq:BD_0}, \eqref{eq:BD_0'} and \eqref{eq:BD_1}
    imply that $s^{(k)}_j(\mathcal{D}') \leq s^{(k)}_j(\mathcal{D})$ for every $j \in \{j_1, \dots, |V(G)|\}$,
    and $s^{(k)}_{j_1}(\mathcal{D}') < s^{(k)}_{j_1}(\mathcal{D})$.
    This proves \ref{item:improvement_step_tree_decomposition_i:bis} and \ref{item:improvement_step_tree_decomposition_ii:bis}.
\end{proof}

As a first application of this lemma, we show the following version of a result of
Cygan, Komosa, Lokshtanov, Pilipczuk, Pilipczuk, Saurabh, and Wahlstr\"om~\cite{Cygan2020}.

\begin{theorem}\label{thm:k_unbreakable_treedecomposition}
    Let $k$ be a positive integer and let $G$ be a graph.
    There is a tree decomposition $\big(T,(W_x \mid x \in V(T)) \big)$ of $G$ of adhesion at most $k-1$
    such that
    \begin{enumerate}[label=(\alph*)]
        \item for every integer $i$ with $1 \leq i \leq k$, \label{item:thm:k_unbreakable_treedecomposition:i}
            for every $x \in V(T)$, for every $Z_1,Z_2 \subseteq W_x$ with $|Z_1|=|Z_2|=i$, 
            there are $i$ pairwise disjoint $(Z_1,Z_2)$-paths in $G$; and
        \item if $G$ is connected, \label{item:thm:k_unbreakable_treedecomposition:ii}
            then for every edge $x_1x_2 \in E(T)$, $\bigcup_{z \in V(T_{x_2 \mid x_1})} W_z$ induces a connected subgraph of $G$. 
    \end{enumerate}
\end{theorem}

\begin{proof}
    Let $\mathcal{D} = \big(T,(W_x \mid x \in V(T))\big)$ be a tree decomposition of $G$ of adhesion less than $k$.
    Note that such a tree decomposition always exists since the tree decomposition of $G$ with a single bag $V(G)$ has adhesion less than $k$.
    We take such a tree decomposition $\mathcal{D}$ such that $s(\mathcal{D})$ is minimum in the lexicographic order.

    First we show \ref{item:thm:k_unbreakable_treedecomposition:i}.
    Suppose for a contradiction that 
    there exists $i$ with $1 \leq i \leq k$, 
    $x \in V(T)$, and $Z_1,Z_2 \subseteq W_x$ with $|Z_1|=|Z_2|=i$, 
    such that there are no $i$ pairwise disjoint $(Z_1,Z_2)$-paths in $G$.
    Then by Lemma~\ref{lemma:improvement_step_tree_decomposition} applied for $S = V(G)$, there exists
    $\mathcal{D}'$ a tree decomposition of $G$ of adhesion at most $\max\{k-1,i-1\} = k-1$
    such that $s(\mathcal{D}')$ is smaller than $s(\mathcal{D})$ in the lexicographic order.
    This contradicts the minimality of $s(\mathcal{D})$ and so \ref{item:thm:k_unbreakable_treedecomposition:i} holds.

    Now suppose that \ref{item:thm:k_unbreakable_treedecomposition:ii} does not hold.
    Then $G$ is connected and there is an edge $x_1x_2 \in E(T)$ such that $G' = G\left[\bigcup_{z \in V(T_{x_2 \mid x_1})} W_z\right]$ is not connected.
    Let $C_1$ be a connected component of $G'$ and let $C_2 = G' - V(C_1)$.
    For each $a \in \{1,2\}$, let $T^a$ be a copy of $T_{x_2 \mid x_1}$ with vertex set $\{z^a \mid z \in V(T_{x_2 \mid x_1})\}$
    and edge set $\{z_1^a z_2^a \mid z_1z_2 \in E(T_{x_2 \mid x_1})\}$.
    Then, for every $z \in V(T_{x_2 \mid x_1})$, let $W'_{z^a} = W_z \cap V(C_a)$, and
    for every $z \in V(T_{x_1 \mid x_2})$, let $W'_z = W_z$.
    Finally, for $T' = T_{x_1 \mid  x_2} \cup T^1 \cup T^2 \cup \{x_1 x^1_2, x_1 x^2_2\}$,
    $\mathcal{D}' = \big(T', (W'_z \mid z \in V(T'))\big)$ is a tree decomposition of $G$ of adhesion at most $k-1$.
    We now prove that $s(\mathcal{D}')$ is smaller than $s(\mathcal{D})$, which contradicts the minimality of $s(\mathcal{D})$ and so proves the theorem.
    To do so,
    consider $z_0 \in V(T_{x_2 \mid x_1})$ be such that $W_{z_0}$ intersects both $V(C_1)$ and $V(C_2)$, and $|W_{z_0}|$ is maximum under this property.
    This is well-defined because $W_{x_2}$ intersects both $V(C_1)$ and $V(C_2)$ since $G$ is connected.
    Now, for every $z \in V(T_{x_2 \mid x_1})$, if $|W_z| > |W_{z_0}|$, then by maximality of $|W_{z_0}|$, one of $W'_{z^1}$ and $W'_{z^2}$ is empty, and the other one is $W_z$.
    Moreover, if $|W_z| = |W_{z_0}|$, then either one of $W'_{z^1}$ and $W'_{z^2}$ is empty, and the other one is $W_z$;
    or $W_z$ intersects both $V(C_1)$ and $V(C_2)$, and it follows that $|W'_{z^1}|,|W'_{z^2}| < |W_z| = |W_{z_0}|$.
    This second case occurs at least once for $z = z_0$.
    Altogether, this implies that $s_j(\mathcal{D}') \leq s_j(\mathcal{D})$ for every $j \in \{|W_{z_0}|+1, \dots, |V(G)|\}$,
    and $s_{|W_{z_0}|}(\mathcal{D}') < s_{|W_{z_0}|}(\mathcal{D})$.
    Therefore, $s(\mathcal{D}')$ is smaller than $s(\mathcal{D})$ in the lexicographic order.
    This contradicts the minimality of $s(\mathcal{D})$ and so \ref{item:thm:k_unbreakable_treedecomposition:ii} holds.
\end{proof}

We can now show the main result of this section.

\begin{lemma}\label{lemma:existence_good_tree_decomposition}
    Let $k,a,t$ be positive integers with $a \geq k$.
    Let $G$ be a graph and let $S \subseteq V(G)$.
    If there exists a tree decomposition of $(G,S)$ of width less than $t$ and adhesion at most $a$, then
    there exists $S' \subseteq V(G)$ with $S \subseteq S'$ and $\mathcal{D} = \big(T,(W_x \mid x \in V(T))\big)$ a tree decomposition of $(G,S')$ such that
    \begin{enumerate}[label=(\alph*)]
        \item $\mathcal{D}$ has adhesion at most $a$; \label{item:lemma:existence_good_tree_decomposition:i}
        \item $\mathcal{D}$ has width less than $\max\big\{t, \frac{3}{2}(k-1)\big\}$; \label{item:lemma:existence_good_tree_decomposition:ii}
        \item for every $x_1,x_2 \in V(T)$, \label{item:lemma:existence_good_tree_decomposition:iii}
            for every $Z_1 \subseteq W_{x_1}$ and $Z_2 \subseteq W_{x_2}$ both of size $k$,
            either
            \begin{enumerate}[label=(\roman*)]
                \item there are $k$ pairwise disjoint $(Z_1,Z_2)$-paths in $G$,
                \item there exists $z_1z_2 \in E(T[x_1,x_2])$ with $|W_{z_1} \cap W_{z_2}| < k$, or
                \item $x_1=x_2$, $|W_{x_1}| \leq \frac{3}{2}(k - 1)$, and $|W_{x_1} \cap W_y|<k$ for every $y \in N_T(x_1)$;
            \end{enumerate}
        \item for every $x_1 x_2 \in E(T)$ with $|W_{x_1} \cap W_{x_2}|<k$,  \label{item:lemma:existence_good_tree_decomposition:iv}
            if there exists $x_3 \in N_T(x_2)$ such that $|W_{x_2} \cap W_{x_3}| \geq k$, then
            for every positive integer $i$, for every $Z_1,Z_2 \subseteq W_{x_1} \cap W_{x_2}$ both of size $i$,
            there are $i$ pairwise disjoint $(Z_1,Z_2)$-paths in 
            \[
                G\left[\bigcup_{z \in V(T_{x_2 \mid x_1})} W_z \cup \bigcup_{C \in \mathcal{C}(x_2 \mid x_1)} V(C) \right] \setminus \binom{W_{x_1} \cap W_{x_2}}{2}
            \]
            where $\mathcal{C}(x_2 \mid x_1)$ is the family of all the components $C$ of $G-S'$ such that $N_{G}(V(C)) \subseteq \bigcup_{z \in V(T_{x_2 \mid x_1})}W_z$ and $N_{G}(V(C)) \not\subseteq W_{x_1} \cap W_{x_2}$.
    \end{enumerate}
\end{lemma}

\begin{figure}[ht]
    \centering
    \includegraphics{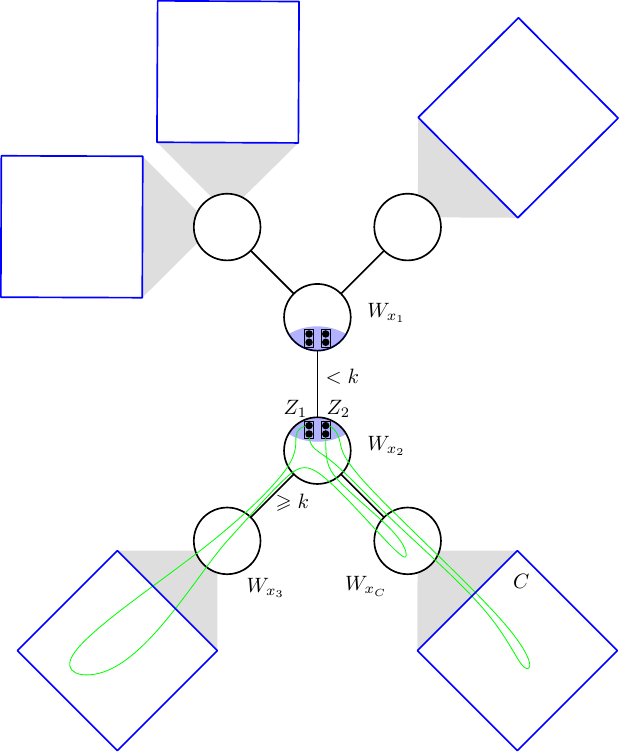}
    \caption{Illustration for Item~\ref{item:lemma:existence_good_tree_decomposition:iv} in Lemma~\ref{lemma:existence_good_tree_decomposition}.}
    \label{fig:good_tree_decomposition}
\end{figure}

See Figure~\ref{fig:good_tree_decomposition} for an illustration of \ref{item:lemma:existence_good_tree_decomposition:iv}.

\begin{proof}
    Let $S' \subseteq V(G)$ contain $S$ and $\mathcal{D}' = \big(T, (W_x \mid x \in V(T))\big)$ be a tree decomposition of $(G,S')$
    of adhesion at most $a$ and width less than $\max\{t,\frac{3}{2}(k-1)\}$.
    Consider such a pair $(S',\mathcal{D}')$ such that
    $s^{(k)}(\mathcal{D}') = \big(s^{(k)}_{|V(G)|}(\mathcal{D}'), \dots, s^{(k)}_0(\mathcal{D}')\big)$ is minimum in the lexicographic order.
    We claim that $\mathcal{D}'$ satisfies \ref{item:lemma:existence_good_tree_decomposition:i}-\ref{item:lemma:existence_good_tree_decomposition:iv}.
    By definition, \ref{item:lemma:existence_good_tree_decomposition:i} and \ref{item:lemma:existence_good_tree_decomposition:ii} hold.

    \medskip

    Now suppose that \ref{item:lemma:existence_good_tree_decomposition:iii} does not hold:
    there exist $x_1,x_2 \in V(T)$, 
    $Z_1 \subseteq W_{x_1}$ and $Z_2 \subseteq W_{x_2}$ both of size $k$,
    such that
    \begin{enumerate}[label=(\roman*)]
        \item there are no $k$ pairwise disjoint $(Z_1,Z_2)$-paths in $G$;  \label{item:existence_good:(c):i}
        \item for every $z_1,z_2 \in E(T[x_1,x_2])$, $|W_{z_1} \cap W_{z_2}| \geq k$; and \label{item:existence_good:(c):ii}
        \item $x_1\neq x_2$, or $|W_{x_1}| > \frac{3}{2}(k - 1)$, or there exists $y \in N_T(x_1)$ such that $|W_{x_1} \cap W_y| \geq k$.  \label{item:existence_good:(c):iii}
    \end{enumerate}
    Then by Lemma~\ref{lemma:improvement_step_tree_decomposition},
    there exists a set $S'' \subseteq V(G)$ containing $S'$ (and so $S$), a tree decomposition $\mathcal{D}''$
    of $(G,S'')$ of adhesion at most $\max\{a,k-1\} = a$,
    and a vertex $z \in V(T[x_1,x_2])$ 
    satisfying the outcome of Lemma~\ref{lemma:improvement_step_tree_decomposition}.
    By minimality of $s^{(k)}(\mathcal{D}')$,
    we have that $s^{(k)}(\mathcal{D}'')$ is at least $s^{(k)}(\mathcal{D}')$ in the lexicographic order.
    Therefore, we can not have both \ref{lemma:improvement_step_tree_decomposition}.\ref{item:improvement_step_tree_decomposition_i:bis} and \ref{lemma:improvement_step_tree_decomposition}.\ref{item:improvement_step_tree_decomposition_ii:bis}.
    It follows from the outcome of Lemma~\ref{lemma:improvement_step_tree_decomposition} that $|W_z| \leq \frac{3}{2}(k-1)$, and for every $z' \in N_T(z)$, we have $|W_z \cap W_{z'}| < k$.
    By \ref{item:existence_good:(c):ii}, and because $z \in V(T[x_1,x_2])$,
    we deduce that $x_1=z=x_2$.
    Altogether, this gives that $x_1=x_2$, $|W_{x_1}| \leq \frac{3}{2}(k-1)$,
    and for every $y \in N_T(x_1)$, $|W_{x_1} \cap W_y| < k$,
    which contradicts \ref{item:existence_good:(c):iii}.
    This proves \ref{item:lemma:existence_good_tree_decomposition:iii}.

    \medskip

    Now suppose that \ref{item:lemma:existence_good_tree_decomposition:iv} does not hold.
    Let $\mathcal{C}(x_2 \mid x_1)$ be the family of all the components $C$ of $G-S'$ such that $N_{G}(V(C)) \subseteq \bigcup_{z \in V(T_{x_2 \mid x_1})}W_z$ and $N_{G}(V(C)) \not\subseteq W_{x_1} \cap W_{x_2}$.
    Hence there exist $x_1 x_2 \in E(T)$ with $|W_{x_1} \cap W_{x_2}|<k$, 
    $x_3 \in N_T(x_2)$ with $|W_{x_2} \cap W_{x_3}| \geq k$,
    a positive integer $i$, 
    and two sets $Z_1,Z_2 \subseteq W_{x_1} \cap W_{x_2}$ both of size $i$,
    such that there are no $i$ pairwise disjoint $(Z_1,Z_2)$-paths in 
    \[
        G_0 = G\left[\textstyle\bigcup_{z \in V(T_{x_2 \mid x_1})} W_z \cup \bigcup_{C \in \mathcal{C}(x_2 \mid x_1)} V(C) \right] 
        \setminus \binom{W_{x_1} \cap W_{x_2}}{2}.
    \]
    By possibly replacing $Z_a$ by $Z_a \cup ((W_{x_1} \cap W_{x_2}) \setminus (Z_1 \cup Z_2))$ for every $a \in \{1,2\}$,
    and $i$ by $i + |(W_{x_1} \cap W_{x_2}) \setminus (Z_1 \cup Z_2)|$, we assume that $Z_1 \cup Z_2 = W_{x_1} \cap W_{x_2}$.

    We call Lemma~\ref{lemma:improvement_step_tree_decomposition} on
    the graph $G_0$,
    the set $S_0 = V(G_0) \cap S'$ of vertices, 
    the tree decomposition $\mathcal{D}_0 = \big(T_{x_2 \mid x_1}, (W_z \mid x \in V(T_{x_2 \mid x_1}))\big)$, 
    the pair of vertices $(x_2,x_2)$,
    and the subsets $(Z_1,Z_2)$ of vertices.
    We deduce that there exists $S'_0 \subseteq V(G_0)$ with $S_0 \subseteq S'_0$ and 
    a tree decomposition $\mathcal{D}'_0 = \big(T_0,(W'_z \mid z \in V(T'))\big)$ of $(G_0,S'_0)$
    of adhesion at most $\max\{a,i-1\}$
    such that
    \begin{enumerate}[label=\ref{lemma:improvement_step_tree_decomposition}.(\alph*)]
        \setcounter{enumi}{2}
        \item there are adjacent vertices $z_1,z_2$ in $T_0$ 
            such that $Z_1 \subseteq W'_{z_1}$, $Z_2 \subseteq W'_{z_2}$, and $W'_{z_1} \cap W'_{z_2}$ is a minimum 
            $(Z_1,Z_2)$-cut in $G$.
            In particular $|W'_{z_1} \cap W'_{z_2}| < i$.
    \end{enumerate}
    Moreover, there exists $j_1 \in \{|W_{x_2}|, \dots, |V(G)|\}$ such that
    \begin{enumerate}[resume*]
        \item $s^{(k)}_j(\mathcal{D}'_0) \leq s^{(k)}_j(\mathcal{D}_0)$ for every $j \in \{j_1, \dots, |V(G)|\}$; and
        \item $s^{(k)}_{j_1}(\mathcal{D}'_0) < s^{(k)}_{j_1}(\mathcal{D}_0)$.
    \end{enumerate}

\begin{figure}[ht]
    \centering
    \includegraphics{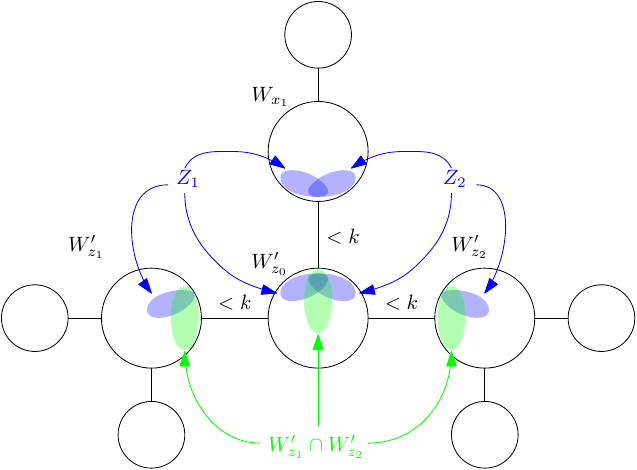}
    \caption{Construction of $\mathcal{D}''$ in the proof of Lemma~\ref{lemma:existence_good_tree_decomposition}.}
    \label{fig:existance_good_tree_decomposition_proof}
\end{figure}

    Now let $T_1$ be obtained from $T_0$ by subdividing once the edge $z_1 z_2$.
    Let $z_0$ be the resulting new vertex.
    Then let $W'_{z_0} = Z_1 \cup Z_2 \cup (W'_{z_1} \cap W'_{z_2})$.
    Since $W'_{z_1} \cap W'_{z_2}$ is a minimum $(Z_1,Z_2)$-cut,
    we have 
    \[
        |(W'_{z_1} \cap W'_{z_2}) \setminus (Z_1 \cap Z_2)| \leq \min\{|Z_1 \setminus Z_2|, |Z_2 \setminus Z_1|\} \leq \tfrac{k-1}{2}.
    \]
    This implies that $|W'_{z_0}| \leq \frac{3}{2}(k-1)$.
    Moreover, for each $b \in \{1,2\}$
    \begin{align*}
        |W'_{z_0} \cap W'_{z_b}| 
        &= |(W'_{z_1} \cap W'_{z_2}) \cup Z_b| \\
        &\leq |(W'_{z_1} \cap W'_{z_2}) \setminus (Z_1 \cap Z_2)| + |Z_b| \\
        &\leq |Z_{3-b} \setminus Z_b| + |Z_b| \\
        &= |Z_1 \cup Z_2| < k
    \end{align*}
    where the last inequality follows from the inclusion $Z_1 \cup Z_2 \subseteq W_{x_1} \cap W_{x_2}$ and $|W_{x_1} \cap W_{x_2}| < k$.
    Then let $S'' = (S' \setminus V(G_0)) \cup S'_0$, and let $\mathcal{D}'' = \big(T'', (W''_z \mid z \in V(T''))\big)$
    where $V(T'') = V(T_{x_1 \mid x_2}) \cup V(T_1)$, $E(T'') = E(T_{x_1 \mid x_2}) \cup E(T_1) \cup \{x_1 z_0\}$,
    $W''_z = W_z$ for every $z \in V(T_{x_1 \mid x_2})$, and $W''_z = W'_z$ for every $z \in V(T_1)$.
    See Figure~\ref{fig:existance_good_tree_decomposition_proof}.
    Then $\mathcal{D}''$ is a tree decomposition of $(G,S'')$ of adhesion at most $\max\{a, i-1\}$ and width less than
    $\max\big\{t, \frac{3}{2}(k-1)\big\}$.
    Moreover, since every neighbor $z$ of $z_0$ is such that $|W''_{z_0} \cap W''_{z}| < k$ and $|W''_{z_0}| \leq \frac{3}{2}(k-1)$,
    the bag $W''_{z_0}$ will not be counted in $s^{(k)}(\mathcal{D}'')$, and it follows that
    $s^{(k)}_{j}(\mathcal{D}'') \leq s^{(k)}_{j}(\mathcal{D}')$ for every $j \in \{j_1, \dots, |V(G)|\}$
    and $s^{(k)}_{j_1}(\mathcal{D}'') < s^{(k)}_{j_1}(\mathcal{D}')$.
    Therefore, $s^{(k)}(\mathcal{D}'')$ is smaller than $s^{(k)}(\mathcal{D}')$ in the lexicographic order, a contradiction.
    This proves that \ref{item:lemma:existence_good_tree_decomposition:iv} holds.
\end{proof}

We finish this section with a proof that torso of nice pairs
in a graph of bounded treewidth have bounded treewidth.

\begin{lemma}\label{lemma:tw_of_torso_of_good_pairs}
    Let $t$ be a positive integer,
    and let $G$ be a graph with $\tw(G) < t$.
    For every nice pair $(U, \mathcal{B})$ in $G$,
    \[
        \tw(\torso_G(U,\mathcal{B})) < 2t.
    \]
\end{lemma}

\begin{proof}
    Let $(U,\mathcal{B})$ be a nice pair in $G$.
    By Theorem~\ref{thm:k_unbreakable_treedecomposition}
    applied for $k$ large enough (say $k = |U|$)
    that there is a tree decomposition $\big(T, (W_x \mid x \in V(T))\big)$
    of $\torso_G(U,\mathcal{B})$ such that (in particular):
    \begin{equation*}
        \begin{array}{c}
            \text{For all $y \in V(T)$, $k' \in \NN$, $Z_1, Z_2 \subseteq W_y$ with $|Z_1|=|Z_2|=k'$,} \\
            \text{there are $k'$ pairwise disjoint $(Z_1,Z_2)$-paths in $\torso_G(U, \mathcal{B})$.}
        \end{array}
    \end{equation*}
    Since $(U, \mathcal{B})$ is a nice pair in $G$
    and by Lemma~\ref{lemma:finding_k_paths_between_sets_connected_in_the_torso_of_a_nice_set},
    this implies that the same holds in $G$:
    \begin{equation}\label{eq:Wy_is_well_linked_in_G'}
        \begin{array}{c}
        \text{For all $y \in V(T)$, $k' \in \NN$, $Z_1, Z_2 \subseteq W_y$ with $|Z_1|=|Z_2|=k'$,} \\
        \text{there are $k'$ pairwise disjoint $(Z_1,Z_2)$-paths in $G$.}
        \end{array}
    \end{equation}

    Consider some $y \in V(T)$,
    and let $\mathcal{F}$ be the family of all the connected subgraphs of $G$
    containing more than $|W_y|/2$ vertices of $W_y$.
    Clearly, there are no two pairwise disjoint members of $\mathcal{F}$,
    and so by Lemma~\ref{lemma:helly_property} applied to a tree decomposition
    of $G$ of width less than $t$,
    there is a set $X \subseteq V(G)$ of size at most $t$ that intersects 
    every member of $\mathcal{F}$.
    Fix such a set $X$.
    For every $u,v \in W_y \setminus X$,
    \eqref{eq:Wy_is_well_linked_in_G'} applied to the pair
    $(Z_1,Z_2) = (X \cup \{u\}, X \cup \{v\})$
    implies that there is a $(u,v)$-path in $G$ disjoint from $X$.
    Therefore, $W_y \setminus X$ lies in a single connected component of $G - X$.
    Since this component is disjoint from $X$,
    it is not a member of $\mathcal{F}$, and so $|W_y \setminus X| \leq |W_y|/2$.
    Therefore,
    \[
        |W_y| \leq |X| + |W_y \setminus X| \leq t + |W_y|/2,
    \]
    and so $|W_y| \leq 2t$.

    This shows that $\big(T,(W_x \mid x \in V(T))\big)$ is a tree decomposition
    of $\torso_G(U,\mathcal{B})$ of width less than $2t$, which proves the lemma.
\end{proof}

\section{Graphs excluding all \texorpdfstring{$k$}{k}-ladders of length \texorpdfstring{$\ell$}{l} have bounded \texorpdfstring{$k$}{k}-treedepth}\label{sec:graphs_excluding_k_ladders_have_bounded_tdk}

In this section, we prove the ``if'' part of Theorem~\ref{thm:main}, 
that is the fact that graphs excluding all $k$-ladders of length $\ell$ as minors have bounded $k$-treedepth.
To do so, we will prove the following statement by induction.

\begin{theorem}\label{thm:main_induction}
    There is a function $f_{\cstref{thm:main_induction}}\colon \NN^4 \to \NN$ such that the following holds.
    Let $k,\ell,a,t$ be positive integers.
    For every graph $G$ that does not contain any $k$-ladder of length $\ell$ as a minor,
    for every nice pair $(U, \mathcal{B})$ in $G$,
    if $S \subseteq U$ is such that there is a tree decomposition of $(\torso_G(U, \mathcal{B}),S)$ of adhesion at most $a$ and width less than $t$,
    then there exists $S'$ with $S \subseteq S' \subseteq U$ such that
    \[
    \td_k(\torso_G(U, \mathcal{B}), S') \leq f_{\cstref{thm:main_induction}}(k, \ell, t, a).
    \]
\end{theorem}

We quickly show how to deduce Theorem~\ref{thm:main} from Theorem~\ref{thm:main_induction}.
First observe that it is enough to consider minor-closed classes.
Proposition~\ref{prop:lower_bound} implies that 
every minor-closed class of graphs of bounded $k$-treedepth excludes
$T \square P_\ell$ for every tree $T$ on $k$ vertices, for some integer $\ell$.
For the other direction, 
if $\mathcal{C}$ is a minor-closed class of graphs that does not contain
$T \square P_\ell$ for every tree $T$ on $k$ vertices,
then Observation~\ref{obs:TsquareP_from_a_long_ladder} implies that 
no graph in $\mathcal{C}$ contains a $k$-ladder of length 
$k^{k-2} (\ell-1)$ as a minor,
and Theorem~\ref{thm:grid_minor_thm} implies that every graph in $\mathcal{C}$
has treewidth less than $f_{\cstref{thm:grid_minor_thm}}(\max\{k, \ell\})$.
Hence Theorem~\ref{thm:main_induction} applied for $U=S=V(G)$ and 
$\mathcal{B} = \emptyset$ implies that a graph in such a class has $k$-treedepth at most
$f_{\cstref{thm:main_induction}}\big(k, k^{k-2} (\ell-1), f_{\cstref{thm:grid_minor_thm}}(\max\{k,\ell\}), f_{\cstref{thm:grid_minor_thm}}(\max\{k,\ell\})\big)$.

Before showing Theorem~\ref{thm:main_induction}, we give some notation and prove a key lemma.
Let $G$ be a graph, let $X \subseteq V(G)$ and let $u \in V(G)$.
The \emph{projection} of $u$ on $X$ is the set
\[
    \Pi_G(X,u) =
    \begin{cases}
        \{u\} &\textrm{if $u \in X$}\\
        N_G(V(C)) &\textrm{if $u \not\in X$ and $C$ is the connected component of $u$ in $G-X$.}
    \end{cases}
\]
Moreover, for every $A \subseteq V(G)$, let $\Pi_G(X,A) = \bigcup_{u \in A} \Pi_G(X,u)$.
Observe that for every connected subgraph $H$ of $G - \Pi_G(X,A)$ intersecting $X$, $V(H)$ is disjoint from $A$.

The following lemma is inspired by Huynh, Joret, Micek, Seweryn, and Wollan's proof of Theorem~\ref{thm:main} for the case $k=2$~\cite[proof of Lemma~11]{Huynh2021}.

\begin{lemma}\label{lemma:connecting_a_small_set}
    There is a function $f_{\cstref{lemma:connecting_a_small_set}}\colon \NN^3 \to \NN$ such that the following holds.
    Let $k,\ell,c$ be positive integers.
    Let $G$ be a connected graph that does not contain any $k$-ladder of length $\ell$ as a minor,
    let $(U, \mathcal{B})$ be a nice pair in $G$,
    and let $R \subseteq U$ with $|R| \leq c$.
    There is a set $Z \subseteq U$ with $R \subseteq Z$ such that
    \begin{enumerate}[label=(\alph*)]
        \item for every connected component $C$ of $\torso_G(U, \mathcal{B})-Z$, $R$ is connected in $\torso_G(U, \mathcal{B})-V(C)$; and \label{item:connecting_lemma_i}
        \item $\td_k(\torso_G(U, \mathcal{B}), Z) \leq f_{\cstref{lemma:connecting_a_small_set}}(k,\ell,c)$. \label{item:connecting_lemma_ii}
    \end{enumerate}
\end{lemma}

\begin{proof}
    Let 
    \begin{align*}
        f_{\cstref{lemma:connecting_a_small_set}}(k,\ell,c) 
        = c \cdot \big(2(k+1)k + 2f_{\cstref{lemma:finding_kladder_from_a_mess}}(k,\ell) \cdot f_{\cstref{thm:grid_minor_thm}}(\max\{k,\ell+1\})^2\big).
    \end{align*}

    When $c=1$, the result is clear for $Z = R$.
    Now suppose $c \geq 2$.
    Since $G$ has no $k$-ladder of length $\ell$,
    the $\max\{k,\ell+1\} \times \max\{k, \ell+1\}$ grid is not a minor of $G$.
    It follows from Theorem~\ref{thm:grid_minor_thm} that
    $\tw(G) < f_{\cstref{thm:grid_minor_thm}}(\max\{k,\ell+1\})$.
    Since $(U, \mathcal{B})$ is a nice pair in $G$,
    Lemma~\ref{lemma:tw_of_torso_of_good_pairs} implies that
    \begin{equation}\label{eq:tw_torsoUB}
        \tw(\torso_G(U,\mathcal{B})) < 2f_{\cstref{thm:grid_minor_thm}}(\max\{k,\ell+1\}).
    \end{equation}
    Fix an ordering $r_1, \dots, r_{|R|}$ of $R$.
    Note that since $G$ is connected, $\torso_G(U,\mathcal{B})$ is connected.
    By Theorem~\ref{thm:k_unbreakable_treedecomposition} applied to $\torso_G(U, \mathcal{B})$,
    there is a tree decomposition $\big(T,(W_x \mid x \in V(T)) \big)$ of $\torso_G(U, \mathcal{B})$
    of adhesion at most $k-1$
    such that
    \begin{enumerate}[label=\ref{thm:k_unbreakable_treedecomposition}.(\alph*)]
        \item for every integer $i$ with $1 \leq i \leq k$,  \label{item:lemma:connecting_a_small_set:thm:k_unbreakable_treedecomposition:i}
            for every $x \in V(T)$, for every $Z_1,Z_2 \subseteq W_x$ with $|Z_1|=|Z_2|=i$, 
            there are $i$ pairwise disjoint $(Z_1,Z_2)$-paths in $\torso_G(U, \mathcal{B})$; and
        \item for every edge $x_1x_2 \in E(T)$, the set $\bigcup_{z \in V(T_{x_1 \mid x_2})} W_z$ 
            induces a connected subgraph of $\torso_G(U, \mathcal{B})$.  \label{item:lemma:connecting_a_small_set:thm:k_unbreakable_treedecomposition:ii}
    \end{enumerate}
    
    Fix some $i \in [|R|-1]$.
    We will build a set $Z_i \subseteq V(G)$ containing $r_i$ and $r_{i+1}$ such that
    \begin{enumerate}[label=(a\arabic*)]
        \item for every connected component $C$ of $\torso_G(U,\mathcal{B})-Z_i$, $r_i$ and $r_{i+1}$ are connected in $\torso_G(U, \mathcal{B})-V(C)$; and \label{item:proof_connecting_lemma_Zi_i}
        \item $\td_k(\torso_G(U, \mathcal{B}), Z_i) \leq 2(k+1)k + 2f_{\cstref{lemma:finding_kladder_from_a_mess}}(k,\ell) \cdot f_{\cstref{thm:grid_minor_thm}}(\max\{k,\ell+1\})^2$.
        \label{item:proof_connecting_lemma_Zi_ii}
    \end{enumerate}
    Let $x_1,x_2 \in V(T)$ be such that $r_i \in W_{x_1}$ and
    and $r_{i+1} \in W_{x_2}$.
    Consider some $y \in V(T[x_1,x_2])$.
    We will build a set $Z_{i,y} \subseteq W_y$ such that
    \begin{enumerate}[label=(b\arabic*)]
        \item for every connected component $C$ of $\torso_G(U, \mathcal{B})-Z_{i,y}$ intersecting $W_y$, $r_i$ and $r_{i+1}$ are connected in $\torso_G(U, \mathcal{B})-V(C)$; \label{item:proof_connecting_lemma_Ziy_i}
        \item $|Z_{i,y}| \leq 2(k+1)k + 2f_{\cstref{lemma:finding_kladder_from_a_mess}}(k,\ell) \cdot f_{\cstref{thm:grid_minor_thm}}(\max\{k,\ell+1\})^2$; 
        and \label{item:proof_connecting_lemma_Ziy_ii}
        \item $Z_{i,y}$ contains the adhesions of $y$ with its neighbors in $T[x_1,x_2]$,
            that is $W_{y} \cap W_{y'} \subseteq Z_{i,y}$ for every $y' \in N_{T[x_1,x_2]}(y)$. \label{item:proof_connecting_lemma_Ziy_iii}
    \end{enumerate}
    Let $A_1$ be the adhesion of $W_y$ with its predecessor in $T[x_1,x_2]$ if $y \neq x_1$, 
    and let $A_1 = \{r_i\}$ if $y=x_1$.
    Similarly, let $A_2$ be the adhesion of $W_y$ with its successor in $T[x_1,x_2]$ if $y \neq x_2$, 
    and let $A_2 = \{r_{i+1}\}$ if $y=x_2$.
    Let $Q$ be a path from $A_1$ to $A_2$ in $\torso_G(U, \mathcal{B})$.
    Such a path exists by \ref{item:lemma:connecting_a_small_set:thm:k_unbreakable_treedecomposition:i}.

    If $|V(Q) \cap W_y| < 2k$, 
    then let 
    \[
        Z_{i,y} = A_1 \cup A_2 \cup \Pi_{\torso_G(U, \mathcal{B})}(W_y, V(Q)).
    \]
    We now quickly justify that $Z_{i,y}$ is small.
    Let $P_1, \dots, P_m$ be a partition of $Q$ into $m = |V(Q) \cap W_y|-1$ vertex disjoint paths which are all internally disjoint from $W_y$.
    Then, for every $j \in [m]$, either $P_j$ has only one vertex,
    or there is a component of $\torso_G(U,\mathcal{B}) - W_y$ containing the interior
    of $P_j$. In the latter case, there is a neighbor $y'$ of $y$ such that
    $\Pi_{\torso_G(U,\mathcal{B})}(W_y, V(P_j)) \subseteq W_y \cap W_{y'}$.
    In both cases, we get that
    $|\Pi_{\torso_G(U,\mathcal{B})}(W_y, V(P_j))| \leq \max\{1,k-1\} \leq k$.
    Summing over all $j \in [m]$, we deduce
    $|\Pi_{\torso_G(U,\mathcal{B})}(W_y,V(Q))| \leq (|V(Q) \cap W_y|-1) \cdot k$ and
    so 
    \[
         |Z_{i,y}| \leq 2\max\{1,k-1\} + (2k-2)k \leq 2k^2.
    \]
    Moreover, for every connected component $C$ of $\torso_G(U, \mathcal{B})-Z_{i,y}$ intersecting $W_y$,
    $V(C)$ is disjoint from $V(Q)$.
    Hence \ref{item:proof_connecting_lemma_Ziy_i}, \ref{item:proof_connecting_lemma_Ziy_ii},
    and \ref{item:proof_connecting_lemma_Ziy_iii} hold.

    Now suppose that $|V(Q) \cap W_y| \geq 2k$.
    Let $B_1$ be the vertex set of the shortest prefix of $Q$ containing $k$ vertices of $W_y$, and let $B_2$ be the vertex set of the shortest suffix of $Q$ containing $k$ vertices in $W_y$.
    By \ref{item:lemma:connecting_a_small_set:thm:k_unbreakable_treedecomposition:i}, there are $k$ pairwise vertex-disjoint $(B_1 \cap W_y, B_2 \cap W_y)$-paths in $\torso_G(U, \mathcal{B})$.
    By Lemma~\ref{lemma:finding_k_paths_between_sets_connected_in_the_torso_of_a_nice_set}, there are $k$ pairwise vertex-disjoint
    $(B_1,B_2)$-paths $Q_1, \dots, Q_k$ in $G$.
    Let $\mathcal{F}$ be the family of all the connected subgraphs $H$ of 
    $G$ such that $V(H) \cap Q_j \neq \emptyset$ 
    for every $j \in [k]$.
    By Lemma~\ref{lemma:finding_kladder_from_a_mess}, there are no $f_{\cstref{lemma:finding_kladder_from_a_mess}}(k,\ell)$ pairwise vertex-disjoint
    members of $\mathcal{F}$.
    Recall that $G$ has no $k$-ladder of size $\ell$,
    and so has no $\max\{k,\ell+1\} \times \max\{k,\ell+1\}$ grid as a minor.
    By Theorem~\ref{thm:grid_minor_thm}, $\tw(G) < f_{\cstref{thm:grid_minor_thm}}(\max\{k,\ell+1\})$.
    Hence by Lemma~\ref{lemma:helly_property} applied to the family $\mathcal{F}$ and
    a tree decomposition of $G$ of width less than $f_{\cstref{thm:grid_minor_thm}}(\max\{k,\ell+1\})$,
    we deduce that there is a set $Z^0_{i,y} \subseteq V(G)$ of size at most 
    $f_{\cstref{lemma:finding_kladder_from_a_mess}}(k,\ell) \cdot f_{\cstref{thm:grid_minor_thm}}(\max\{k,\ell+1\})$
    that intersects every member of $\mathcal{F}$.
    Then let 
    \[
        Z_{i,y} = A_1 \cup A_2 \cup \Pi_{\torso_G(U,\mathcal{B})}(W_y, \Pi_G(U, Z^0_{i,y}) \cup B_1 \cup B_2).
    \]
    See Figure~\ref{fig:connecting_lemma}.
        
    \begin{figure}[ht]
        \centering
        \includegraphics{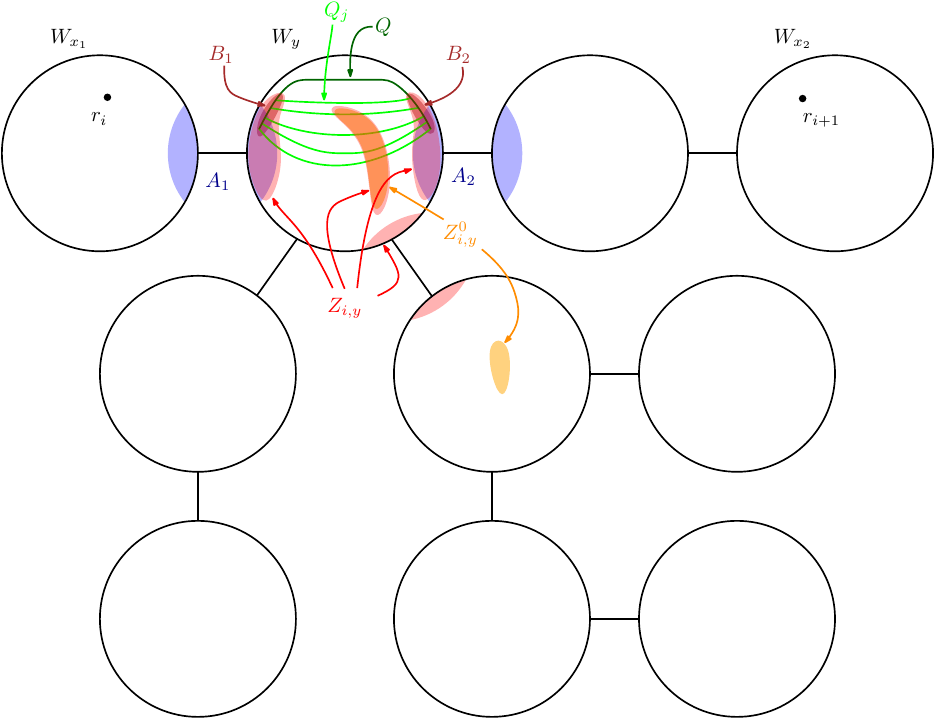}
        \caption{Illustration for the proof of Lemma~\ref{lemma:connecting_a_small_set}. 
            For sake of clarity, we assume here that $U=V(G)$, and $Q$, $B_1$, and $B_2$ are drawn inside $W_y$,
            but in general this is not the case, and so $Z_{i,y}$ contains only the projection of $B_1 \cup B_2$ on $W_y$.}
        \label{fig:connecting_lemma}
    \end{figure}

    Note that \ref{item:proof_connecting_lemma_Ziy_iii} holds by construction, and $Z_{i,y} \subseteq W_y$.
    We now argue that $Z_{i,y}$ is small.
    First, $|A_1 \cup A_2| \leq 2 \max\{1, k-1\} \leq 2k$.
    We now show an upper bound for the size of $\Pi_{\torso_G(U,\mathcal{B})}(W_y,\Pi_G(U,Z^0_{i,y}))$.
    Recall that $|Z^0_{i,y}| \leq f_{\cstref{lemma:finding_kladder_from_a_mess}}(k,\ell) \cdot f_{\cstref{thm:grid_minor_thm}}(\max\{k,\ell+1\})$.
    Let $u \in Z^0_{i,y}$.
    The set $\Pi_G(U, u)$ induces a clique in $\torso_G(U,\mathcal{B})$.
    Since $\big(T,(W_x \mid x \in V(T))\big)$ is a tree decomposition of $\torso_G(U,\mathcal{B})$,
    there exists $z \in V(T)$ such that $\Pi_G(U,u) \subseteq W_z$.
    If $z \neq y$, we deduce that 
    \[
        \Pi_{\torso_G(U,\mathcal{B})}(W_y, \Pi_G(U,u)) \subseteq W_y \cap W_z,
    \]
    and so $|\Pi_{\torso_G(U,\mathcal{B})}(W_y,\Pi_G(U,u))| \leq k-1$.
    Otherwise, if $z=y$, then 
    \[
        \Pi_{\torso_G(U,\mathcal{B})}(W_y, \Pi_G(U,u)) = \Pi_G(U,u),
    \]
    which is a clique in $\torso_G(U,\mathcal{B})$,
    and so size at most $1+\tw(\torso_G(U,\mathcal{B}))\leq 2f_{\cstref{thm:grid_minor_thm}}(\max\{k,\ell+1\})$ by \eqref{eq:tw_torsoUB}.
    In both cases, 
    we have $|\Pi_{\torso_G(U,\mathcal{B})}(W_y,\Pi_G(U,u))| \leq \max\{k-1,2 f_{\cstref{thm:grid_minor_thm}}(\max\{k,\ell+1\})\} = 2f_{\cstref{thm:grid_minor_thm}}(\max\{k,\ell+1\})$.
    Therefore,
    \begin{align*}
        |\Pi_{\torso_G(U,\mathcal{B})}(W_y, \Pi_G(U,Z^0_{i,y}))| 
        &\leq 2f_{\cstref{thm:grid_minor_thm}}(\max\{k,\ell+1\}) \cdot |Z^0_{i,y}| \\
        &\leq 2 f_{\cstref{lemma:finding_kladder_from_a_mess}}(k,\ell) \cdot f_{\cstref{thm:grid_minor_thm}}(\max\{k,\ell+1\})^2.
    \end{align*}
    Moreover,
    $B_1 \cup B_2$ can be covered by the vertex sets of $2k$ subpaths of $Q$
    each of them internally disjoint from $W_y$ 
    and intersecting $W_y$ in at most one vertex.
    Each of these subpath $Q'$ has a projection
    $\Pi_{\torso_G(U,\mathcal{B})}(W_y, V(Q'))$
    of size at most $\max\{1,k-1\}\leq k$
    since $\big(T,(W_x \mid x \in V(T))\big)$ has adhesion at most $k-1$.
    We deduce that $|\Pi_{\torso_G(U,\mathcal{B})}(W_y, B_1 \cup B_2)| \leq 2k^2$.
    Altogether, this shows 
    \begin{align*}
        |Z_{i,y}| 
        &\leq |A_1 \cup A_2| + |\Pi_{\torso_G(U,\mathcal{B})}(W_y, \Pi_G(U,Z^0_{i,y}))| + |\Pi_{\torso_G(U,\mathcal{B})}(W_y, B_1 \cup B_2)| \\
        &\leq 2k + 2f_{\cstref{lemma:finding_kladder_from_a_mess}}(k,\ell) \cdot f_{\cstref{thm:grid_minor_thm}}(\max\{k,\ell+1\})^2
        + 2k^2 \\
        &= 2(k+1)k + 2f_{\cstref{lemma:finding_kladder_from_a_mess}}(k,\ell) \cdot f_{\cstref{thm:grid_minor_thm}}(\max\{k,\ell+1\})^2,
    \end{align*}
    which proves \ref{item:proof_connecting_lemma_Ziy_ii}.
    Consider now a connected component $C$ of $\torso_G(U, \mathcal{B})-Z_{i,y}$ intersecting $W_y$.
    First, 
    since $C$ is disjoint from $\Pi_{\torso_G(U,\mathcal{B})}(W_y, B_1 \cup B_1)$,
    $C$ is disjoint from $B_1 \cup B_2$.
    Second, since $C$ is disjoint from $\Pi_{\torso_G(U,\mathcal{B})}(W_y, \Pi_G(U,Z^0_{i,y}))$, 
    $C$ is disjoint from $\Pi_G(U,Z^0_{i,y})$.
    Since $(U, \mathcal{B})$ is a nice pair
    and by Lemma~\ref{lemma:paths_in_G_implies_paths_in_torso} applied for $k'=1$,
    for every edge $uv \in E(C) \setminus E(G)$, 
    there is a $(u,v)$-path in $G$ internally disjoint from $U$,
    and it follows that there is a connected subgraph $C'$ of $G$
    such that $V(C) = V(C') \cap U$.
    In particular, $V(C')$ is disjoint from $\Pi_G(U, Z^0_{i,y})$ 
    and intersects $W_y$. 
    Therefore, $V(C')$ is disjoint from $Z^0_{i,y}$, 
    and so is not a member of $\mathcal{F}$.
    As a consequence,
    there is $j \in [k]$ such that $V(Q_j)$ is disjoint from $V(C')$,
    and so from $V(C) \subseteq V(C')$.
    Fix such a $j \in [k]$.
    Then, 
    $V(Q_j) \cap U$ induces a connected subgraph of $\torso_G(U,\mathcal{B}) - V(C)$
    which intersects $B_1$ and $B_2$.
    Therefore, we obtain a path from $A_1$ to $A_2$ in $\torso_G(U, \mathcal{B}) - V(C)$.
    If $y = x_1$, then recall that $A_1 = \{r_i\}$.
    If $y \neq x_1$, let $y'$ the predecessor of $y$ along $T[x_1,x_2]$.
    Then $\bigcup_{z \in V(T_{y' \mid y})}W_z$ induces a connected subgraph in $\torso_G(U, \mathcal{B})$, and this subgraph is disjoint from $V(C)$ and contains $r_i$.
    In both cases, there is path from $r_i$ to every vertex in $A_1$ in $\torso_G(U, \mathcal{B})-V(C)$.
    Similarly, there is a path between $r_{i+1}$ and every vertex in $A_2$ in $\torso_G(U, \mathcal{B})-V(C)$.
    Altogether, we deduce that there is a path from $r_i$ to $r_{i+1}$
    in $\torso_G(U, \mathcal{B})-V(C)$.
    This proves \ref{item:proof_connecting_lemma_Ziy_i}.

    \medskip

    Now let 
    \[
        Z_i = \textstyle\bigcup_{y \in V(T[x_1,x_2])} Z_{i,y}.
    \]
    By \ref{item:proof_connecting_lemma_Ziy_iii},
    $\big(T[x_1,x_2], (Z_{i,y} \mid y \in V(T[x_1,x_2]))\big)$ 
    is a tree decomposition of $(\torso_G(U,\mathcal{B}),Z_i)$
    of adhesion at most $k-1$ and of width less than
    $2(k+1)k + 2f_{\cstref{lemma:finding_kladder_from_a_mess}}(k,\ell) \cdot f_{\cstref{thm:grid_minor_thm}}(\max\{k,\ell+1\})^2$.
    In particular, it is a $k$-dismantable tree decomposition,
    and therefore 
    \begin{align*}
        \td_k(\torso_G(U, \mathcal{B}), Z_i) 
        \leq 2(k+1)k + 2f_{\cstref{lemma:finding_kladder_from_a_mess}}(k,\ell) \cdot f_{\cstref{thm:grid_minor_thm}}(\max\{k,\ell+1\})^2,
    \end{align*}
    and so \ref{item:proof_connecting_lemma_Zi_ii} holds.
    We now show \ref{item:proof_connecting_lemma_Zi_i}.
    Let $C$ be a connected component of $\torso_G(U, \mathcal{B})-Z_i$.
    If $V(C)$ intersects $W_y$ for some $y \in V(T[x_1,x_2])$, then
    by \ref{item:proof_connecting_lemma_Ziy_i},
    there is a path from $r_i$ to $r_{i+1}$ in $\torso_G(U, \mathcal{B})-V(C)$.
    Otherwise, there is an edge $yy' \in E(T)$ with $y \in V(T[x_1,x_2])$ and $y' \not\in V(T[x_1,x_2])$
    such that $V(C) \subseteq \bigcup_{z \in V(T_{y' \mid y})} W_z \setminus W_{y}$.
    Then, the set $\bigcup_{z \in V(T_{y \mid y'})} W_z$ induces a connected subgraph in $\torso_G(U, \mathcal{B})$
    which contains $r_i$ and $r_{i+1}$, and which is disjoint from $V(C)$.
    Therefore, there is a path from $r_i$ to $r_{i+1}$ in $\torso_G(U, \mathcal{B})-V(C)$.
    This proves \ref{item:proof_connecting_lemma_Zi_i}.

    \medskip

    Finally, let
    \[
        Z = \textstyle\bigcup_{i \in [|R|-1]} Z_i.
    \]
    Then,
    for every connected component $C$ of $\torso_G(U, \mathcal{B}) - Z$,
    for every $i \in [|R|-1]$, 
    $C$ is disjoint from $Z_i$ and so by \ref{item:proof_connecting_lemma_Zi_i},
    there is a path from $r_i$ to $r_{i+1}$ in $\torso_G(U, \mathcal{B}) - V(C)$.
    Therefore, $R$ is connected in $\torso_G(U, \mathcal{B}) - V(C)$, which proves \ref{item:connecting_lemma_i}.
    Moreover, by Lemma~\ref{lemma:combine_tree_decompositions_of_S},
    \begin{align*}
        \td_k(\torso_G(U, \mathcal{B}),Z) 
        &\leq (|R|-1) \cdot \left(2(k+1)k + 2f_{\cstref{lemma:finding_kladder_from_a_mess}}(k,\ell) \cdot f_{\cstref{thm:grid_minor_thm}}(\max\{k,\ell+1\})^2\right)\\
        &\leq f_{\cstref{lemma:connecting_a_small_set}}(k,\ell,c),
    \end{align*}
    and so \ref{item:connecting_lemma_ii} holds.
    This proves the lemma.
\end{proof}

We can now give the proof of Theorem~\ref{thm:main_induction}.

\begin{proof}[Proof of Theorem~\ref{thm:main_induction}]
    Let $k,\ell$ be fixed positive integers.
    We proceed by induction on $a$.
    If $a \leq k-1$, then $(\torso_G(U, \mathcal{B}),S)$ has a tree decomposition of width less than $t$ and adhesion at most $k-1$, 
    and so $\td_k(\torso_G(U, \mathcal{B}),S) \leq t \leq f_{\cstref{thm:main_induction}}(k, \ell, t, a)$ for $f_{\cstref{thm:main_induction}}(k, \ell, t, a) = t$.
    Now suppose $a \geq k$, and that the result holds for $a-1$.

    First, we prove the following claim.
    \begin{claim}\label{claim:thm:main_induction}
        There is a function $f_{\cstref{claim:thm:main_induction}} \colon \NN^4 \to \NN$ such that the following holds.
        Let $\ell',c,c'$ be positive integers.
        Let $G$ be a
        graph that does not contain any $k$-ladder of length $\ell$ as a minor.
        Let $(U, \mathcal{B})$ be a nice pair in $G$.
        Let $S \subseteq U$ be such that there is a tree decomposition of $(\torso_G(U, \mathcal{B}),S)$ of adhesion at most $a$ and width less than $t$.
        Let $(V_0, \dots, V_{\ell'})$ be a path partition of $\torso_G(U, \mathcal{B})$ such that
        \begin{enumerate}[label=(\alph*)]
            \item $N_{\torso_G(U, \mathcal{B})}(V_{i}) \cap V_{i-1}$ is connected in 
                $\torso_G(U, \mathcal{B})[(N_{\torso_G(U, \mathcal{B})}(V_{i}) \cap V_{i-1}) \cup V_i]$, for every $i \in [\ell'-1]$;\label{item:claim:main_induction:i}
            \item $|N_{\torso_G(U, \mathcal{B})}(V_{\ell'})| \leq c$; and \label{item:claim:main_induction:ii}
            \item for every $W \subseteq V_{\ell'}$, either \label{item:claim:main_induction:iii}
                \begin{enumerate}[label=(\roman*)]
                    \item there are $k$ pairwise disjoint $(V_0,W)$-paths in $\torso_G(U, \mathcal{B})$, or
                    \item $\td_k(\torso_G(U, \mathcal{B})[W], S \cap W) \leq c'$.
                \end{enumerate}
        \end{enumerate}
        Then there exists $S' \subseteq V_{\ell'}$ with $S \cap V_{\ell'} \subseteq S'$ such that
        \[
        \td_k(\torso_G(U, \mathcal{B})[V_{\ell'}], S') \leq f_{\cstref{claim:thm:main_induction}}(t,\ell', c, c').
        \]
    \end{claim}

    Before proving it, we motivate the statement of this claim.
    Its proof will consists in a downward induction on $\ell'$.
    The path partition $(V_0, \dots, V_{\ell'})$ will ensure that the induction will eventually stop: if this path partition is long enough, 
    and if we have $k$ pairwise disjoint paths from $V_0$ to $V_{\ell'}$,
    then we will find a long $k$-ladder by Lemma~\ref{lemma:finding_kladder_from_a_mess}, since each part is morally connected by \ref{item:claim:main_induction:i}.
    The problem with this approach is that there are graphs, even with large $k$-treedepth, that have no path partition with more than $3$ parts.
    This is for example the case when there is a universal vertex.
    But in this specific case, removing the universal vertex will decrease the treewidth by one.
    In general, we will use the fact mentioned in Section~\ref{sec:outline} that for every graph $G$, for every $u \in V(G)$,
    there exists $S \subseteq V(G) \setminus \{u\}$ containing $N_G(u)$ with $\tw(G-u,S)<\tw(G)$.
    Hence, by induction, we can decompose a superset of $N_G(u)$ in $G-u$,
    which will enable us to construct a path partition.
    This is why we are working in the setting of decomposition ``focused'' on a subset $S$ of vertices.

    \medskip

    \begin{proofclaim}
        Let $L = f_{\cstref{lemma:finding_kladder_from_a_mess}}(k,\ell)$, where $f_{\cstref{lemma:finding_kladder_from_a_mess}}$ 
        is the function from Lemma~\ref{lemma:finding_kladder_from_a_mess}.
        If $\ell'\geq 2L$, then let
        \begin{align*}
            f_{\cstref{claim:thm:main_induction}} (t, \ell',c,c')
            &= c'. \\
        \intertext{Otherwise, if $\ell' \leq 2L-1$, let}
            f_{\cstref{claim:thm:main_induction}} (t, \ell',c,c')
            &= f_{\cstref{lemma:connecting_a_small_set}}(k,\ell,c) \\
            & \hspace{1cm}
                + c \cdot t \\
            & \hspace{1cm}
                + c \cdot f_{\cstref{thm:main_induction}}(k, \ell, t-1, a-1)  \\
            & \hspace{1cm}
                + f_{\cstref{claim:thm:main_induction}}\big(t, \ell'+1, 
                f_{\cstref{lemma:connecting_a_small_set}}(k,\ell,c)
                    + c \cdot t 
                    + c \cdot f_{\cstref{thm:main_induction}}(k, \ell, t-1, a-1), c'\big).
        \end{align*}

        We proceed by induction on $(\max\{0,2L-\ell'\},|V_{\ell'}|)$, in the lexicographic order.
        Here, $\ell'$ will play the role of a variable, and so we will apply the induction hypothesis on instances where $\ell'$ is larger, or equal but with $|V_{\ell'}|$ smaller.
        Let $\mathcal{D}_0 = \big(T_0, (W^0_x \mid x \in V(T_0))\big)$ be a tree decomposition of $(\torso_G(U, \mathcal{B}),S)$ of adhesion at most $a$ and width less than $t$.
        
        First suppose $\ell' \geq 2L$.
        By~\ref{item:claim:main_induction:iii} applied to $W=V_{\ell'}$, either
        \begin{enumerate}[label=(\roman*)]
            \item there are $k$ pairwise disjoint $(V_0,V_{\ell'})$-paths in $\torso_G(U, \mathcal{B})$, or
            \item  $\td_k(\torso_G(U, \mathcal{B})[V_{\ell'}], S \cap V_{\ell'}) \leq c'$.
        \end{enumerate}
        In the second case,
        we set $S' = S \cap V_{\ell'}$ and we are done. 
        Now assume that there are $k$ pairwise disjoint $(V_0,V_{\ell'})$-paths in $\torso_G(U, \mathcal{B})$.
        Since $(V_0, \dots, V_{\ell'})$ is a path partition of $\torso_G(U,\mathcal{B})$, each of these paths intersects $V_i$ for every $i \in \{0, \dots, \ell'\}$.
        Therefore, there are $k$ pairwise disjoint $(V_0,V_{2L})$-paths in $\torso_G(U, \mathcal{B})$.
        By Lemma~\ref{lemma:finding_k_paths_between_sets_connected_in_the_torso_of_a_nice_set},
        there are $k$ pairwise disjoint $(V_0,V_{2L})$-paths $Q_1, \dots, Q_k$ in $G$.
        For every $i \in [2L-1]$, $N_{\torso_G(U, \mathcal{B})}(V_{i}) \cap V_{i-1}$ is connected
        in $\torso_G(U, \mathcal{B})[(N_{\torso_G(U, \mathcal{B})}(V_{i}) \cap V_{i-1}) \cup V_i]$ by~\ref{item:claim:main_induction:i}.
        Let $C_i$ be the vertex set of the connected component of $\torso_G(U, \mathcal{B})[(N_G(V_{i}) \cap V_{i-1}) \cup V_i]$ containing $N_G(V_{i}) \cap V_{i-1}$.
        Let $\mathcal{C}_i$ be the family of all the members $B$ of $\mathcal{B}$ with $B \cap V(C_i) \neq \emptyset$.
        Note that $\bigcup \mathcal{C}_i$ and $\bigcup \mathcal{C}_j$ are pairwise disjoint for every $i,j \in [2L-1]$ with $|i-j|\geq 2$, 
        because $B \cap U$ induces a clique in $\torso_G(U, \mathcal{B})$ for every $B \in \mathcal{B}$.
        Moreover, $N_{\torso_G(U,\mathcal{B})}(V_{i}) \cap V_{i-1}$ is connected in $G[C_i \cup \bigcup \mathcal{C}_i]$.
        Let $C'_i$ be the connected component of $G[C_i \cup \bigcup \mathcal{C}_i]$ containing $N_{\torso_G(U, \mathcal{B})}(V_i) \cap V_{i-1}$.
        Then $C'_1, C'_3, \dots, C'_{2L-1}$ is a family of $L$ pairwise disjoint connected subgraphs of $G$.
        We claim that each of these subgraphs intersects $V(Q_j)$ for every $j \in [k]$.
        Indeed, by Lemma~\ref{lemma:paths_in_G_implies_paths_in_torso}, 
        for all $i \in \{1,3,\dots, 2L-1\}, j \in [k]$,
        $\torso_G(U,\mathcal{B})[V(Q_j) \cap U]$
        contains a $(V_0,V_{2L})$-path.
        Since $(V_0, \dots, V_{\ell'})$ is a path partition of $\torso_G(U,\mathcal{B})$, we deduce that $Q_j$ intersects
        $N_{\torso_G(U,\mathcal{B})}(V_i) \cap V_{i-1} \subseteq V(C'_i)$ for every $j \in [k]$.
        Therefore, by Lemma~\ref{lemma:finding_kladder_from_a_mess}
        and the definition of $L$,
        we have that $G$ contains a $k$-ladder of length $\ell$ as a minor,
        a contradiction.

        \medskip
        
        Now suppose $\ell' \leq 2L-1$.
        If $\torso_G(U, \mathcal{B})[V_{\ell'}]$ is not connected, then let $\mathcal{C}$ be the family of all the connected components of $\torso_G(U,\mathcal{B})[V_{\ell'}]$.
        For every $C \in \mathcal{C}$, we call induction on the path partition $\big(V_0, \dots, V_{\ell'-2}, V_{\ell'-1} \cup (V_{\ell'} \setminus V(C)), V(C)\big)$.
        We deduce that there exists $S'_C \subseteq V(C)$ containing $S \cap V(C)$ such that 
        $\td_k(C,S'_C) \leq f_{\cstref{claim:thm:main_induction}}(t, \ell', c, c')$.
        Then for $S' = \bigcup_{C \in \mathcal{C}} S'_C$, 
        we have $S \cap V_{\ell'} \subseteq S' \subseteq V_{\ell'}$ and 
        $\td_k(\torso_G(U, \mathcal{B})[V_{\ell'}], S') = \max_{C \in \mathcal{C}} \td_k(C, S' \cap V(C)) \leq f_{\cstref{claim:thm:main_induction}}(t, \ell', c, c')$.
        Now assume that $\torso_G(U, \mathcal{B})[V_{\ell'}]$ is connected.

        \medskip

        \begin{figure}[p]
            \centering
            \includegraphics{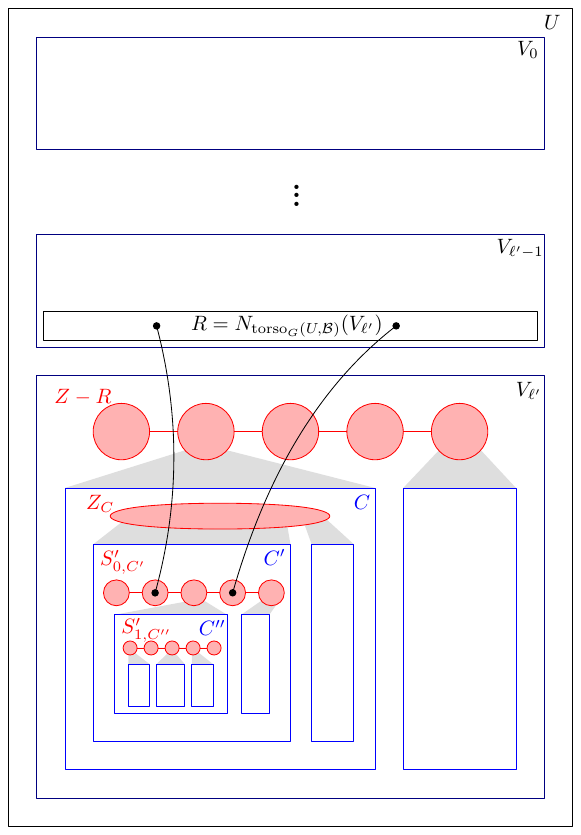}
            \caption{Illustration for the proof of Claim~\ref{claim:thm:main_induction}.
                First we remove the set $Z$ given by Lemma~\ref{lemma:connecting_a_small_set} to ensure that $R$ is connected 
                in $\torso_{G'}(U',\mathcal{B}')-V(C)$ for every connected component $C$ of $\torso_{G'}(U',\mathcal{B}')-Z$.
                Fix a connected component $C$ of $\torso_{G'}(U',\mathcal{B}')$.
                Then, we remove $Z_C$ to ensure that connected components of $C-Z_C$ are either disjoint from $S$ (and so we are done), or have their neighborhood included in $S$.
                Fix a connected component $C'$ of $C-Z_C$ which intersects $S$.
                We want to call induction on $\ell'$ by considering the partition $\big(V_0, \dots, V_{\ell'-1}, V_{\ell'} \setminus V(C'), V(C')\big)$.
                However, this might be not a path partition, since $V(C')$ can have neighbors in $R$.
                Hence we need to decompose $N_{\torso_{G'}(U', \mathcal{B}')}(u) \cap V(C')$ for every $u \in R$.
                This is done by calling induction on $a$, the adhesion of the given tree decomposition of $(\torso_G(U,\mathcal{B}),S)$.
                Now that the neighborhood of $R$ is covered by a set $S'_{0,C'}$, we can indeed call induction on $\ell'$,
                which yields $S'_{1,C''}$ for each component $C''$ of $C' - S'_{0,C'}$ intersecting $S$.
                Then, combining the $k$-dismantable tree decompositions obtained at each step, we deduce the claim.}
            \label{fig:main_claim}
        \end{figure}
        
        A general picture of the following proof is provided in Figure~\ref{fig:main_claim}.
        Let $X = V_0 \cup \dots \cup V_{\ell'-2} \cup (V_{\ell'-1} \setminus N_{\torso_G(U, \mathcal{B})}(V_{\ell'}))$,
        and let $U' = U \setminus X$.
        Note that $U' = N_{\torso_G(U, \mathcal{B})}(V_{\ell'}) \cup V_{\ell'}$.
        By Lemma~\ref{lemma:U_nice_impliesU-X_nice}, 
        there exists a nice pair $(U',\mathcal{B}'_0)$ in $G-X$ with 
        \[
            \torso_{G-X}(U',\mathcal{B}'_0) = \torso_G(U,\mathcal{B}) - X = \torso_G(U,\mathcal{B})[N_{\torso_G(U, \mathcal{B})}(V_{\ell'}) \cup V_{\ell'}].
        \]
        Moreover, $U'$ induces a connected subgraph of $\torso_{G}(U, \mathcal{B})$ disjoint from $X$,
        and so by Lemma~\ref{lemma:finding_k_paths_between_sets_connected_in_the_torso_of_a_nice_set}, $U'$ is connected in $G-X$.
        Let $G'$ be the connected component of $G-X$ containing $U'$.
        Let $\mathcal{B}' = \{B \cap V(G') \mid B \in \mathcal{B}'_0\} \setminus \{\emptyset\}$.
        Then, $(U', \mathcal{B}')$ is a nice pair in $G'$ and $\torso_{G'}(U', \mathcal{B}') = \torso_{G-X}(U',\mathcal{B}'_0)$.
        By Lemma~\ref{lemma:connecting_a_small_set} applied to 
        \[
            R = N_{\torso_G(U, \mathcal{B})}(V_{\ell'}),
        \]
        there is a set $Z \subseteq U'$ containing $R$ such that     
        \begin{enumerate}[label=\ref{lemma:connecting_a_small_set}.(\alph*)]
            \item for every connected component $C$ of $\torso_{G'}(U', \mathcal{B}')-Z$, the set $R$ is connected in $\torso_{G'}(U', \mathcal{B}')-V(C)$; and
            \item $\td_k(\torso_{G'}(U', \mathcal{B}'),Z) \leq f_{\cstref{lemma:connecting_a_small_set}}(k,\ell,c)$.
        \end{enumerate}
 
        Recall that $\torso_{G'}(U', \mathcal{B}')-Z = \torso_{G}(U, \mathcal{B})-(X \cup Z)$.
        Fix a connected component $C$ of $\torso_{G}(U, \mathcal{B})- (X \cup Z)$.
        For every $u \in (N_{\torso_{G}(U, \mathcal{B})}(V(C)) \cap R) \setminus S$, 
        let $C_u$ be the connected component of $u$ in $\torso_G(U, \mathcal{B})-S$.
        Since $\mathcal{D}_0$ is a tree decomposition of $(\torso_G(U, \mathcal{B}),S)$,
        there exists $x_u \in V(T_0)$ such that $N_{\torso_G(U, \mathcal{B})}(V(C_u)) \subseteq W^0_{x_u}$.
        Now let 
        \[
            Z_C = \bigcup_{u \in \big(N_{\torso_G(U, \mathcal{B})}(V(C)) \cap R\big) \setminus S} \left(W^0_{x_u} \cap V(C)\right).
        \]
        Note that $Z_C \subseteq V(C) \cap S$ and $|Z_C| \leq c \cdot t$.
        Now, for every connected component $C'$ of $C-Z_C$, either $V(C') \cap S = \emptyset$, or $N_{\torso_G(U, \mathcal{B})}(V(C')) \cap R \subseteq S$.
        Let $\mathcal{C}_C$ be the family of all the connected components of $C-Z_C$ intersecting $S$.
        Fix $C' \in \mathcal{C}_C$, and consider some
        vertex $u \in N_{\torso_{G}(U, \mathcal{B})}(V(C')) \cap R$.
        Recall that $u \in S$.
        Let $T_{0,u}$ be the subtree of $T_0$ induced by $\{x \in V(T_0) \mid u \in W^0_x\}$,
        and let $S_u = \bigcup_{x \in V(T_{0,u})} W^0_x \cap V(C')$.
        Then $\big(T_{0,u}, (W^0_x \cap V(C') \mid x \in V(T_{0,u}))\big)$ is a tree decomposition of $(C', S_u)$
        of width less than $t-1$ and adhesion at most $a-1$ since $u \in W^0_x \setminus V(C')$ for every $x \in V(T_{0,u})$.
        Moreover, by Lemma~\ref{lemma:U_nice_impliesU-X_nice}, 
        since $V(C') = U \setminus (U \setminus V(C'))$,
        there is a nice pair in $G-(U \setminus V(C'))$ of the form $(V(C'), \mathcal{B}_{C'})$ such that
        $C' = \torso_G(U,\mathcal{B}) - (U \setminus V(C')) = \torso_G(V(C'), \mathcal{B}_{C'})$.
        Hence, by the induction hypothesis applied to $G-(U \setminus V(C')), V(C'), \mathcal{B}_{C'}, S_u$, there exists $S'_u \subseteq V(C')$ with $S_u \subseteq S'_u$ such that
        \[
        \td_k(C', S'_u) \leq \td_k\big(\torso_{G-(U \setminus V(C'))}(V(C'), \mathcal{B}_{C'}), S'_u\big) \leq f_{\cstref{thm:main_induction}}(k, \ell, t-1, a-1).
        \]       
        Let 
        \[
            S'_{0,C'} = \bigcup_{u \in N_{\torso_{G}(U, \mathcal{B})}(V(C')) \cap R} S'_u.
        \]
        Then, by Lemma~\ref{lemma:combine_tree_decompositions_of_S},
        \[
            \td_k(C', S'_{0,C'}) \leq c \cdot f_{\cstref{thm:main_induction}}(k, \ell, t-1, a-1).
        \]
        Consider now a connected component $C''$ of $C' - S'_{0,C'}$ with $V(C'') \cap S \neq \emptyset$.
        Since $C''$ is disjoint from $\bigcup_{u \in N_{\torso_{G'}(U', \mathcal{B}')}(V(C')) \cap R} S_u$, this implies that there is no edge between $R$ and $V(C'')$ in $\torso_G(U, \mathcal{B})$.
        Moreover,
        \begin{align*}
            |N_{\torso_G(U, \mathcal{B})}(V(C''))| 
            &\leq |N_{\torso_G(U, \mathcal{B})}(V(C'')) \cap Z| + |N_{\torso_G(U, \mathcal{B})}(V(C'')) \cap Z_C| \\
            & \qquad\qquad + |N_{\torso_G(U, \mathcal{B})}(V(C'')) \cap S'_{0,C'}| \\
            &\leq f_{\cstref{lemma:connecting_a_small_set}}(k,\ell,c) + c \cdot t + c \cdot f_{\cstref{thm:main_induction}}(k, \ell, t-1, a-1) \\
        \end{align*}

        Then $\big(V_0, \dots, V_{\ell'-1}, V_{\ell'} \setminus V(C''), V(C'')\big)$ is a path partition of $\torso_G(U, \mathcal{B})$.
        By the properties of $Z$, $N_{\torso_G(U,\mathcal{B})}(V_{\ell'} \setminus V(C'')) \cap V_{\ell'-1} = N_{\torso_G(U,\mathcal{B})}(V_{\ell'}) = R$ is connected in 
        $\torso_{G'}(U', \mathcal{B}') - V(C'') = \torso_G(U, \mathcal{B}) - X - V(C'') = \torso_G(U,\mathcal{B})\big[(N_{\torso_G(U,\mathcal{B})}(V_{\ell'} \setminus V(C'')) \cap V_{\ell'-1}) \cup (V_{\ell'} \setminus V(C''))\big]$.
        Hence $\big(V_0, \dots, V_{\ell'-1}, V_{\ell'} \setminus V(C''), V(C'')\big)$ satisfies the assumptions of the claim
        for the parameters $(\ell'+1, f_{\cstref{lemma:connecting_a_small_set}}(k,\ell,c) + c \cdot t + f_{\cstref{thm:main_induction}}(k, \ell, t-1, a-1), c')$.
        Hence, by induction on $\ell'$,
        there exists $S'_{1,C''} \subseteq V(C'')$ with $S \cap V(C'') \subseteq S'_{1,C''}$ such that
        \[
        \td_k(C'', S'_{1,C''}) 
              \leq f_{\cstref{claim:thm:main_induction}}\big(t, \ell'+1, f_{\cstref{lemma:connecting_a_small_set}}(k,\ell,c)
                  + c \cdot t 
                  + c \cdot f_{\cstref{thm:main_induction}}(k, \ell, t-1, a-1), c'\big).
        \]

        Now, let 
        \[
            S'_{C'} = S'_{0,C'} \cup \bigcup_{\substack{\text{$C''$ connected component of $C'-S'_{0,C}$} \\\text{$V(C'') \cap S \neq \emptyset$}}} S'_{1,C''}.
        \]
        Then, by Lemma~\ref{lemma:combine_tree_decompositions_of_S},
        \begin{align*}
            \td_k(C', S'_{C'}) 
            &\leq \td_k(C', S'_{0,C'}) 
            + \max_{\substack{\text{$C''$ connected component of $C'-S'_{0,C}$} \\\text{$V(C'') \cap S \neq \emptyset$}}} \td_k(C'',S'_{1,C''}) \\
            &\leq c \cdot f_{\cstref{thm:main_induction}}(k, \ell, t-1, a-1) +\\
            & \hspace{1cm}
                + f_{\cstref{claim:thm:main_induction}}\big(t, \ell'+1, \\
            & \hspace{1cm} \hspace{1.5cm} f_{\cstref{lemma:connecting_a_small_set}}(k,\ell,c) 
                + c \cdot t 
                + c \cdot f_{\cstref{thm:main_induction}}(k, \ell, t-1, a-1),c'\big).
        \end{align*}

        Now let 
        \[
            S_C = Z_C \cup \bigcup_{C' \in \mathcal{C}_C} S'_{C'}.
        \]
        Then, by Lemma~\ref{lemma:combine_tree_decompositions_of_S},
        \begin{align*}
            \td_k(C, S_{C}) 
            &\leq |Z_C| + \max_{C' \in \mathcal{C}_C} \td_k(C',S'_{C'}) \\
            &\leq c \cdot t
                + c \cdot f_{\cstref{thm:main_induction}}(k, \ell, t-1, a-1) + \\
            & \hspace{1cm}
                + f_{\cstref{claim:thm:main_induction}}\big(t, \ell'+1, \\
            & \hspace{1cm} \hspace{1.5cm} f_{\cstref{lemma:connecting_a_small_set}}(k,\ell,c)
                + c \cdot t 
                + c \cdot f_{\cstref{thm:main_induction}}(k, \ell, t-1, a-1),c'\big).
        \end{align*}
    
        Finally, let 
        \[
            S' = Z \cup \bigcup_{\substack{\text{$C$ connected component of}\\\text{$\torso_{G'}(U', \mathcal{B}') - Z$}}} S_C.
        \]
        Recall that $\torso_{G'}(U', \mathcal{B}') - Z = \torso_G(U,\mathcal{B})-X-Z = \torso_G(U,\mathcal{B})[V_{\ell'}]-Z$.
        Then, by Lemma~\ref{lemma:combine_tree_decompositions_of_S},
        \begin{align*}
            \td_k(\torso_G(U, \mathcal{B})[V_{\ell'}], S') 
            &\leq \td_k(\torso_G(U, \mathcal{B})[V_{\ell'}], Z)
                + \max_{\substack{C \text{ connected component of }\\ \torso_{G'}(U', \mathcal{B}')-Z}} \td_k(C,S_C)\\
            &\leq f_{\cstref{lemma:connecting_a_small_set}}(k,\ell,c) \\
            & \hspace{1cm}
                + c \cdot t\\
            & \hspace{1cm}
                + c \cdot f_{\cstref{thm:main_induction}}(k, \ell, t-1, a-1) +  \\
            & \hspace{1cm}
                + f_{\cstref{claim:thm:main_induction}}\big(t, \ell'+1, \\
            & \hspace{1cm} \hspace{1.5cm} f_{\cstref{lemma:connecting_a_small_set}}(k,\ell,c)
                    + c \cdot t 
                    + c \cdot f_{\cstref{thm:main_induction}}(k, \ell, t-1, a-1), c'\big) \\
            &\leq f_{\cstref{claim:thm:main_induction}}(t, \ell', c, c').
        \end{align*}
        Since $S \cap V_{\ell'} \subseteq S'$, this proves the claim.
    \end{proofclaim}

    \bigskip

    We now define $f_{\cstref{thm:main_induction}}(k, \ell, t, a)$ by
    \[
        f_{\cstref{thm:main_induction}}(k, \ell, t, a) 
        = \max\Big\{t,\ 2k-1 + f_{\cstref{claim:thm:main_induction}}\left(\max\left\{t,\tfrac{3}{2}(k-1)\right\}, 1, k, \max\left\{t, \tfrac{3}{2}(k-1)\right\}\right) \Big\}.
    \]

    \medskip

    Let $G$ be a graph that does not contain any $k$-ladder of length $\ell$ as a minor,
    let $(U, \mathcal{B})$ be a nice pair in $G$,
    and let $S \subseteq U$ be such that there is a tree decomposition of $(\torso_G(U, \mathcal{B}),S)$ of adhesion at most $a$ and width less than $t$.
    To deduce the theorem, we will decompose $\torso_G(U, \mathcal{B})$ into nice pairs using Lemma~\ref{lemma:existence_good_tree_decomposition}.
    By Lemma~\ref{lemma:existence_good_tree_decomposition} applied to $\torso_G(U, \mathcal{B})$ and $S$, 
    there exists $S'_0 \subseteq U$ containing $S$ and $\mathcal{D} = \big(T,(W_x \mid x \in V(T))\big)$ a tree decomposition of $(\torso_G(U, \mathcal{B}),S'_0)$ such that
    if $\mathcal{C}$ is the family of all the connected components of $\torso_G(U, \mathcal{B})-S'_0$, then
    \begin{enumerate}[label=\ref{lemma:existence_good_tree_decomposition}.(\alph*)]
        \item $\mathcal{D}$ has adhesion at most $a$; \label{item:lemma:existence_good_tree_decomposition:i:main_proof}
        \item $\mathcal{D}$ has width less than $\max\big\{t, \frac{3}{2}(k-1)\big\}$; \label{item:lemma:existence_good_tree_decomposition:ii:main_proof}
        \item for every $x_1,x_2 \in V(T)$, \label{item:lemma:existence_good_tree_decomposition:iii:main_proof}
            for every $Z_1 \subseteq W_{x_1}$ and $Z_2 \subseteq W_{x_2}$ both of size $k$,
            either
            \begin{enumerate}[label=(\roman*)]
                \item there are $k$ pairwise disjoint $(Z_1,Z_2)$-paths in $\torso_G(U, \mathcal{B})$,
                \item there exists $z_1,z_2 \in E(T[x_1,x_2])$ with $|W_{z_1} \cap W_{z_2}| < k$, or
                \item $x_1=x_2$, $|W_{x_1}| \leq \frac{3}{2}(k - 1)$ and $|W_{x_1} \cap W_y|<k$ for every $y \in N_T(x_1)$; and
            \end{enumerate}
        \item for every $x_1 x_2 \in E(T)$ with $|W_{x_1} \cap W_{x_2}|<k$,  \label{item:lemma:existence_good_tree_decomposition:iv:main_proof}
            if there exists $x_3 \in N_T(x_2)$ such that $|W_{x_2} \cap W_{x_3}| \geq k$, then
            for every positive integer $i$, for every $Z_1,Z_2 \subseteq W_{x_1} \cap W_{x_2}$ both of size $i$,
            there are $i$ pairwise disjoint $(Z_1,Z_2)$-paths in 
            \[
                \torso_G(U, \mathcal{B})\left[\bigcup_{z \in V(T_{x_2 \mid x_1})} W_z \cup \bigcup_{C \in \mathcal{C}(x_2 \mid x_1)} V(C) \right] 
                    \setminus \binom{W_{x_1} \cap W_{x_2}}{2},
            \]
            where $\mathcal{C}(x_2 \mid x_1)$ is the family of all the components $C \in \mathcal{C}$ such that $N_{\torso_G(U, \mathcal{B})}(V(C)) \subseteq \bigcup_{z \in V(T_{x_2 \mid x_1})} W_z$ and $N_{\torso_G(U, \mathcal{B})}(V(C)) \not\subseteq W_{x_1} \cap W_{x_2}$.
    \end{enumerate}
    For every connected component $C \in \mathcal{C}$, fix $x_C \in V(T)$ such that $N_{\torso_G(U, \mathcal{B})}(V(C)) \subseteq W_{x_C}$.

    Let $E_{< k}$ be the set of all the edges $zz'$ of $T$ such that $|W_z \cap W_{z'}| < k$.
    If $E_{<k} = E(T)$, then $\mathcal{D}$ has adhesion at most $k-1$ and so
    it is a $k$-dismantable tree decomposition of $(\torso_G(U, \mathcal{B}),S'_0)$ of width less than $\max\big\{t, \frac{3}{2}(k-1)\big\}$.
    This implies that $\td_k(\torso_G(U, \mathcal{B}),S'_0) \leq \max\big\{t, \frac{3}{2}(k-1)\big\} \leq f_{\cstref{thm:main_induction}}(k, \ell, t, a)$, and so the result follows for $S'=S'_0$.
    Now suppose $E_{<k} \neq E(T)$.
    
    We root $T$ on a vertex $r$ which is incident to an edge in $E(T) \setminus E_{<k}$.
    Let $T_1, \dots, T_m$ be all the connected components of $T \setminus E_{< k}$ that contain at least one edge.

    Consider some $i \in [m]$.
    Let $A_i$ be the adhesion between the root of $T_i$ and its parent if $r \not\in V(T_i)$, or the empty set if $r \in V(T_i)$.
    Let $T'_i$ be the connected component containing $T_i$ in $T - \bigcup_{j \in [m], j \neq i} V(T_j)$,
    see Figure~\ref{fig:main_proof_finding_nice_sets}.
    Then, let $T''_i$ be the subtree of $T'_i$ rooted on the root of $T_i$.
    Let 
    \[
        U_i = \bigcup_{z \in V(T'_i)} \left(W_z \cup \bigcup_{C \in \mathcal{C}, N_{\torso_G(U, \mathcal{B})}(V(C)) \subseteq W_z} V(C)\right),
    \]
    and let $\mathcal{B}_i$ be the family of all sets
    \[
        \bigcup_{z \in V(F)} \left(W_z \cup \bigcup_{\substack{C \in \mathcal{C}, N_{\torso_G(U, \mathcal{B})}(V(C)) \subseteq W_z,\\ N_{\torso_G(U, \mathcal{B})}(V(C)) \not\subseteq U_i}} V(C) \right)
    \]
    for $F$ connected component of $T-T'_i$.

    Let 
    \[
        V_i = 
        \bigcup_{z \in V(T''_i)} 
            \left(W_z \cup \bigcup_{C \in \mathcal{C}, x_C = z} V(C)\right)
    \]
    and let $\mathcal{B}''_i$ be the family of all the sets
    \[
        \bigcup_{z \in V(F)} 
            \left(
                W_z \cup 
                \bigcup_{\substack{C \in \mathcal{C}, N_{\torso_G(U, \mathcal{B})}
                           (V(C)) \subseteq W_z,\\ x_C \not\in V(T''_i)}} 
                    V(C)
            \right)
    \]
    for $F$ connected component of $T-T''_i$.
    Then let $G_i = \torso_{\torso_G(U,\mathcal{B})}\left(V_i, \mathcal{B}''_i \right)$.
    The motivation for the introduction of these graphs $G_1, \dots, G_m$ is that,
    at the end of the proof,
    we will use the fact that $\torso_G(U, \mathcal{B})$ 
    can be obtained by $(<k)$-clique-sums from $G_1, \dots, G_m$.
    Hence it will be enough to bound the $k$-treedepth
    of $(G_i,S'_i)$ for a suitable $S'_i$, for each $i \in [m]$.

    \begin{figure}[p]
        \centering
        \includegraphics{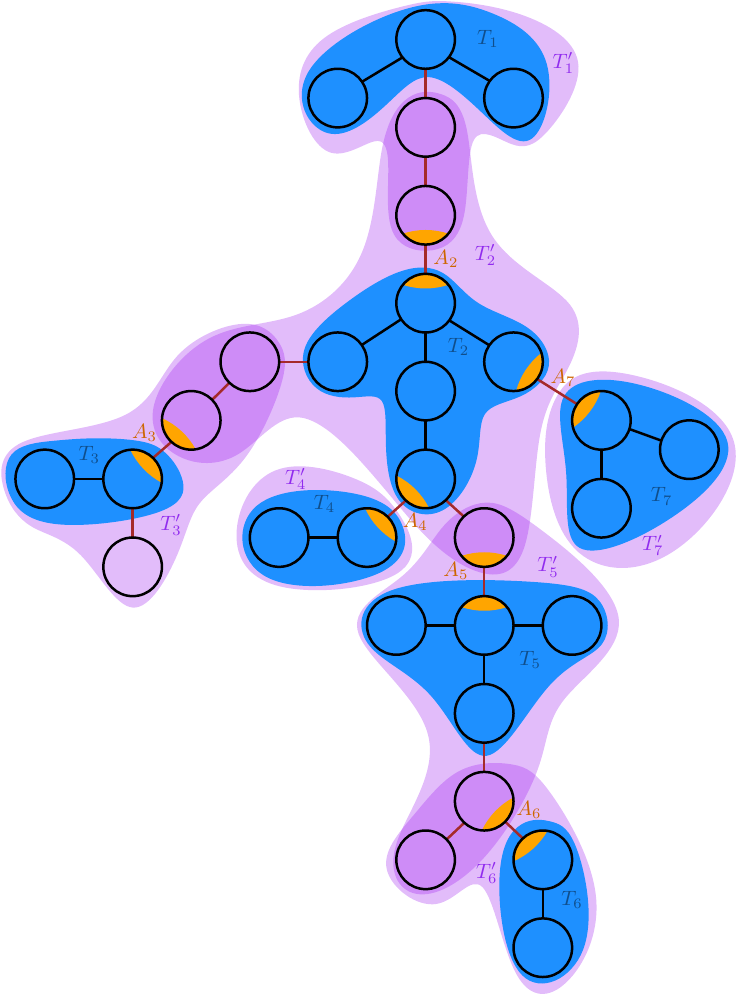}
        \caption{Application of Lemma~\ref{lemma:existence_good_tree_decomposition} in the proof of Theorem~\ref{thm:main_induction}.
            This yields a tree decomposition of $(\torso_G(U, \mathcal{B}),S'_0)$ for some $S'_0 \subseteq U$ containing $S$.
            For sake of clarity, the connected components of $\torso_G(U, \mathcal{B})-S'_0$ are not depicted.
            We want to decompose $\torso_G(U,\mathcal{B})$ using this tree decomposition as a clique sum of some ``nice'' $G_1, \dots, G_m$.
            To do so, we rely on the ``well-connected'' adhesions given by \ref{item:lemma:existence_good_tree_decomposition:iv:main_proof}.
            In orange are depicted the clique-sums we will perform.
            In black are depicted the edge $zz' \in E(T)$ with $|W_z \cap W_{z'}| \geq k$, and in brown those for which $|W_z \cap W_{z'}| < k$.
            By the properties of the tree decomposition, and in particular \ref{item:lemma:existence_good_tree_decomposition:iv:main_proof},
            each $T'_i$ corresponds
            to a nice pairs $(U_i, \mathcal{B}'_i)$ in $G$,
            and so we can apply Claim~\ref{claim:thm:main_induction} to decompose it.
            }
        \label{fig:main_proof_finding_nice_sets}
    \end{figure}

    Let $i \in [m]$.
    By construction, $(U_i, \mathcal{B}_i)$ is a good pair in $\torso_G(U, \mathcal{B})$.
    Moreover, for every connected component $F$ of $T-T'_i$, if the edge between $T'_i$ and $F$ is $x_1x_2$ with $x_1 \in V(T'_i)$,
    then $|W_{x_1} \cap W_{x_2}| < k$ and $x_2$ belongs to some $T_j$ for $j \in [m] \setminus \{i\}$, 
    which implies that $x_2$ has a neighbor $x_3$ with $|W_{x_2} \cap W_{x_3}| \geq k$.
    By \ref{item:lemma:existence_good_tree_decomposition:iv:main_proof},
    this implies that
    for every positive integer $j$,
    for every $Z_1,Z_2 \subseteq W_{x_1} \cap W_{x_2}$ both of size $j$,
    there are $j$ pairwise disjoint $(Z_1,Z_2)$-paths in 
    $\torso_G(U, \mathcal{B})\left[
        \bigcup_{z \in V(T_{x_2 \mid x_1})} 
            W_z \cup \bigcup_{C \in \mathcal{C}(x_2 \mid x_1)} V(C) 
        \right] \setminus \binom{W_{x_1} \cap W_{x_2}}{2}$.
    This implies that $(U_i, \mathcal{B}_i)$ is a nice pair in $\torso_G(U, \mathcal{B})$.
    Now, since $(U,\mathcal{B})$ is a nice pair in $G$, we deduce by Lemma~\ref{lemma:nice_in_torso_of_nice_are_nice} that 
    there exists $\mathcal{B}'_i$ such that $(U_i,\mathcal{B}'_i)$ is a nice pair in $G$
    and $\torso_{G}(U_i,\mathcal{B}'_i) = \torso_{\torso_G(U,\mathcal{B})}(U_i, \mathcal{B}_i)$.
    Now, observe that $G_i \subseteq \torso_{G}(U_i, \mathcal{B}'_i) \cup \binom{A_i}{2}$.
    We will now apply Claim~\ref{claim:thm:main_induction} to find a suitable set $S'_i$
    with $\td_k(G_i,S'_i)$ bounded. 

    \medskip

    Let $V_0$ be a subset of $k$ vertices of the bag of the root of $T_i$.
    Consider a subset $W \subseteq U_i$,
    and suppose there are no $k$ pairwise disjoint $(V_0,W)$-paths in $\torso_G(U_i, \mathcal{B}'_i)$.
    Then,
    by Lemma~\ref{lemma:finding_k_paths_between_sets_connected_in_the_torso_of_a_nice_set},
    there are no $k$ pairwise disjoint $(V_0,W)$-paths in $\torso_G(U, \mathcal{B})$, and so by \ref{item:lemma:existence_good_tree_decomposition:iii:main_proof},
    we have $|W \cap W_z| < k$ for every $z \in V(T_i)$.
    This implies that $|W_z \cap W_{z'} \cap W| <k$ for every $zz' \in E(T'_i)$, and
    so $\big(T'_i ,(W_z \cap W \mid z \in V(T'_i))\big)$ is a tree decomposition of $\big(\torso_G(U_i, \mathcal{B}'_i)[W],S'_0 \cap W\big)$
    of adhesion at most $k-1$ and width less than $\max\big\{t, \frac{3}{2}(k-1)\big\}$.
    This proves that $\td_k\big(\torso_G(U_i, \mathcal{B}'_i)[W], S'_0 \cap W\big) \leq \max\big\{t, \frac{3}{2}(k-1)\big\}$.
    Hence, we can apply Claim~\ref{claim:thm:main_induction} to
    the graph $G$,
    the nice pair $(U_i, \mathcal{B}'_i)$,
    the path partition $(V_0,U_i \setminus V_0)$, $\ell'=1$, $c=k$, 
    $c' = \max\big\{t,\frac{3}{2}(k-1)\big\}$:
    there exists $S'_{0,i} \subseteq U_i \setminus V_0$ with $S'_0 \cap (U_i\setminus V_0) \subseteq S'_{0,i}$ such that
    \begin{align*}
        \td_k\big(\torso_{G}(U_i, \mathcal{B}'_i)[U_i \setminus V_0], S'_{0,i}\big) 
          &\leq f_{\cstref{claim:thm:main_induction}}\left(\max\left\{t,\tfrac{3}{2}(k-1)\right\}, 1, k, \max\left\{t, \tfrac{3}{2}(k-1)\right\}\right).
    \end{align*}
    Let $S'_i = (S'_{0,i} \cup V_0) \cap V(G_i)$.
    Then, $S \cap U_i \subseteq S'_0 \cap U_i \subseteq S'_i$ and
    \begin{align*}
        \td_k(G_i, S'_i)
        &\leq \td_k\big(\torso_G(U_i, \mathcal{B}'_i), S'_{0,i} \cup V_0\big) + |A_i| \\
        &\leq \td_k\big(\torso_G(U_i, \mathcal{B}'_i)[U_i \setminus V_0], S'_{0,i}\big) + |V_0| + |A_i| \\
        &\leq f_{\cstref{claim:thm:main_induction}}\left(\max\left\{t,\tfrac{3}{2}(k-1)\right\}, 1, k, \max\left\{t, \tfrac{3}{2}(k-1)\right\}\right) + k + (k-1)\\
        &\leq f_{\cstref{thm:main_induction}}(k, \ell, t, a).
    \end{align*}
    Let $\mathcal{D}_i = \big(F_i, (W_z \mid z \in V(F_i))\big)$ be a $k$-dismantable tree decomposition 
    of $(G_i, S'_i)$ of width less than $f_{\cstref{thm:main_induction}}(k, \ell, t, a)$.
    
    Assume that the trees $F_i$ for $i \in [m]$ have pairwise disjoint vertex sets.
    Let 
    \[
        S' = \textstyle\bigcup_{i \in [m]} S'_i.
    \]
    Let $I$ be the set of all the pairs $(i,j)$ of distinct $i,j \in [m]$ 
    such that there is an edge between $V(T'_j)$ and $V(T_i)$ in $T$,
    with the root of $T_i$ being closer to the root $r$ in $T$ than the root of $T_j$.
    Note that $([m], \{ij \mid (i,j) \in I\})$ is a tree.
    For every $(i,j) \in I$, the set $A_i = V(G_i) \cap V(G_j)$ induces a clique in
    both $G_i$ and $G_j$, and $A_i \subseteq S'$.
    This implies that there exist $z_{(i,j)} \in V(F_i)$ and $z'_{(i,j)} \in V(F_j)$
    such that $V(G_i) \cap V(G_j) \subseteq W_{z_{(i,j)}}, W_{z'_{(i,j)}}$.
    Now, let 
    \[
        F' = \left(\textstyle\bigcup_{i \in [m]} F_i\right) \cup \{z_{(i,j)} z'_{(i,j)} \mid (i,j) \in I\}
        \quad \text{and} \quad
        S' = \textstyle\bigcup_{i \in [m]} S'_i.
    \]
    Then, $S \subseteq S'_0 \subseteq S'$
    and $\big(F', (W_z \mid z \in V(F'))\big)$ is a tree decomposition of $(G,S')$
    of width less than 
    $f_{\cstref{thm:main_induction}}(k, \ell, t, a)$.
    Recall that $|A_i| \leq k-1$ for every $i \in [m]$.
    Therefore, since $\mathcal{D}_i$ is a $k$-dismantable tree decomposition for each $i \in [m]$,
    applying \ref{item:kdismantble_def:k-cliquesum} from the definition
    of $k$-dismantable decomposition,
    we deduce that
    $\big(F', (W_z \mid z \in V(F'))\big)$ is a $k$-dismantable 
    tree decomposition of $(G,S')$.
    This proves the theorem.
\end{proof}

\section{\texorpdfstring{$k\times \ell$}{k x l} grid minor of long \texorpdfstring{$(2k-1)$}{(2k-1)}-ladders}\label{sec:token_sliding}

In this section, we show how to find a $k \times \ell$ grid in $T \square P_L$ for any tree $T$ on $2k-1$ vertices, and some large enough integer $L$.
Together with Theorem~\ref{thm:main}, this implies Corollary~\ref{corollary:excluding_only_the_kl_grid}.

\begin{lemma}\label{lemma:kl_grid_in_2k-1_ladder}
    Let $k, \ell$ be positive integers.
    For every tree $T$ on $2k-1$ vertices,
    there exists an integer $L$ such that the $k \times \ell$ grid is a minor of $T \square P_L$.
\end{lemma}

\begin{figure}[H]
    \centering
    \vspace{5mm}
    \includegraphics{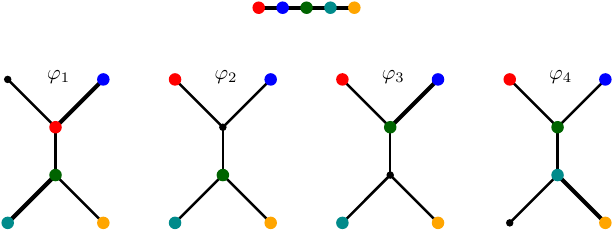}
    \caption{A sequence $(\phi_1, \phi_2,\phi_3, \phi_4)$ witnessing that $P_5 \sqsubseteq T$, where $T$ is the tree on $6$ vertices having two internal vertices of degree $3$.
        Informally, we place the vertices of $P_5$ as tokens on the graph $T$, and we move them one by one in such a way that for every edge $xy$ of $P_5$,
        the tokens of $x$ and $y$ are adjacent in $T$ at least once.
        The fat edges represent these realizations $\phi_i(x)\phi_i(y) \in E(T)$ of the edges $xy \in E(P_5)$.}
    \label{fig:sliding_token}
\end{figure}

First, let us introduce some notation.
Given two graphs $H,G$, we write $H \sqsubseteq G$ if there is a sequence $\phi_1, \dots, \phi_m$ of injective functions from $V(H)$ to $V(G)$ such that
\begin{enumerate}
    \item for every $i \in [m-1]$, $\phi_i$ and $\phi_{i+1}$ differ on exactly one vertex $x \in V(H)$, and $\phi_i(x)\phi_{i+1}(x)$ is an edge in $G$; and
    \item for every $xy \in E(H)$, there exists $i \in [m]$ such that $\phi_i(x) \phi_i(y) \in E(G)$.
\end{enumerate}
Less formally, we can see the vertices of $H$ as tokens placed on the vertices of $G$, and we successively move these tokens along the edges of $G$
in such a way that every edge of $H$ is such that the corresponding tokens are adjacent in $G$ at least once in this sequence.
See Figure~\ref{fig:sliding_token}.

This definition is motivated by the following observation.
\begin{observation}\label{observation:sliding_token_to_model}
    Let $G,H$ be two graphs.
    If $H \sqsubseteq G$, then for every positive integer $\ell$, there exists a positive integer $L$ such that
    $H \square P_\ell$ is a minor of $G \square P_L$.
\end{observation}

\begin{figure}[H]
    \centering
    \includegraphics{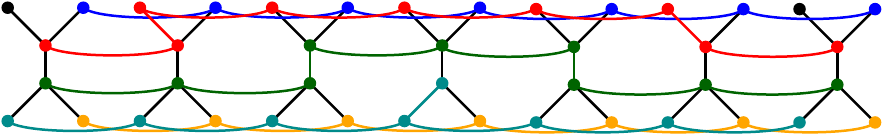}
    \caption{Illustration for Observation~\ref{observation:sliding_token_to_model}. Using the sequence $(\phi_1, \dots, \phi_4)$ given in Figure~\ref{fig:sliding_token},
        we construct a model of $P_5 \square P_\ell$ in $T \square P_{7\cdot\ell}$ by gluing together $\ell$ copies of the model $(B^1_x \mid x \in V(P_5))$ depicted here.
        For sake of clarity, some edges of $T \square P_7$ which are not used by this model are not depicted.}
    \label{fig:sliding_token_model}
\end{figure}

\begin{proof}
    Since $H \sqsubseteq G$, 
    there is a sequence $\phi_1, \dots, \phi_m$ of injective functions from $V(H)$ to $V(G)$ such that
    \begin{enumerate}
        \item for every $i \in [m-1]$, $\phi_i$ and $\phi_{i+1}$ differ on exactly one vertex $x \in V(H)$, and $\phi_i(x)\phi_{i+1}(x)$ are adjacent in $G$; and
        \item for every $xy \in E(H)$, there exists $i \in [m]$ such that $\phi_i(x) \phi_i(y) \in E(G)$.
    \end{enumerate}
    Let $L = \ell \cdot (2m-1)$.
    We label the vertices of $G \square P_L$ by $(x,i)$ for $x \in V(G)$ and $i \in [L]$. 
    For every $x \in V(H)$, let $B^0_x = \{(\phi_i(x),i) \mid i \in [m]\} \cup \{(\phi_{i}(x),i+1) \mid i \in [m-1], \phi_i(x) \neq \phi_{i+1}(x)\}$.
    Then let $B^1_x = B^0_x \cup \{(u,2m-i) \mid (u,i) \in B^0_x\}$,
    and finally $B_{(x,j)} = \{(u,i+(2m-1)(j-1)) \mid (u,i) \in B^1_x\}$ for every $j \in [\ell]$.
    See Figure~\ref{fig:sliding_token_model}.
    It is not hard to check that $(G \square P_L)[B_{(x,j)}]$ is connected,
    that the sets $(B_{(x,j)})_{x \in V(H),j \in [\ell]}$ are pairwise disjoint, and
    for every $xy \in E(H)$, there is an edge between $B_{(x,j)}$ and $B_{(y,j)}$ in $G \square P_L$, for every $j \in [\ell]$.
    Moreover, there is an edge between $B_{(x,j)}$ and $B_{(x,j+1)}$ in $G \square P_L$ for every $j \in [\ell-1]$.
    Hence $(B_{(x,j)} \mid (x,j) \in V(H \square P_\ell))$ is a model of $H \square P_\ell$ in $G \square P_L$.
    This proves the observation.
\end{proof}

Using Observation~\ref{observation:sliding_token_to_model}, it is enough to show the following to prove Lemma~\ref{lemma:kl_grid_in_2k-1_ladder}.

\begin{lemma}\label{lemma:token_sliding_2k-1}
    For every positive integer $k$, for every tree $T$ on $2k-1$ vertices,
    $P_k \sqsubseteq T$.
\end{lemma}

\begin{figure}[H]
    \centering
    \includegraphics{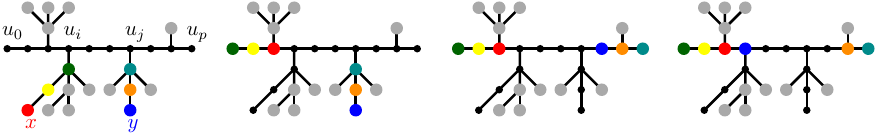}
    \caption{The main steps in the sequence given in the proof of Lemma~\ref{lemma:token_sliding_2k-1} to make adjacent the red token $x$ and the blue token $y$.
    In gray are depicted the tokens which are not moved.}
    \label{fig:token_sliding_proof}
\end{figure}

\begin{proof}
    Let $P = (u_0, \dots, u_p)$ be a longest path in $T$.
    If $p \geq k-1$, then set $m=1$ and $\phi_1(i)=u_{i-1}$, for every vertex $i$ of the path $P_k = (1, \dots, k)$,
    and it follows that $P_k \sqsubseteq T$.
    Now we assume that $p \leq k-2$, and so $|V(T) \setminus V(P)| \geq k$.
    Let $\phi_0$ be any injection from $[k] = V(P_k)$ to $V(T) \setminus V(P)$.
    We claim that it is enough to show that for every distinct $x,y \in [k] = V(P_k)$,
    there is a sequence $(\phi^{x,y}_1, \dots, \phi^{x,y}_{m_{x,y}})$ of injections from $[k]$ to $V(T)$ such that
    \begin{enumerate}
       \item for every $i \in [m_{x,y}-1]$, $\phi^{x,y}_i$ and $\phi^{x,y}_{i+1}$ differ on exactly one vertex $x \in V(H)$, and $\phi^{x,y}_i(x), \phi^{x,y}_{i+1}(x)$ are adjacent in $G$; and \label{item:token_sliding_2k-1:proof:i}
       \item $\phi^{x,y}_{m_{x,y}}(x)\phi^{x,y}_{m_{x,y}}(y) \in E(T)$;
       \item $\phi^{x,y}_1 = \phi_0$.
    \end{enumerate}
    Indeed, the sequence $\phi^{1,2}_1, \dots, \phi^{1,2}_{m_{1,2}-1}, \phi^{1,2}_{m_{1,2}}, \phi^{1,2}_{m_{1,2}-1}, \dots, \phi^{1,2}_2, \phi^{2,3}_1, \dots, \dots, \phi^{(k-1),k}_{m_{(k-1),k}}$
    then witnesses that $P_k \sqsubseteq T$.
    Note that this argument also yields $K_k \sqsubseteq T$.

    For sake of clarity, we fix $x, y \in [k] = V(P_k)$ distinct, and we will write $\phi_i$ for $\phi^{x,y}_i$.
    First set $\phi_1 = \phi_0$.
    Let $u_i$ (resp. $u_j$) be the vertex of $P$ closest to $\phi_0(x)$ (resp. $\phi_0(y)$).
    Without loss of generality, assume that $i \leq j$.
    By maximality of $P$, $|V(P[u_0, u_i[)| \geq |V(T[\phi_0(x),u_i[)|$ and $|V(P]u_j, u_p])| \geq |V(T[\phi_0(y),u_j[)|$.

    Then we move the vertices of $\phi_0([k]) \cap V(T[\phi_0(x),u_i[)$ into $\{u_0, \dots, u_{i-1}\}$ by a sequence $\phi_1, \dots, \phi_{m_0}$
    satisfying~\ref{item:token_sliding_2k-1:proof:i} and such that $\phi_{m_0}(x) = u_{i-1}$, and the vertices not in $T[\phi_0(x),u_i[$ are not moved.
    Then, we move the vertices of $\phi_0([k]) \cap V(T[\phi_0(y),u_j[)$ into $\{u_{j+1}, \dots, u_{p}\}$ by a sequence $\phi_{m_0+1}, \dots, \phi_{m_1}$
    satisfying~\ref{item:token_sliding_2k-1:proof:i} and such that $\phi_{m_1}(y) = u_{j+1}$, and the vertices not in $T[\phi_0(y),u_j[$ are not moved.
    Finally, we move $y$ from $u_{j+1}$ to $u_{i}$ by a sequence $\phi_{m_1+1}, \dots, \phi_m$ satisfying~\ref{item:token_sliding_2k-1:proof:i}
    and such that $\phi_i(z)=\phi_{m_1}(z)$ for every $z \in [k] \setminus \{y\}$.
    This sequence is depicted on Figure~\ref{fig:token_sliding_proof}.
    Then $\phi_m(x)\phi_m(y) = u_{i-1}u_i$ is an edge of $T$, and so $\phi_1, \dots, \phi_m$ is as claimed.
    This proves the lemma.
\end{proof}

\section{Conclusion}\label{sec:conclusion}

To conclude, we describe a large family of parameters including treewidth, $k$-treedepth, 
and other classical minor-monotone parameters bigger than treewidth.
The $k$-treedepth describes how a graph can be built by alternatively adding a vertex and doing some $(<k)$-clique-sums.
Another similar example is the vertex cover number, 
which counts how many vertices must be added to an edgeless graph
to obtain a given graph.
In these two examples, the value of the obtained parameter is given by the number of times a vertex is added. %
We propose to define a parameter for every such sequence of vertex additions and $(<k)$-clique-sums.

Let $\Sigma = \{a\} \cup \{s_k \mid k \geq 1\}$ be an alphabet. We denote by $\epsilon$ the empty word in $\Sigma^*$.
Intuitively, the letter $a$ corresponds to a vertex addition and the letter $s_k$ to $(<k)$-clique-sums.
We will identify regular expressions over $\Sigma$ with their corresponding languages.
For every $w \in \Sigma^*$, we define $\p_w$ as the graph parameter taking values in $\mathbb{N}_{\geq 0} \cup \{+\infty\}$
defined inductively as follows.
For every graph $G$, for every $w \in \Sigma^*$, for every positive integer $k$,
\begin{align*}
    \p_\epsilon(G) &= 
    \begin{cases}
        0 & \textrm{if $G = \emptyset$,} \\
        +\infty & \textrm{otherwise,}
    \end{cases} \\
    \p_{w a}(G) &= \min \big(\{\p_w(G)\} \cup \{\p_{w}(G-u)+1 \mid u \in V(G)\}\big), \\
    \p_{w s_k}(G) &= \min \big(\{ \p_w(G) \} \cup \{\max\{\p_{w s_k}(G_1), \p_{w s_k}(G_2)\} \mid (G_1,G_2)\}\big)
\end{align*}
where in the last case, $(G_1,G_2)$ ranges over all the pairs of graphs such that $G$ is a $(<k)$-clique-sum of $G_1$ and $G_2$ with $|V(G_1)|,|V(G_2)| > |V(G_1) \cap V(G_2)|$.
Less formally, we read from left to right $w$ as follows.
When reading the letter $a$, we can add a vertex (with arbitrary neighborhood) to a previously obtained graph, 
and when reading the letter $s_k$, we can make arbitrarily many $(<k)$-clique-sums of previously obtained graphs.
Then $\p_w(G)$ is the minimum number of times we add a vertex (when reading the letter $a$) in such a construction of $G$.
If it is not possible to construct $G$ this way, then $\p_w(G)=+\infty$.
Now, for every $L \subseteq \Sigma^*$, let 
\[
    \p_L(G) = \min_{w \in L} \p_w(G),
\]
with the convention $\min \emptyset = +\infty$.
An induction as in the proof of Lemma~\ref{lemma:ktd_minor_monotone} shows that for every $L \subseteq \Sigma^*$, the parameter $\p_L$ is minor-monotone.
Note that $\p_L(G) \geq \tw(G)+1$ for every graph $G$, for every $L \subseteq \Sigma^*$.
We denote by $s_{+\infty}$ the language $\{s_k \mid k \geq 1\}$.
Several classical minor-monotone graph parameters are of the form $\p_L$ for $L \subseteq \Sigma^*$:
for every graph $G$ with at least one edge, for every positive integer $k$,
\begin{align*}
    \p_{a^*}(G) &= |V(G)|, \\
    \p_{a^* s_1}(G) &= \max_{\substack{\text{$C$ connected} \\\text{component of $G$}}} |V(C)|, \\
    \p_{a^* s_2}(G) &= \max_{\text{$B$ block of $G$}} |V(B)|, \\
    \p_{a^* s_k}(G) &= \tw_k(G) + 1, \\
\intertext{where $\tw_k(G)$ is the $k$-treewidth of $G$~\cite{geelen2016_a_generalization_fo_the_grid_minor_theorem}, 
which is the minimum width of a tree decomposition of $G$ of adhesion less than $k$,}
    \p_{a^* s_{+\infty}}(G) & = \tw(G)+1, \\
    \p_{(a s_k)^*}(G) &= \td_k(G), \\
    \p_{(as_1)^* s_2}(G) &= \max_{\text{$B$ block of $G$}}\td(B), \\
    \p_{a s_1 a^*}(G) &= \vc(G)+1, &&\textrm{where $\vc(G)$ is the vertex cover number of $G$}, \\
    \p_{a^2 s_2 a^*}(G) &= \fvs(G)+2, &&\textrm{where $\fvs(G)$ is the feedback vertex set number of $G$}.
\end{align*}
The last two examples follows from the following property:
for every graph $G$,
\begin{align*}
    \p_{a s_1}(G) &=
    \begin{cases}
        0 &\textrm{if $G= \emptyset$,} \\
        1 &\textrm{if $E(G) = \emptyset$ and $V(G) \neq \emptyset$,} \\
        +\infty &\textrm{if $E(G) \neq \emptyset$}
    \end{cases} \\
\intertext{and}
    \p_{a^2 s_2}(G) &=
    \begin{cases}
        0 &\textrm{if $G= \emptyset$,} \\
        1 &\textrm{if $E(G) = \emptyset$ and $V(G) \neq \emptyset$,} \\
        2 &\textrm{if $E(G) \neq \emptyset$ and $G$ is a forest,} \\
        +\infty &\textrm{otherwise.}
    \end{cases}
\end{align*}
Another noteworthy example is $\p_{a^2s_2 (as_1)^*}$ which is (two plus) the elimination distance to a forest,
a parameter used in the context of parameterized complexity, see~\cite{Dekker2024}.

\begin{problem}
    Given a language $L \subseteq \Sigma^*$,
    characterize classes of graphs having bounded $\p_L$ in terms of excluded minors.
\end{problem}

There is little hope that a general answer can be given to this problem.
Indeed, for every nonnegative integer $k$,
for every graph $G$, 
\[
\p_{a^{k+1}s_{+\infty}}(G) =
\begin{cases}
    \tw(G)+1 &\textrm{if $\tw(G) \leq k$} \\ 
    +\infty &\textrm{otherwise}
\end{cases}
\]
and so a class of graphs has bounded $\p_{a^{k+1} s_{+\infty}}$ if and only if it excludes all the minimal obstructions to treewidth at most $k$,
which have not been determined yet, even for $k=4$ (see~\cite{Ramachandramurthi1997}).
Nonetheless, as highlighted in~\cite{Graph_parameter_obstructions_wqo_Paul, universal_obstructions_Paul}, the problem of determining universal obstructions for some general families of graph parameters is a challenging line of research.

An important observation is that it is enough to consider only very specific languages $L \subseteq \Sigma^*$, 
namely the downward closed languages under a specific order $\preccurlyeq^*$ on $\Sigma^*$, that we now define.
Let $\preccurlyeq$ be the order on $\Sigma$ whose relations are $a \preccurlyeq a$ and $s_i \preccurlyeq s_j$ for every positive integers $i,j$ with $i \leq j$,
and let $\preccurlyeq^*$ be the order on $\Sigma^*$ defined by,
for all $x_1, \dots, x_\ell, y_1, \dots, y_m \in \Sigma$,
$x_1 \dots x_\ell \preccurlyeq^* y_1 \dots y_m$ if and only if there exist $1 \leq i_1 < \dots < i_\ell \leq m$ such that
$x_j \preccurlyeq y_{i_j}$ for every $j \in [\ell]$.
It follows from the definition that $\p_w(G) \geq \p_{w'}(G)$ for every graph $G$ and for every $w,w' \in \Sigma^*$ with $w \preccurlyeq^* w'$.
Hence, for every $L \subseteq \Sigma^*$, $\p_L = \p_{\downarrow L}$, where $\downarrow L$ is the language $\{u \in \Sigma^* \mid \exists w \in L, u \preccurlyeq^* w\}$.
Moreover, by Higman's Lemma~\cite{Higman1952}, $(\Sigma^*, \preccurlyeq^*)$ is a well-quasi-ordering.
Hence for every $L \subseteq \Sigma^*$,
there exists a finite set $F \subseteq \Sigma^*$ (namely the set of all the minimal elements of $\Sigma^* \setminus \downarrow L$) such that 
$\downarrow L = \{w \in \Sigma^* \mid \forall u \in F, u \not\preccurlyeq^* w\}$.
In particular, $\p_L$ is computable for every fixed $L \subseteq \Sigma^*$.

\subsection*{Acknowledgements.}
Thanks to Frédéric Havet for his comments on an earlier version of this work,
to Gwenaël Joret and Piotr Micek for helpful discussions,
and to Jędrzej Hodor for a suggestion of notation.
I am also very grateful to two anonymous reviewers for their careful reading and
their valuable feedback that lead to 
improvements on the presentation of the paper.

\bibliographystyle{alpha}
\bibliography{biblio.bib}

\end{document}